\documentclass{siamart0516}

\usepackage{textcomp}
\usepackage{multirow}
\usepackage{dsfont}
\usepackage{amsfonts}
\usepackage{graphicx}
\usepackage{subfigure}
\usepackage{epstopdf}
\usepackage{algorithmicx}
\usepackage{algpseudocode}
\usepackage{filecontents}
\usepackage{enumitem}
\usepackage{bm}
\usepackage{tikz-cd}
\usepackage{booktabs}
\usepackage{siunitx}

\usepackage{amsmath,amscd, amssymb}
\tikzcdset{
arrow style=tikz,
diagrams={>={Stealth[scale=0.9]}}
}
\def\D{\mathrm{d}}
\newcommand{\vertiii}[1]{{\left\vert\kern-0.25ex\left\vert\kern-0.25ex\left\vert #1 
    \right\vert\kern-0.25ex\right\vert\kern-0.25ex\right\vert}}

\ifpdf
  \DeclareGraphicsExtensions{.eps,.pdf,.png,.jpg}
\else
  \DeclareGraphicsExtensions{.eps}
\fi

\newcommand{\TheTitle}{Stochastic Primal-Dual Method on Riemannian Manifolds with Bounded Sectional Curvature} 
\newcommand{\TheAuthors}{Masoud Badiei Khuzani, Na Li}

\headers{Constrained Optimization on Riemannian Manifolds}{\TheAuthors}

\title{{\TheTitle}
}

\author{
  Masoud Badiei Khuzani, Na Li
\thanks{Harvard University, Cambridge, MA \email{(mbadieikhuzani@g.harvard.edu,nali@seas.harvard.edu)}
 }
}

\usepackage{amsopn}

\ifpdf
\hypersetup{
  pdftitle={\TheTitle},
  pdfauthor={\TheAuthors}
}
\fi

\usepackage{multirow}
\usepackage{url}
\newtheorem{problem}{Problem}
\theoremstyle{definition}
\newtheorem{assumption}{Assumption}
\newtheorem{remark}{\normalfont\scshape Remark}

\newcommand{\expect}{{\rm I\!E}}
\newcommand{\real}{{\rm I\!R}}
\newcommand{\prob}{{\rm I\!P}}

\newcommand{\ball}{{\rm I\!B}}

\newcommand{\RN}[1]{%
  \textup{\uppercase\expandafter{\romannumeral#1}}%
}

\DeclareMathOperator{\sinc}{sinc}
\DeclareMathOperator{\grad}{grad}
\DeclareMathOperator{\Exp}{exp}
\DeclareMathOperator{\topdoteq}{\dot{=}\ }
\allowdisplaybreaks[4]

\begin{document}

\maketitle

\begin{abstract}
  We study a stochastic primal-dual method for constrained optimization over Riemannian manifolds with bounded sectional curvature. We prove non-asymptotic convergence to the optimal objective value. More precisely, for the class of hyperbolic manifolds, we establish a convergence rate that is related to the sectional curvature lower bound. To prove a convergence rate in terms of sectional curvature for the elliptic manifolds, we leverage Toponogov's comparison theorem. In addition, we provide convergence analysis for the asymptotically elliptic manifolds, where the sectional curvature at each given point on manifold is locally bounded from below by the distance function. We demonstrate the performance of the primal-dual algorithm on the sphere for the \textit{non-negative} principle component analysis (PCA). In particular, under the non-negativity constraint on the principle component and for the symmetric spiked covariance model, we empirically show that the primal-dual approach outperforms the spectral method. We also examine the performance of the primal-dual method for the anchored synchronization from partial noisy measurements of relative rotations on the Lie group ${\rm SO(3)}$. Lastly, we show that the primal-dual algorithm can be applied to the weighted MAX-CUT problem under constraints on the admissible cut. Specifically, we propose different approximation algorithms for the weighted MAX-CUT problem based on optimizing a function on the manifold of direct products of the unit spheres as well as the manifold of direct products of the rotation groups.  
\end{abstract}



\section{Introduction}
Modern optimization problems emerging in statistics and signal processing are increasingly high dimensional and thus scalable numerical techniques are needed to analyze these problems. Optimization methods on Riemannian manifolds are promising techniques for constrained optimization problems that exploit the geometry of the feasible set to improve the computational efficiency. To use these techniques, the feasible set must admit the structure of a Riemannian manifold,
\begin{subequations}
\begin{align}
\label{Eq:Informal-1}
&\min_{x\in \mathcal{M}} f(x)\topdoteq\expect_{P}[F(x;\xi)]=\int_{\Xi}F(x;\xi)dP(\xi)      \\ \label{Eq:Informal-2} 
&\text{subject to}:h_{k}(x)\leq 0, \quad k=1,2,\cdots,m,
\end{align}
\end{subequations}
where $\{F(\cdot;\xi),\xi\in \Xi \}$ and $\{h_{k}\}$ are collections of real valued functions defined on the Riemannian manifold $\mathcal{M}$, and $P$ is a probability measure on $\Xi$. Here, by including the explicit inequality constraints \eqref{Eq:Informal-2}, the objective function $f$ is optimized over a sub-manifold of $\mathcal{M}$.

The problem in eqs. \eqref{Eq:Informal-1}-\eqref{Eq:Informal-2} can alternatively be framed as an Euclidean optimization using the Whitney embedding theorem \cite{whitney1936differentiable}. Specifically, the problem in eqs. \eqref{Eq:Informal-1}-\eqref{Eq:Informal-2} can be treated as an optimization problem in the Euclidean space with the implicit constraint set $\mathcal{X}$, where $\mathcal{X}$ is the image of the embedding map $\iota:\mathcal{M}\hookrightarrow {\rm I\!R}^{2n+1}$ for an $n$-dimensional manifold $\mathcal{M}$. To solve such constrained optimization problems in the Euclidean space, many algorithms have been devised over recent years, see \cite{bertsekas1976multiplier} for a survey of those results. However, the Euclidean optimization techniques are inadequate for optimizing functions on manifolds since (\textit{i}) in most cases the Euclidean optimization methods rely on the convexity of the feasible set $\mathcal{X}$ to ensure convergence, (\textit{ii}) optimization over the ambient space must be carried out on a much larger dimension than that of the low-dimensional manifold, and (\textit{iii}) the embedding map $\iota$ is often difficult to compute.

To circumvent these issues, we propose a primal-dual algorithm to directly optimize the objective function on the manifold $\mathcal{M}$, where hereafter we refer to as the Riemannian primal-dual algorithm. In the proposed method, the primal variables are restricted to the manifold $\mathcal{M}$ hence the underlying method is projection free. Consequently, the proposed algorithm scales adequately to high dimensional problems. 


\subsection{Contributions}
We prove non-asymptotic convergence of the primal-dual algorithm for optimizing functions on manifolds with a bounded sectional curvature. Specifically, for geodesically convex objective functions and for the class of manifolds with the negative sectional curvature (\textit{i.e.} the hyperbolic manifolds), we prove a global convergence in the objective value that involves the sectional curvature of the underlying manifold. 

Similarly, for the class of manifolds with the positive sectional curvature (\textit{i.e.} the elliptic manifolds), we prove a global convergence to the optimal objective value. The approach we take for this class of manifolds is based on the comparative geometry techniques. More precisely, using Toponogov's comparison theorem, we prove a convergence bound that includes the sectional curvature. To demonstrate the consistency of the upper bounds for elliptic and hyperbolic manifolds, we show that when the sectional curvature diminishes (\textit{i.e.} flat manifold), the convergence bounds coincide up to a constant multiplicative factor. We also extend our global convergence analysis to asymptotically elliptic manifolds, where the sectional curvature is lower bounded by the quadratic distance function.

We consider three applications for the Riemannian primal-dual algorithm, namely (\textit{i}) the non-negative principle component analysis (PCA), (\textit{ii}) the anchored state synchronization on the rotation group ${\rm SO}(3)$, and (\textit{iii}) the weighted MAX-CUT problem under graph constraints.

First, we consider the non-negative PCA problem, where we consider a symmetric spiked covariance model. We compare the Riemannian primal-dual algorithm to the spectral method for the traditional PCA problem. Interestingly, we show empirically that the Riemannian primal-dual method outperforms the spectral method in low SNR regimes under the cosine similarity metric, due to the non-negativity constraint on the principle component. The same observation has been made in \cite{montanari2016non} using the approximate message passing, and thus our algorithm can be viewed as an alternative approach. We further show that our proposed framework scales very gracefully with the sample size while the semi-definite programming (SDP) approach is limited to sample sizes of a few hundreds.

Next, we apply the Riemannian primal-dual algorithm to the anchored synchronization problem, where the objective is to identify the elements of the rotation group $\rm{SO}(3)$ from partial noisy observations of their relative rotations. In particular, we consider the problem of minimizing the mean squared error (MMSE) associated with the noisy measurements. Due to the invariance of the mean squared error under rotation, we include certain inequality constraints (anchors) to the MMSE problem as a symmetry breaking approach. To handle the inequality constraints, we in turn use a primal-dual method and show that it converges to the true rotation states for the measurements that are defined on the Erd\"{o}s-R\'{e}yni random graphs.

Lastly, we use the primal-dual algorithm for the weighted MAX-CUT problem under certain restrictions on admissible cuts, see \cite{lee2016max}. The weighted MAX-CUT problem belongs to the class of integer programming problems, and hence it is intractable in general, unless $P=NP$. However, by using the standard `lifting' argument \cite{goemans1995improved}, the MAX-CUT problem can be relaxed to optimizing a function on the manifold of direct product of unit $d$-spheres ${\rm S}^{d}\otimes {\rm S}^{d}\otimes \cdots \otimes {\rm S}^{d}$, where $d$ is the dimension of the lifted vector. We thus apply optimization methods on the manifold to this problem. Remarkably, regardless of initialization point, the optimization techniques on manifolds achieves solutions that matches that of SDP on dense  Erd\"{o}s-R\'{e}yni random graphs.

\subsection{Previous Works} 
Optimization methods on Riemannian manifolds generalize the Euclidean optimization techniques, and thus they are classified similarly \cite{ji2007optimization}:
\begin{itemize}[leftmargin=*]
\vspace*{2mm}
\item \textit{Gradient descent methods}. This is the class of first-order methods that were studied for manifolds by Luenberger \cite{luenberger1972gradient,luenberger1973introduction} and Gabay \cite{gabay1982minimizing}. These methods were further expanded for systems and control theory applications by Brockett \cite{brockett1993differential}, Smith \cite{smith2013geometric}, and Mahony \cite{mahony1994optimization}. A stochastic gradient descent for the Hadamard-Cartan hyperbolic manifolds is also studied \cite{bonnabel2013stochastic}. However, the convergence analysis of \cite{bonnabel2013stochastic} for the elliptic manifolds does not involve the sectional curvature.
\vspace*{2mm}
\item \textit{Newton and quasi-Newton's methods}. In Newton's method, the Hessian information of a function is used to accelerate the convergence speed. However, the complexity of computing the Hessian matrix as well as the matrix inversion in Newton's method restricts its use to low dimensional problems. Quasi-Newton's methods are first-order methods with a super-linear local convergence rate that update the Hessian through analyzing successive gradient vectors.
\vspace*{2mm}
The Riemannian Newton's method was studied by Smith \cite{smith1994optimization} and Mahony \cite{mahony1996constrained} for the compact Lie groups, and Riemannian quasi-Newton's methods on the Grassmann
manifold were studied in \cite{brace2006improved,savas2010quasi}.
\vspace*{2mm}
\item \textit{Conjugate gradient method}. Conjugate gradient method is an algorithm for solving systems of linear equations whose matrix is symmetric and positive semi-definite. This method can also be used as an iterative method for unconstrained minimizination \cite{fletcher1964function}. Conjugate gradient method was generalized to the Riemannian manifolds by Smith \cite{smith1994optimization} and was applied to problems on the Stiefel and Grassmann manifolds by Edelman, \textit{et al.} \cite{edelman1998geometry}.

\end{itemize} 
 
Note that with the exception of \cite{bonnabel2013stochastic}, these optimization methods are deterministic. Further, these methods do not include the inequality constraints in \eqref{Eq:Informal-2} in the optimization problem. However, incorporating the inequality constraints is useful as it extends the Riemannian optimization techniques to sub-manifolds of well studied manifolds such as Lie groups or a sphere. 
 
The remainder of the paper is organized as follows. In Section \ref{sec:Perl}, we review basic materials of differential geometry that are necessary for our analysis. In Section \ref{sec:problem_statement}, we formulate the optimization problem and we also present a primal-dual algorithm for solving it. We also study the convergence rate of the proposed algorithm on two different classes of manifolds, namely manifolds with negative bounded curvature, and the elliptic manifolds of positive bounded curvature. In Section \ref{Applications}, we develop Riemannian primal-dual method for the non-negative PCA problem. Further, we use the primal-dual algorithm for synchronization of rotations from noisy observations over the Lie group ${\rm SO}(3)$. We also study the application of the primal-dual method for the MAX-CUT problem on random graphs. We provide the proofs of the main theorems in Section \ref{sec:Proofs}, while the more technical parts of the proofs are deferred to the appendices. Lastly, we discuss our results and conclude the paper in Section \ref{Sec:discussion}.
 
\section{Preliminaries}
\label{sec:Prel}

To provide a self-contained presentation, we review relevant Riemannian geometry definitions here. A brief review of the Lie groups is also provided in Appendix \ref{sec:Lie Groups}. A more comprehensive treatise on this subject can be found in \cite{lee2009manifolds,do1976differential}.

\subsection{Riemannian manifold} A Riemannian manifold $(\mathcal{M}, g)$ is a real smooth manifold $\mathcal{M}$ equipped with a Riemannain metric $g$. In particular, the metric $g$ is defined at each point $p\in \mathcal{M}$ by the bi-linear map $g_{p}:T_{p}\mathcal{M}\times T_{p}\mathcal{M}\rightarrow \real$ and induces an inner product structure on the tangent space $T_{p}\mathcal{M}$. For any two vector fields $X,Y\in \mathfrak{X}(\mathcal{M})$ on the space of smooth time invariant vector fields $\mathfrak{X}(\mathcal{M})$, we denote the inner product by $\langle X_{p},Y_{p} \rangle\ \dot{=}\ g_{p}(X_{p},Y_{p})$. Accordingly, the norm of a vector $X_{p}\in T_{p}\mathcal{M}$ is denoted by $\|X_{p}\|\topdoteq(g_{p}(X_{p},X_{p}))^{1\over 2}$.
\begin{definition}
The affine connection on a manifold $\mathcal{M}$ is a bilinear map $\mathfrak{X}(\mathcal{M})\times \mathfrak{X}(\mathcal{M})\rightarrow \mathfrak{X}(\mathcal{M})$ and $(X,Y)\mapsto \nabla_{X}Y$, such that for all functions $f\in C^{2}(\mathcal{M},\real)$ it satisfies
\begin{itemize}
\item $\nabla_{fX}Y=f\nabla_{X}Y$
\item $\nabla_{X}(fY)=df(X)Y+f\nabla_{X}Y$,
\end{itemize}
where $df(X)$ is the directional derivative of function $f$ in the direction of $X$.
\end{definition}

Define a coordinate chart $(\mathcal{U},\varphi)$ as $\varphi:\mathcal{U}\rightarrow \mathcal{V}$, where $\mathcal{U}\subset \mathcal{M}$ is an open set and $\mathcal{V}\subset \real^{n}$. The affine connection on an $n$-dimensional manifold is completely determined by $n^{3}$ real valued smooth functions on $\mathcal{U}$, namely the Christoffel symbols $\Gamma_{ij}^{k}$ on local coordinates $(u_{1},\cdots,u_{n})$. Let $\partial_{k}$ denotes the vector field $\partial/\partial u_{k}$ on $\mathcal{U}$. Then, in local coordinates, the affine connection can be characterized in terms of the covariant derivative of basis vectors
\begin{align*}
\nabla_{\partial_{i}}\partial_{j}=\sum_{k=1}^{n}\Gamma_{ij}^{k}\partial_{k}.
\end{align*}

The \textit{Levi-Civita} connection on $\mathcal{M}$ is an affine connection that (i) preserves the metric $\nabla g=0$, and (ii) is torsion free, \textit{i.e.}, for any two vector fields $X,Y\in \mathfrak{X}(\mathcal{M})$ we have $\nabla_{X}Y-\nabla_{Y}X=[X,Y]$ where $[\cdot,\cdot]$ is the Lie bracket of the vector fields. Throughout the paper, we adopt this connection. In conjunction with the Levi-Civita connection, the \textit{Christoffel} symbols of the second kind take the following particular form
\begin{align*}
\Gamma_{ij}^{k}\topdoteq {1\over 2}\sum_{l=1}^{n}g^{kl}(\partial_{j}g_{il}+\partial_{i}g_{jl}-\partial_{l}g_{ij})
\end{align*}
where $g_{rj}=\langle\partial_{r}, \partial_{j}\rangle$, and $[g^{rj}]$ is the inverse of the matrix $[g_{rj}]$ such that $g^{rj}g_{rj}=\delta_{rj}$, where $\delta_{rj}$ is the Dirac delta function.

\begin{definition}
\label{Definition:gradient}
The gradient of a real-valued function $f\in C^{2}(\mathcal{M})$ at the point $p\in \mathcal{M}$, denoted
by $\grad_{p} f$, is the unique vector in $T_{p}\mathcal{M}$ such that $\langle \grad_{p} f,X_{p} \rangle=d_{p}f(X_{p})$
for all $X_{p}\in T_{p}\mathcal{M}$. 
\end{definition}

\begin{definition}\textsc{\cite[Def. 13.3, p. 358]{gallier2012notes}}
Given a real-valued function $f\in C^{2}(\mathcal{M})$ on a Riemannian manifold, the Riemannian Hessian of $f$ at a point $p\in \mathcal{M}$ is a $(0,2)$-tensor field $\text{Hess} f(p): T_{p}\mathcal{M}\times
 T_{p}\mathcal{M}\rightarrow \real$ defined as
\begin{align*}
\text{Hess}f(p)[X_{p},Y_{p}]\topdoteq\langle\nabla_{X_{p}} \grad_{p} f, Y_{p}\rangle=(\nabla^{2}f(X_{p},Y_{p}))_{p},
\end{align*}
where $\nabla^{2}f\topdoteq\nabla\nabla f$. Further, let $\nabla_{X_{p}}^{2}f\topdoteq (\nabla^{2}f(X_{p},X_{p}))_{p}$.
\end{definition}

In the coordinate chart $(\mathcal{U},\phi)$ where $p\in \mathcal{U}\subset \mathcal{M}$, the Hessian tensor takes the following form
\begin{align*}
(\nabla^{2}f)_{p}=\sum_{m,\ell=1}^{n}\left(\left({\partial_{m\ell}^{2} f}\right)_{p} -\sum_{k=1}^{n}\Gamma_{\ell m}^{k}\left({\partial_{k}^{2} f}\right)_{p}\right)du_{m}\otimes du_{\ell},
\end{align*}
where $\otimes$ denotes the tensor product on the product of the cotangent spaces $T_{p}^{*}\mathcal{M}\otimes T_{p}^{*}\mathcal{M}$. 

The Riemann curvature tensor in terms of the Levi-Civita connection is defined as follows
\begin{align*}
R(X,Y)Z\topdoteq\nabla_{X}\nabla_{Y}Z-\nabla_{Y}\nabla_{X}Z-\nabla_{[X,Y]}Z, \quad X,Y,Z\in \mathfrak{X}(\mathcal{M}).
\end{align*} 
The sectional curvature at point $p\in \mathcal{M}$ is then defined as 
\begin{align}
\label{Eq:sectional_curvature}
K_{p}(X_{p},Y_{p})\topdoteq\dfrac{\langle R(X_{p},Y_{p})X_{p},Y_{p}\rangle}{\langle X_{p},X_{p}\rangle\langle Y_{p},Y_{p}\rangle-\langle X_{p},Y_{p}\rangle^{2}},
\end{align}
where $X_{p},Y_{p}\in T_{p}\mathcal{M}$. In particular, when $X_{p}$ and $Y_{p}$ are orthonormal, then we simply have $K_{p}(X_{p},Y_{p})=\langle R(X_{p},Y_{p})X_{p},Y_{p}\rangle$.

\subsection{Geodesic curve} A geodesic on a smooth manifold $\mathcal{M}$ equipped with an affine connection $\nabla$ is defined as a curve $\gamma:[0,1]\rightarrow \mathcal{M}$ such that the parallel transport along the curve preserves the tangent vector to the curve. Specifically,
\begin{align*}
\nabla_{\dot{\gamma}(t)}\dot{\gamma}(t)=0,
\end{align*}
at each point along the curve, where $\dot{\gamma}(t)$ is the derivative with respect to $t$. Alternatively, the geodesics of the Levi-Civita connection can be defined as the locally distance-minimizing paths. Specifically, in a Riemannian manifold $\mathcal{M}$ with the metric tensor $g$, the length of a path $\gamma:[0,1]\rightarrow \mathcal{M}$ is the following functional,
\begin{align}
\label{Eq:functional}
L[\gamma]\topdoteq\int _{0}^{1}\sqrt{g_{\gamma(t)}({\dot {\gamma }}(t),\dot {\gamma }(t))}\,dt.
\end{align}
Accordingly, the distance $d(p,q)$ between two points $p,q\in \mathcal{M}$ is defined as
\begin{align*}
d(p,q)\topdoteq \inf\left\{L[\gamma]:\gamma:[0,1]\rightarrow \mathcal{M},\gamma(0)=p,\gamma(1)=q, \gamma\  \text{is piecewise smooth} \right\}.
\end{align*}
Using the Euler-Lagrange equations to minimize the functional in eq. \eqref{Eq:functional} yields the following set of differential equations for geodesics in local coordinates,
\begin{align}
\label{Eq:Computing_Geodesics}
\ddot{\gamma}_{i}(s)+\sum_{j,k=1}^{n}\Gamma_{jk}^{i}\dot{\gamma}_{j}(s)\dot{\gamma}_{k}(s)=0.
\end{align}

\begin{definition}
\label{Def:exponential_map}
The (geodesic) exponential map $\Exp_{p}:T_{p}\mathcal{M} \rightarrow \mathcal{M}$ is a mapping from $v\in T_{p}\mathcal{M}$
to $q\in \mathcal{M}$. In particular, 
\begin{align*}
\Exp_{p}(v)=\gamma_{v}(1), \quad v\in T_{p}\mathcal{M},
\end{align*}
where $\gamma:[0,1]\rightarrow \mathcal{M}$ is a geodesic such that $\gamma_{v}(0)=p$, $\gamma_{v}(1)=q$, and $\dot{\gamma}_{v}(0)=v$.
\end{definition}
In the case that the exponential map admits an inverse, we denote it by $\log_{p}(q):\widetilde{\mathcal{M}}\rightarrow T_{p}{\mathcal{M}},\widetilde{\mathcal{M}}=\exp_{p}(T_{p}\mathcal{M})\subset \mathcal{M}$, and we refer to it as the Riemannian logarithmic map or simply the logarithmic map.

\begin{definition}
\label{Def:geodesic_convexity}
 A real-valued function $f:\mathcal{M}\rightarrow \real$ is geodesically convex if for any two points $p,q\in \mathcal{M}$ we have  
\begin{align}
\label{Eq:In_Terms}
f(\gamma(t))\leq (1-t)f(p)+tf(q), \quad \text{for all}\ t\in [0,1],
\end{align} 
where $\gamma(t):[0,1]\rightarrow \mathcal{M}$ is the geodesic with the initial and end points of $\gamma(0)=p$ and $\gamma(1)=q$, respectively. Equivalently, in terms of the logarithmic map $\log_{p}(\cdot)$, the condition \eqref{Eq:In_Terms} is described by
\begin{align}
\label{Eq:geodesically_convex}
f(p)+\langle \grad_{p} f,\log_{p}(q) \rangle\leq f(q), \quad \text{for all}\ p,q\in \mathcal{M}.
\end{align}
\end{definition}

Now, we consider the geodesics of the product manifold $\mathcal{N}={\mathcal{M}}\otimes {\mathcal{M}}\otimes\cdots\otimes {\mathcal{M}}$ with the product metric tensor. In particular, we argue that $\gamma(t)=(\gamma_{1}(t),\cdots,\gamma_{n}(t))$ determines a geodesic for $\mathcal{M}$, where each $\gamma_{i}(t)$ is a geodesic for the $i$-th manifold. To see this, consider the inclusion embedding $\iota:\mathcal{M}\rightarrow \mathcal{N}$. Given a vector $X\in T_{p}\mathcal{M}$, we denote the push-forward of $X$ to $T_{\iota(p)}\mathcal{N}$ via $(d\iota)_{p}:T_{p}\mathcal{M}\rightarrow T_{\iota(p)}\mathcal{N}$ by $\widehat{X}$. Note that this extension $\widehat{X}$ is not unique. Let $\nabla^{\mathcal{N}}$ and $\nabla^{\mathcal{M}}$ denote the Levi-Civita connection for $\mathcal{N}$ and $\mathcal{M}$, respectively. By the Gauss formula \cite[Thm. 8.2]{lee1997riemannian_C}, the relation between the Levi-Civita connections on $\mathcal{N}$ and $\mathcal{M}$ are established by  
\begin{align}
\label{Eq:Second_Fundamental_Form}
\nabla^{\mathcal{N}}_{\widehat{X}}\widehat{Y}=\nabla^{\mathcal{M}}_{X}Y+\RN{2}(\widehat{X},\widehat{Y}),
\end{align}
where $\RN{2}(\widehat{X},\widehat{Y})\topdoteq \left(\nabla^{\mathcal{N}}_{\widehat{X}}\widehat{Y}\right)^{\perp}\topdoteq \Pi_{(T_{p}\mathcal{M})^{\perp}}(\nabla^{\mathcal{N}}_{\widehat{X}}\widehat{Y})$ is \textit{the second fundamental form}, and $(T_{p}\mathcal{M})^{\perp}$ denotes the orthogonal complement of $T_{p}\mathcal{M}$.

From eq. \eqref{Eq:Second_Fundamental_Form}, we observe that there exists a natural inclusion of $X$ in $\mathcal{N}$, such that $\RN{2}(\widehat{X},\widehat{Y})=0$. Specifically, consider the local coordinates $(x_{1},\cdots,x_{d})$ for $p\in \mathcal{M}$. We locally have $X=\sum_{i=1}^{d}X_{i}(p){\partial_{i}}$, where $X_{i}(p)$ are smooth functions. Now, consider the local coordinate $(x_{1},\cdots,x_{d},\hat{x}_{1},\cdots,\hat{x}_{(n-1)d})$ of $\iota(p)\in \mathcal{M}$. Then, for the canonical inclusion 
\begin{align}
\label{Eq:construction}
\widehat{X}=\sum_{i=1}^{d}X_{i}(\iota(p)){\partial_{i}}+\sum_{i=1}^{(n-1)d}0\times {\partial_{i}},
\end{align}
we have $\RN{2}(\widehat{X},\widehat{Y})=0$, see \cite{taringoo2015synchronization}.

We now notice that due to eq. \eqref{Eq:Second_Fundamental_Form} and the choice of the extension in eq. \eqref{Eq:construction},  we have $\nabla^{\mathcal{N}}_{X_{1}+\cdots+X_{n}}(Y_{1}+\cdots+Y_{n})=\nabla^{\mathcal{M}}_{X_{1}}Y_{1}+\cdots+\nabla^{\mathcal{M}}_{X_{n}}Y_{n}$, where $X_{i},Y_{i}\in \mathfrak{X}(\mathcal{M}),i=1,2,\cdots,n$; see \cite[Excer. 1, Chap. 6]{do1976differential}. The claim now follows by letting $X_{i}=Y_{i}=\dot{\gamma}_{i}(t)$ and noting that for a geodesic $\nabla^{\mathcal{M}}_{\dot{\gamma}_{i}(t)}\dot{\gamma}_{i}(t)=0$. 
 
\subsection{Parallel Transport}

Parallel transport provides a way to compare different tangent spaces on manifolds. For a curve $\gamma:[0,1]\rightarrow \mathcal{M},t\mapsto \gamma(t)$ with the local expression $(\gamma_{1}(t),\cdots,\gamma_{n}(t))$ and a vector field $X(t)=\sum_{i}x_{i}(t)(\partial_{i})_{\gamma(t)}$, we define
\begin{align*}
\dfrac{\D X}{\D t}=\sum_{i}\dfrac{\D x_{i}}{\D t}\partial_{i}+\sum_{i}x_{i}\sum_{j}\dfrac{\D \gamma_{j}}{\D t}\nabla_{\partial_{j}}\partial_{i}.
\end{align*}
The vector field $X(t)$ is called \textit{parallel} if ${\D X}/{\D t}=0$. For a parallel vector field, we use the notation $P_{\gamma(s)\rightarrow \gamma(t)}X(s)$ to denote the transport of the vector $X(s)\in T_{\gamma(s)}\mathcal{M}$ to the tangent space at $\gamma(t)$, \textit{i.e.}, $P_{\gamma(s)\rightarrow \gamma(t)}X(s)\in T_{\gamma(t)}\mathcal{M}$.

\section{Problem Statement and Algorithm}
\label{sec:problem_statement}

\subsection{Primal-Dual Algorithm}
In this paper, we are concerned with the following optimization problem:
\begin{problem}
\label{Problem}
Suppose $\mathcal{M}$ is a smooth, complete, connected
Riemannian manifold, equipped with a smooth metric tensor, and $F(\cdot,\xi),h_{1},\cdots,h_{m}\in C^{2}(\mathcal{M},\real)$ are functions that are obtained through restriction of their domains to the manifold $\mathcal{M}$. Compute
\begin{align*}
&\min_{x\in \mathcal{M}} f(x)\topdoteq\expect_{P}[F(x;\xi)]=\int_{\Xi}F(x;\xi)dP(\xi),     \\ 
&\text{subject to}:  h(x)\preceq 0, 
\end{align*}
where $h\topdoteq(h_{1},h_{2},\cdots,h_{m})$.
\end{problem}

Notice that since $\mathcal{M}$ is a finite dimensional complete metric space by the problem statement assumption, the Hopf-Rinow theorem \cite{do1976differential},\cite{milnor2016morse} ensures that it is also a geodesically complete manifold, \textit{i.e.}, any two points in $\mathcal{M}$ are joined by \textit{at least} one minimal geodesic. 

To solve Problem \ref{Problem}, we define the augmented Lagrangian functions $\mathfrak{L}(x,\lambda;\xi)\in C^{2}(\mathcal{M}\times \real_{+}^{m}\times \Xi,\real)$ on the manifold $\mathcal{M}$ as follows
\begin{align}
\label{Eq:Label_d}
\mathfrak{L}(x,\lambda;\xi)&\topdoteq F(x;\xi)+\langle\lambda, h(x)\rangle-{\alpha\over 2}\|\lambda\|^{2},\\ \nonumber
\mathfrak{L}(x,\lambda)&\topdoteq\expect_{P}[\mathfrak{L}(x,\lambda;\xi)]\\ \label{Eq:Label_d1}
&=f(x)+\langle\lambda, h(x)\rangle-{\alpha\over 2}\|\lambda\|^{2},
\end{align}
where $\lambda\in \real^{m}_{+}$ is the Lagrangian dual vector for inequality constraints, and $\alpha>0$ is the regularizer parameter that controls the norm of the dual variables. Due to augmenting the Lagrangian function $\mathfrak{L}(x_{t},\lambda_{t};\xi_{t})$ with the regularizer term ${\alpha\over 2}\|\lambda\|^{2}$ in \eqref{Eq:Label_d}, the norm of the Lagrangian dual variables $\|\lambda\|$ are upper bounded. This in turn ensures that the gradients of the deterministic Lagrangian function $\mathfrak{L}(x_{t},\lambda_{t})$ are also upper bounded (cf. Lemma \ref{Lemma:3}). 
We defer the discussion about how to choose the parameter $\alpha$ to the next section. We only emphasize here that the regularization term is mainly included for our proof techniques to work. In the numerical experiments, we did not find any performance loss in setting $\alpha=0$.

A solution to Problem \ref{Problem} is the saddle point $(x_{\ast},\lambda_{\ast})\in \mathcal{M}\times \real_{+}$ for the Lagrangian function $\mathfrak{L}(x,\lambda)$. In particular, the saddle (min-max) point $(x_{\ast},\lambda_{\ast})$ satisfies
\begin{align}
\label{Eq:saddle_point}
\mathfrak{L}(x_{\ast},\lambda)\leq \mathfrak{L}(x_{\ast},\lambda_{\ast})\leq \mathfrak{L}(x,\lambda_{\ast}), \quad \forall x\in \mathcal{M},\forall \lambda\in \real_{+}.
\end{align}
Based on eq. \eqref{Eq:Label_d}, we define
\begin{align}
\label{Eq:G_R0}
\grad_{x}\mathfrak{L}(x,\lambda;\xi)\ &\dot{=}\ \grad_{x}F(x;\xi)+\langle\lambda, \grad_{x} h(x)\rangle\\
\label{Eq:G_R1}
\grad_{\lambda}\mathfrak{L}(x,\lambda)\ &\dot{=}\ h(x)-\alpha \lambda.
\end{align}

Algorithm \ref{alg:1} now describes the primal-dual method on the Riemannian manifold $\mathcal{M}$. At Step 3 of the algorithm, a point from the manifold $x_{t}\in \mathcal{M}$ is queried from the stochastic oracle, and an estimate $\grad_{x_{t}}F(x_{t},\xi_{t})$ of the true gradient $\grad_{x_{t}}f(x_{t})$ is received.

\subsection{Assumptions}
We impose the following assumptions throughout the paper:
\begin{assumption}
\label{Assumption:1}
There exists a sub-manifold $\mathcal{U}(x_{\ast})\subset \mathcal{M}$ with diameter $R=\sup_{x,y\in \mathcal{U}(x_{\ast})} d(x,y)<\infty$ such that for all $x\in \mathcal{U}(x_{\ast})$, the gradients are bounded by $\|\grad_{x}f(x)\|\leq M_{f}$, $\|\grad_{x}h_{k}(x)\|\leq M_{h}$, and $|h_{k}(x)|\leq G,k\in \{1,2,\cdots,m\}$ for some positive bounded values $M_{f},M_{h}$, and $G$.
\end{assumption}
In the case that $\mathcal{M}$ is compact, we simply take $\mathcal{U}(x_{\ast})=\mathcal{M}$ in Assumption \ref{Assumption:1}.  

For a non-asymptotic convergence result, we also require the following assumption:
\begin{assumption}
\label{Assumption:g-convex}
The functions $f(x)$ and $h_{k}(x),k=1,2,\cdots,m$ are geodesically convex on the manifold $\mathcal{M}$ in the sense of Definition \ref{Def:geodesic_convexity}.
\end{assumption}

Based on Assumption \ref{Assumption:g-convex}, any local minima of a geodesically convex function over $\mathcal{M}$ is a global minima. The Kronecker products and the logarithms of determinants are two examples of functions that are geodesically convex on the cone of positive definite matrices, see \cite{wiesel2012geodesic} for some applications of this result in the context of covariance estimation.  

We note that the more familiar notion of convexity in the Euclidean space follows naturally from the geodesic convexity condition described in Definition \ref{Def:geodesic_convexity}. In particular, the geodesic path between two arbitrary points $p$ and $q$ in the Euclidean space is a straight line $\gamma(t)=(1-t)p+tq$ for $t\in [0,1]$, which leads to the usual notion of convexity condition, namely the line segment joining any two points on the graph of the function $(p,f(p))$ lies above the graph. 

\begin{algorithm}[t!]
\caption{Stochastic Primal-Dual Method on the Riemannian Manifolds}
\label{alg:1}
\begin{algorithmic}[1]
\State{\textbf{initialize}: $x_{0}\in \mathcal{U}(x_{\ast})\subseteq \mathcal{M}$, $\lambda_{0}=0$, and a decreasing sequence $\{\eta_{t}\}_{t=0}^{T}$ as the step size $\eta_{t}\in \real_{+}$.}
\For{$t=0,1,\cdots,T$}
\State{Query $x_{t}$ and receive $\grad_{x_{t}} F(x_{t};\xi_{t})$ where $\xi_{t}\sim P$.}
\State{Update $x_{t+1}$ and $\lambda_{t+1}$ as follows}
\State{$x_{t+1}\leftarrow \Exp_{x_{t}}\left(-\eta_{t}\grad_{x_{t}} \mathfrak{L}(x_{t},\lambda_{t};\xi_{t})\right).$}
\State{$\lambda_{t+1}\leftarrow \Pi_{\real_{+}^{m}}\left(\lambda_{t}+\eta_{t}\grad_{\lambda_{t}} \mathfrak{L}(x_{t},\lambda_{t}) \right) $.}
\EndFor

\hspace{-7mm}\Return $\{x_{t}\}_{t=1}^{T}$.
\end{algorithmic}
\end{algorithm}

Associated with the stochastic gradient vector, we define the estimation error
\begin{align*}
e_{t}\ \dot{=} \ \grad_{x_{t}} f (x_{t})-\grad_{x_{t}} F(x_{t},\xi_{t}).
\end{align*}

\begin{assumption}
\label{Assumption:unbias_error}
We assume that the estimations are unbiased, \textit{i.e.}, $\expect[e_{t}|\mathfrak{F}_{t}]=0$ for all $t\in [T]$, where $\mathfrak{F}_{t}$ is the $\sigma$-algebra of the random variables $\xi_{0},\cdots,\xi_{t-1}$. 
\end{assumption}

Moreover, when we prove the high probability bounds, the estimation error $e_{t}$
must have a distribution whose tail decays sufficiently fast. To make this statement more precise, we require the following definition:

\begin{definition}
\label{Def:Orlicz}
The Young-Orlicz modulus is a convex non-decreasing function $\psi:{\rm I\!R}_{+}\rightarrow {\rm I\!R}_{+}$ such that $\psi(0)=0$ and $\psi(x)\rightarrow \infty$ when $x\rightarrow \infty$. Accordingly, the Orlicz norm of an integrable random variable $X$ with respect to the modulus $\psi$ is defined as
\begin{align*}
\|X\|_{\psi}\topdoteq\inf\left\{\beta >0:{\rm I\!E}\left[\psi(|X|/\beta)\right]\leq 1\right\}.
\end{align*}
\end{definition}

Among the classical examples of the Young-Orlicz moduli are $\psi_{p}(x)\topdoteq x^{p},p\geq 1$ and $\psi_{\nu}(x)\topdoteq\exp(x^{\nu})-1$. Throughout the paper, we consider the modulus $\psi_{\nu}(x)$ $\topdoteq\exp(x^{\nu})-1$. Accordingly, we say $X$ is \textit{sub-Gaussian} when $\|X\|_{\psi_{2}}<\infty$, and  \textit{sub-exponential} when $\|X\|_{\psi_{1}}<\infty$. See Appendix \ref{sec:Heavy_Tail}.

\vspace{4mm}
\textsc{Notation}. A list of symbols used in this paper is provided in Table \ref{Table:1} at the end of this paper. We use the standard asymptotic notation. In particular, if $x_{n}$ and $y_{n}$ are positive sequences, then $x_{n} = \mathcal{O}(y_{n})$ means that $\limsup_{n\rightarrow \infty} {x_{n}\over y_{n}} < \infty$. We denote both vectors and matrices by $\bm{X}\equiv (X_{1},X_{2},\cdots,X_{n})$ and $\bm{X}=[X_{ij}]$, respectively. For two vectors $\bm{X}$ and $\bm{Y}$, the vector inequality $\bm{X}\preceq \bm{Y}$ means the element-wise inequality, \textit{\textit{i.e.}}, $X_{i}\leq Y_{i}$ for all $i=1,2,\cdots,n$. Moreover, for (non-negative) sequences $x_{i}$ and $y_{i}$, $x_{i}\lesssim y_{i}$ means that there exists a constant $c<\infty$ such that $x_{i}\leq c\cdot y_{i}$ for all $i\in \{1,2,\cdots,n\}$. For a scalar $a\in \real$ we define $[a]_{+}\topdoteq\max\{0,a\}$. We denote the projection of the vector $\bm{X}$ onto the closed set $\mathcal{X}$ by $\Pi_{\mathcal{X}}(\bm{X})\equiv\arg \min_{\bm{X}\in \mathcal{X}}\|\bm{X}-\bm{Y}\|_{2}$. We often use the shorthand $[T]\topdoteq\{1,2,\cdots,T\}$ to denote a set. Lastly, for vectors $\bm{X}$ defined in the Euclidean space and manifolds, $\|\bm{X}\|$ denotes the Euclidean $2$-norm, and the norm defined by the metric on the manifold, respectively.

\subsection{Convergence on Hyperbolic Manifolds}
\label{subsec:Hyperbolic}
We first consider the case of hyperbolic manifolds. In particular, throughout this subsection, we make the following assumption:

\begin{assumption}
\label{Assumption:hyper}
At every point on the manifold $p\in \mathcal{M}$, the sectional curvature is negative and lower bounded by $\kappa$, \textit{i.e.}, $\kappa\leq K(X_{p},Y_{p})\leq 0$ for all $X_{p},Y_{p}\in T_{p}\mathcal{M}$.
\end{assumption}

The Euclidean space $\real^{n}$ is a simple example of a manifold that has zero sectional curvature everywhere $\kappa=0$ and hence satisfies the condition of Assumption \ref{Assumption:hyper}. A non-trivial example of a hyperbolic manifold is the statistical manifold given by $\{p(\cdot|\theta)\}_{\theta\in \Theta}$ under the Fisher information metric \cite{lafferty2005diffusion}. Here, $\theta=(\mu,\sigma)$, $\Theta=\real^{n-1}\times \real_{+}$, and $p(\cdot|\mu,\sigma)=\mathsf{N}(\mu,\sigma^{2}I_{n-1})$ is the multivariate Gaussian having mean $\mu\in \real^{n-1}$, and variance $\sigma^{2}I_{n-1}$ with $\sigma>0$.

Under Assumption \ref{Assumption:hyper}, the following results are well-known:

\begin{itemize}
\item[(\textsc{A}1)] For any point $p\in \mathcal{M}$, there exists a unique distance minimizing geodesic between $p$ and $q$ for all $q\in \mathcal{M}$. 
\item[(\textsc{A}2)] For any two points $p,q\in \mathcal{M}\times \mathcal{M}$, the squared Riemannian distance function $d^{2}(p,q)$ is a continuously differentiable function for all $q\in \mathcal{M}/\{p\}$, \textit{i.e.}, $d^{2}\in C^{2}(\mathcal{M}\times \mathcal{M},\real_{+})$.

\item[(\textsc{A}3)] The exponential map is a diffeomorphism from the tangent bundle $T\mathcal{M}$ $\topdoteq\cup_{p}\{p\}\times T_{p}\mathcal{M}$ onto the manifold $\mathcal{M}$. Hence, the logarithmic map $\log_{p}:\mathcal{M}\rightarrow T_{p}\mathcal{M}$ exists. 
\end{itemize}

The proofs of these assertions are based on the fact that on manifolds with the negative sectional curvature, the cut-locus\footnote{The cut-locus is the set of cut points of $p$, namely the points $q$ along the emanating geodesic $\gamma:[0,\infty)\rightarrow \mathcal{M}, \gamma(0)=p$ such that if $\gamma(t_{0})=q$ and $\gamma(t)=\tilde{q}, t>t_{0}$, there exists a geodesic between $p$ and $\tilde{q}$ shorter than $\gamma|_{[0,t]}$. In this case, a geodesic between $p$ and $q$ stops to be distance minimizing beyond $q$.} as well as the conjugate-locus\footnote{The conjugate locus is the set of conjugate points of $p\in \mathcal{M}$, \textit{i.e.}, the points $q\in \mathcal{M}$ on the manifold that are connected by a 1-parameter family of geodesics to $p$. In this case, a geodesic between $p$ and $q$ is not \textit{locally} distance minimizing.} of $p$ are empty; cf. \cite[Lemma 19.1]{milnor2016morse}.

With these technical details in place, we are now ready to state the first main theorem of this paper:
\begin{theorem}
\label{Thm:1}
Suppose Assumptions \ref{Assumption:1},\ref{Assumption:g-convex}, \ref{Assumption:unbias_error}, and \ref{Assumption:hyper} hold. Let $\{x_{t},\lambda_{t}\}_{t=1}^{T}$ be a sequence of the primal-dual points generated by Algorithm \ref{alg:1} with the step size $\eta_{t}<1/\alpha$, such that $(x_{t},\lambda_{t})\in \mathcal{U}(x_{\ast})\times \real_{+}$ for all $t=0,1,\cdots,T$. Then,
\begin{align*}
&\min_{t\in [T]}\expect[f(x_{t})]-f(x_{\ast})\\ \leq &\dfrac{1}{2}R^{2}+{1\over 2}\sum_{t=0}^{T-1}\Big(A+(1+{2R\sqrt{|\kappa|/2}})\expect[\|e_{t}\|^{2}]\Big){\eta_{t}^{2}}+{4 \sqrt{|\kappa|}\over 3\sqrt{2}}\sum_{t=0}^{T-1}(B+\expect[\|e_{t}\|^{3}])\eta_{t}^{3}.
\end{align*}
and the constants $A$ and $B$ are defined as
\begin{align}
\label{Eq:Def_A}
A\ &\dot{=}\ 4mG^{2}+M^{2}\left(1+\dfrac{\sqrt{m}G}{\alpha}\right)^{2},\\
\label{Eq:Def_B}
B\ &\dot{=}\ M^{3}\left(1+\dfrac{\sqrt{m}G}{\alpha}\right)^{3}.
\end{align}
\end{theorem}
Notice that in Theorem \ref{Thm:1}, it is essential that the sequence of the primal-dual points must remain inside the sub-manifold $\mathcal{U}(x_{\ast})$, as otherwise, the gradients of the Lagrangian function can become unbounded. Similar to the Euclidean optimization algorithms, projection of the primal variables to the sub-manifold $\mathcal{U}(x_{\ast})$ guarantees the boundedness of the gradients.

It is worthwhile to compare the convergence bound described in Theorem \ref{Thm:1} with that of the stochastic optimization over the Euclidean space $(\kappa=0)$. We observe that the convergence rate on the hyperbolic manifolds depends on the second and third moments of the gradient estimation error, \textit{i.e.}, $\expect[\|e_{t}\|^{2}]$ and $\expect[\|e_{t}\|^{3}]$, respectively. In comparison, the convergence rate of the Euclidean optimization techniques only depends on the variance of the estimation error, also see \cite{duchi2016variance}.

Moreover, similar to the regularized primal-dual method in the Euclidean space, the value of the regularization parameter $\alpha$ determines the rate at which the constraints violation decays to zero, see \cite{khuzani2016distributed}. In particular, we have the following corollary:
\begin{corollary}
\label{Cor:2}
Consider the conditions of Thm. \ref{Thm:1} and let $e_{t}=0,t\in [T]$ for simplicity (deterministic algorithm). Further, choose the step size as $\eta_{t}={1\over \sqrt{t+1}}$ and let $\alpha\eta_{t}\leq 1$ for all $t=0,1,\cdots,T$ (thus $\alpha\leq 1$). Using the asymptotic notations, the convergence rate is  
\begin{align*}
\min_{t\in [T]}\expect[f(x_{t})]-f(x_{\ast})=\mathcal{O}\left(\dfrac{\log(T)}{\alpha^{2}(\sqrt{T}-1)}+\dfrac{1}{\alpha^{3}(\sqrt{T}-1)} \right),
\end{align*}
and the constraint violation is
\begin{align*}
\left\| \left[\sum_{t=0}^{T-1}\widehat{\eta}_{t}h(x_{t})\right]_{+}\right\|_{2}^{2}=\mathcal{O}\left(\alpha\right),
\end{align*}
where $\widehat{\eta}_{t}\topdoteq \eta_{t}/{\sum_{t=0}^{T-1}\eta_{t}}$ is the normalized step size. In particular, when $\alpha=1/T^{{1-\beta\over 6}}$ for $\beta\in (0,1)$, the rate of convergence and constraint violation are $\mathcal{O}\big(1/{T^{\beta\over 2}}\big)$ and $\mathcal{O}\big({1/T^{{1-\beta\over 6}}}\big)$, respectively.
\end{corollary}

From Corollary \ref{Cor:2}, we observe that there is a tension between the rate of convergence and the decay rate of the constraint violation. Specifically, increasing $\alpha$ results in smaller dual variables. This in turn results in a larger constraint violation since the magnitude of the dual variables controls the penalty of violating the inequality constraints. 

Also, observe that in the asymptotic convergence analysis of stochastic algorithms, the step size of the algorithm must satisfies $\sum_{t=1}^{\infty}\eta_{t}=\infty$ and $\sum_{t=1}^{\infty}\eta_{t}^{2}<\infty$, \cite{zhang2008impact}. As we observe from Corollary \ref{Cor:2}, the step size in Algorithm \ref{alg:1} only needs to have a converging ratio $\sum_{t=1}^{T}\eta^{2}_{t}/\sum_{t=1}^{T}\eta_{t}\rightarrow 0$ as $T\rightarrow \infty$. For instance, the step size $\eta_{t}=1/\sqrt{t+1}$ in Corollary \ref{Cor:2} is diverging $\sum_{t=1}^{\infty}\eta_{t}=\infty$ and $\sum_{t=1}^{\infty}\eta_{t}^{2}=\infty$, but the ratio is converging.

\subsection{Convergence on the Elliptic and Asymptotically Elliptic Manifolds}
In this section, we focus on the class of manifolds with positive and asymptotically non-negative sectional curvatures. First, we consider the following assumption:

\begin{assumption}
\label{Assumption:4}
At every point on the manifold $p\in \mathcal{M}$, the sectional curvature is positive and bounded from below by $\kappa>0$, \textit{i.e.}, $K(X_{p},Y_{p})\geq \kappa$ for all $X_{p},Y_{p}\in T_{p}\mathcal{M}$.
\end{assumption}

The compact Lie groups, such as $\rm{SU}(2)$ and $\rm{SO}(3)$, the compact Grassmann manifold ${\rm{G}}(n,1)={{\rm O}(n)}/{{{\rm O}(1)}\times {{\rm O}(n-1)}}$ (one dimensional subspaces
of $\real^{n}$), and a sphere are the notable examples of this class of manifolds. 

In contrast to the case of hyperbolic manifolds, the exponential map fails to be a local diffeomorphism on the entire tangent plan $T_{p}\mathcal{M}$ on the elliptic manifolds. Nevertheless, for any given point $p\in \mathcal{M}$, there exists a neighborhood ${\rm I\!B}_{r}(0_{p})\topdoteq\{X_{p}\in T_{p}\mathcal{M}:g_{p}(X_{p},X_{p})<r\}\subset T_{p}\mathcal{M}$ around $0_{p}\in T_{p}\mathcal{M}$ for some radius $r\in \real_{+}$ such that the restriction $\Exp_{p}\big\vert_{{\rm I\!B}_{r}(0_{p})}$ is a diffeomorphism onto its range $\Exp_{p}({\rm I\!B}_{r}(0_{p}))\subset \mathcal{M}$, cf. \cite[Lemma. 5.10]{lee2009manifolds}. 

The following definition is concerned with the largest value of such $r$ over the entire manifold $\mathcal{M}$.
\begin{definition}
The injectivity radius of $\mathcal{M}$ is defined as
\begin{align*}
i(\mathcal{M})\topdoteq\inf_{p\in \mathcal{M}}i(p),
\end{align*}
where 
\begin{align*}
i(p)\topdoteq\sup \{r\in \real_{+}:\Exp_{p}\big\vert_{{\rm I\!B}_{r}(0_{p})}\text{is a diffeomorphism}\}.
\end{align*}
\end{definition}

The next theorem establishes the convergence bound of the stochastic primal-dual method on the elliptic manifolds. 
Since the elliptic manifolds of interest in this paper are compact, in the following we assume that the gradients of the the Lagrangian function in Assumption \ref{Assumption:1} are bounded on the entire manifold, that is $\mathcal{U}(x_{\ast})=\mathcal{M}$.

\begin{theorem}
\label{Thm:2}
\textsc{(Non-Asymptotic Convergence)} Suppose Assumptions \ref{Assumption:1}, \ref{Assumption:g-convex}, \ref{Assumption:unbias_error}, and \ref{Assumption:4} hold. Let $\{x_{t},\lambda_{t}\}_{t=1}^{T}$ be the sequence of the primal-dual points generated by Algorithm \ref{alg:1} with the step size $\eta_{t}\leq  \min\{1/\alpha, i(\mathcal{M})/\|X_{t}\|\}$. Then, for all the primal-dual pairs $(x_{t},\lambda_{t})\in \mathcal{M}/\{\mathcal{C}(x_{\ast})\}\times \real_{+},t=1,2,\cdots,T$, we have
\begin{align*}
&\min_{t\in [T]}\expect\left[\left(f(x_{t})-f(x_{\ast})\right)\cdot \sinc(\sqrt{\kappa}d(x_{t},x_{\ast}) ) \right]\\ &\leq \dfrac{1}{\sum_{t=0}^{T-1}\eta_{t}}\bigg[ \dfrac{2}{\kappa}\sin^{2}\left({\sqrt{\kappa}\over 2} d(x_{0},x_{*})\right)+\sum_{t=0}^{T-1}(A+\expect[\|e_{t}\|^{2}])\eta_{t}^{2}\bigg],
\end{align*}
where $\sinc(x)\topdoteq\sin(x)/x$, $X_{t}\topdoteq -\grad_{x_{t}}\mathfrak{L}(x_{t},\lambda_{t};\xi_{t})$, and $\mathcal{C}(x_{\ast})\subset \mathcal{M}$ denotes the conjugate locus of the point $x_{\ast}$ (cf. Section \ref{subsec:Hyperbolic}).
\end{theorem}

Let us emphasize a few important points about Theorem \ref{Thm:2}.

First, note that due to Myers' theorem (cf. Thm. \ref{Thm:Meyer}), $d(x_{t},x_{*})\leq \pi/\sqrt{\kappa}$. Due to the assumption $x_{t}\not\in \mathcal{C}(x_{\ast})$, we indeed have the strict inequality $d(x_{t},x_{*})< \pi/\sqrt{\kappa}$. To observe this, note that when $d(x_{t},x_{\ast})=\pi/\sqrt{\kappa}$, the manifold $\mathcal{M}$ has a diameter $\pi/\sqrt{\kappa}$ and thus by Cheng's rigidity result \cite{cheng1975eigenvalues}, $\mathcal{M}$ is isometric to the sphere of radius $\sqrt{\kappa}$. In this case, $x_{t}$ and $x_{\ast}$ are conjugate (antipodal points) which is excluded by the statement of Theorem \ref{Thm:2}. Hence, $d(x_{t},x_{\ast})<\pi/\sqrt{\kappa}$ and $\sinc(\sqrt{\kappa}d(x_{t},x_{*}))\in (0,1]$. 

Second, as the curvature lower bound diminishes $\kappa\downarrow 0$, $\sinc(\kappa d(x_{t},x_{*}))\uparrow 1$ and $\dfrac{2}{\kappa}\sin^{2}\left({\sqrt{\kappa}\over 2} d(x_{0},x_{*})\right)\rightarrow \dfrac{1}{2}d^{2}(x_{0},x_{*})$. We thus observe that the convergence bounds we obtained in Theorems \ref{Thm:1} and \ref{Thm:2} are in agreement (up to a constant factor of 1/2 in the second term) in the asymptotic of vanishing curvature $|\kappa|\rightarrow 0$. This consistency is appealing in light of the fact that different proof techniques are used to analyze the convergence rates for elliptic and hyperbolic manifolds.

Third, note that although computing the exact value of the injectivity radius $i(\mathcal{M})$ is challenging in general, if the sectional curvature is further bounded from above $0<\kappa\leq K(X_{p},Y_{p})\leq \delta$ for all $X_{p},Y_{p}\in T_{p}\mathcal{M}$, then Klingenberg \cite{klingenberg1995riemannian} has shown that $\mathcal{M}$ is compact and either $i(\mathcal{M})\geq \pi/\sqrt{\delta}$, or there is a closed geodesic $\gamma$ of minimal length among all closed geodesics in $\mathcal{M}$ such that $i(\mathcal{M})={1\over 2}L[\gamma]$. By establishing a lower bound on $L[\gamma]$,  the Cheeger-Gromov-Taylor inequality provides a more refined lower bound (cf. \cite[Theorem 4.3.]{cheeger1982finite}),
\begin{align*}
i(p)\geq \min\left\{\pi/\sqrt{\delta}, \sup_{0<r<{\pi/4\sqrt{\delta}}} r\left(1+\dfrac{\text{Vol }({\rm I\!B}_{2r}(0_{p}))}{\text{Vol }(\widetilde{\rm I\!B}_{r}(p))} \right)^{-1}\right\},\quad \text{for all}\ p\in \mathcal{M},
\end{align*}
where $\text{Vol}({\rm I\!B}_{2r}(0_{p}))$ denotes the volume of the the ball ${\rm I\!B}_{2r}(0_{p})\in T_{p}\mathcal{M}$ with respect to the pull-back metric via the exponential map, \textit{i.e.}, $g_{p}^{\ast}(X_{p},Y_{p})$ $\topdoteq\Exp^{\ast} g_{p}(X_{p},Y_{p})$ $\dot{=}\ g_{p}(d\Exp(X_{p}), d\Exp(Y_{p}))$ for all $X_{p},Y_{p}\in T_{p}\mathcal{M}$. 

Lastly, we note that in the Euclidean optimization methods, the distance between the initial and optimal points $d(x_{0},x_{\ast})$ often appears as a quadratic term in the upper bound and hence greatly influences the convergence speed, \textit{e.g.} see \cite{nedic2009subgradient}. In contrast, from the upper bounds in Theorem \ref{Thm:2} we observe that the primal-dual optimization on the elliptic manifolds is less sensitive to the initialization point. 

In the following theorem, we characterize a high probability convergence bound when the tail of the estimation error $\|e_{t}\|$ decays as a sub-Gaussian random variable.
\begin{theorem}
\label{Thm:3}
\textsc{(High-Probability Convergence)} Consider the assumptions and the step size $\eta_{t}$ of Thm. \ref{Thm:2}. Further, suppose $\vertiii{e_{t}}_{\psi_{2}}=\beta<\infty$. Then, there exists a constant $c>0$, such that with the probability of at least $1-2\varrho$,
\small
\begin{align*}
&\min_{t\in [T]}\left(f(x_{t})-f(x_{\ast})\right)\cdot \sinc(\sqrt{\kappa}d(x_{t},x_{\ast})) \\ &\leq \dfrac{1}{\sum_{t=0}^{T-1}\eta_{t}}\Bigg[ \dfrac{2}{\kappa}\sin^{2}\left({\sqrt{\kappa}\over 2} d(x_{0},x_{*})\right)
+\left(A+4\beta^{2}\right)\sum_{t=1}^{T-1}\eta_{t}^{2}+{2\sqrt{\kappa} \beta \log^{1\over 2}({2/\varrho)}\over c^{1\over 2}}\left(\sum_{t=0}^{T-1}\eta_{t}^{2}\right)^{1\over 2}\\
&\hspace{20mm} +12\beta^{2}\max\bigg\{{\log^{1\over 2}({2/\varrho})\over c^{1\over 2}}\left(\sum_{t=0}^{T-1}\eta_{t}^{4}\right)^{1\over 2},{\log({2/\varrho})\over c}\sum_{t=0}^{T-1}\eta_{t}^{2}\bigg\}\Bigg],
\end{align*}
\normalsize
for all $x_{t}\in \mathcal{M},t=1,2,\cdots,T$.
\end{theorem}

Although the statement of Theorem \ref{Thm:3} only includes the estimation errors with a sub-Gaussian tail, the proof we present is somewhat more general and includes those cases where the norm of the estimation error $\|e_{t}\|$ has a tail decaying faster than a sub-Gaussian, \textit{i.e.}, when $\vertiii{e_{t}}_{\psi_{\nu}}= \beta<\infty$ for all $\nu\geq 2$. In those cases, a tighter bound can be proved, cf. Section \ref{subsec:Proof_of_Thm3}.

We now replace Assumption \ref{Assumption:4} with a milder assumption. In particular, Assumption \ref{Assumption:4} is rather stringent as it imposes a strict positivity constraint on the sectional curvature. We thus relax this assumption by requiring that the sectional curvature to be locally bounded from below by the inverse of the quadratic distance function. 

\begin{assumption}
\label{Assumption:5}
At every point on the manifold $p\in \mathcal{M}$, the sectional curvature is lower bounded from below $K(X_{p},Y_{p})\geq -\kappa/d^{2}(p,q)$ for all $X_{p},Y_{p}\in T_{p}\mathcal{M}$, where $q\in \mathcal{M}$ is an arbitrary point on the manifold, and $\kappa>0$ is a constant.
\end{assumption}

The class of manifolds that satisfy Assumption \ref{Assumption:5} is broader than Assumption \ref{Assumption:4} and subsumes the elliptic manifolds. In the sequel, we refer to the class of manifolds satisfying Assumption \ref{Assumption:5} as the \textit{asymptotically elliptic} since $K(X_{p},Y_{p})\geq 0$ as $d(p,q)\rightarrow \infty$. 

It is meaningful to choose a stationary point of Problem \ref{Problem} as the reference point of Assumption \ref{Assumption:5}, \textit{i.e.}, $q=x_{\ast}$. In this case, the condition of Assumption \ref{Assumption:5} is locally satisfied by all the smooth manifolds in the sense that if the algorithm outputs $x_{t}$ remain in a small neighborhood of the stationary point $x_{\ast}$ for all $t\in [T]$, then $d(x_{t},x_{\ast})$ is very small and thus the lower bound $-\kappa/d^{2}(x_{t},x_{\ast})$ is also a small negative number. From this point of view, Assumption \ref{Assumption:5} is more appealing than Assumption \ref{Assumption:4}. We now have the following result:

\begin{theorem}
\label{Thm:5}
\textsc{(Non-Asymptotic Convergence)} Suppose Assumptions \ref{Assumption:1}, \ref{Assumption:g-convex}, and \ref{Assumption:5} hold. Let $\{x_{t}\}_{t=1}^{T}$ be a finite sequence of iterates generated by Algorithm \ref{alg:1} with a step size satisfying $\eta_{t}\alpha<1$. Then,  
\begin{align*}
&\min_{t\in [T]}\expect[f(x_{t})]-f(x_{\ast})\\ &\leq \dfrac{1}{\sum_{t=0}^{T-1}\eta_{t}}\Bigg[\dfrac{1}{2}d^{2}(x_{0},x_{\ast})+{1\over 2}(1+\sqrt{1+4\kappa^{2}}) \sum_{t=0}^{T-1}(A+\expect[\|e_{t}\|^{2}]){\eta_{t}^{2}}\Bigg],
\end{align*}
for all $x_{t}\in \mathcal{M},t=1,2,\cdots,T$.
\end{theorem}

Notice that when $\kappa\downarrow 0$, the condition of Assumption \ref{Assumption:5} is the same as Assumption \ref{Assumption:4}, namely $K\geq 0$. In this limit, the upper bounds in Theorems \ref{Thm:3} and Theorem \ref{Thm:5} are the same up to a constant factor 1/2 in the second term, which confirms our results.

\section{Applications}
\label{Applications}
In this section, we evaluate the performance of Algorithm \ref{alg:1} for a few applications. \footnote{The MATLAB codes for all simulations can be found in our Github repository: \url{https://github.com/CoNG-Harvard/Optimization_on_Manifolds.}} In particular, we consider (\textit{i}) the non-negative online principle component analysis (PCA), (\textit{ii}) the anchored synchronization on the rotation group ${\rm SO}(3)$, and (\textit{iii}) the weighted MAX-CUT problem under graph constraints. In each example, the optimization is carried over a sub-manifold of a positively curved Riemannian manifold. 

\subsection{Non-negative Online PCA}
\label{sec:Non-negative_PCA}
\subsubsection{Problem Description}

The problem formulation we consider here is due to \cite{montanari2016non}. The principal component analysis (PCA) is a popular dimension-reduction technique that deals with estimating the direction of maximal variability in given zero-mean samples $\bm{\xi}_{0},\cdots,\bm{\xi}_{T-1}\in \real^{d}$ from a $d$-dimensional probability distribution $P(\bm{\xi})$. In particular, PCA can be formulated as the following non-convex optimization on the unit sphere,  
\begin{subequations}
\begin{align}
\label{Eq:population}
&\max_{\bm{x}}\ \expect_{P}[\langle\bm{x},\bm{\xi} \rangle^{2}]=\int_{\Xi}\langle\bm{x},\bm{\xi} \rangle^{2}dP(\bm{\xi})\\
\label{Eq:population1}
&\text{Subject to:} \bm{x}\in {\rm S}^{d},
\end{align}
\end{subequations}
where $\bm{\xi}$ is a random vector distributed with the distribution of $P(\cdot)$ and the support $\Xi$, and ${\rm S}^{d}\subset \real^{d+1}$ is the unit sphere
\begin{align*}
{\rm S}^{d}\topdoteq\left\{(x_{0},\cdots,x_{d})\in \real^{d+1}: \sum_{i=0}^{d}x_{i}^{2}=1\right \}.
\end{align*}
The solution of eqs. \eqref{Eq:population}-\eqref{Eq:population1} corresponds to the principal eigenvector of the covariance matrix $\expect(\bm{\xi}\bm{\xi}^{T})$ which we denote by $\bm{\xi}^{\ast}$.

In practice, the underlying distribution $P(\bm{\xi})$ is unknown and the problem in eqs. \eqref{Eq:population}-\eqref{Eq:population} cannot be solved directly. The standard approach is thus to approximate the population objective \eqref{Eq:population} with the following empirical cost function 
\begin{subequations}
\begin{align}
\label{Eq:emperical}
&\max_{\bm{x}}\ \expect_{\hat{P}_{T}}[\langle\bm{x},\bm{\xi} \rangle^{2}]={1\over T}\sum_{t=0}^{T-1}\langle \bm{x},\bm{\xi}_{t} \rangle^{2}, \\
\label{Eq:emperical2}
&\text{Subject to:} \bm{x}\in {\rm S}^{d},
\end{align}
\end{subequations}
where $\hat{P}_{T}$ is the empirical measure. The solution of eqs. \eqref{Eq:emperical}-\eqref{Eq:emperical2} corresponds to the principal eigenvector of the sample covariance matrix $(1/T)\sum_{t=0}^{T-1}\bm{\xi}_{t}\bm{\xi}_{t}^{T}$ which we denote by $\hat{\bm{\xi}}_{T}\topdoteq\hat{\bm{\xi}}(\bm{\xi}_{0},\cdots,\bm{\xi}_{T-1})$. In the low dimension regime $T/d\rightarrow \infty$, PCA is successful in recovering the population eigenvector in the sense that $\|\hat{\bm{\xi}}_{T}-\bm{\xi}_{\ast}\|_{2}\rightarrow 0$ with a high probability as $T\rightarrow \infty$ \cite{anderson1963asymptotic}. However, it is well-known that when the dimension is comparable to the sample size $T=\mathcal{O}(d)$, PCA undergoes a phase transition. Specifically, in the limit of $T,d\rightarrow \infty,T/d\rightarrow \alpha\in (0,1)$ the phase transition threshold is $\sqrt{\alpha}$ below which $\lim_{T\rightarrow \infty}|\langle \hat{\bm{\xi}}_{T},\bm{\xi}_{\ast} \rangle|=0$. That is, the empirical eigenvector $\hat{\bm{\xi}}_{T}$ is uninformative about the population eigenvector $\bm{\xi}_{\ast}$. Clearly, without any further structure, recovering the high-dimensional principle component is challenging. Thus, to increase the phase transition threshold, some side information about the population component $\bm{\xi}_{\ast}$ is normally considered. For example, $\bm{\xi}_{\ast}$ is either sparse \cite{zou2006sparse,ma2013sparse,journee2010generalized,grbovic2012sparse}, resides inside a cone \cite{deshpande2014cone}, or is simply a non-negative vector \cite{montanari2016non}. In the sequel, we focus on the simple case of non-negative PCA which arises naturally in the non-negative matrix factorization problem. In particular, we consider the non-negative PCA problem under the \textit{symmetric spiked covariance model} with a single spike \cite{johnstone2012consistency,montanari2016non},
\begin{subequations}
\begin{align}
\label{Eq:spike1}
&\max_{\bm{x}}\ \expect[\langle\bm{x},\bm{\xi} \rangle^{2}]\\
\label{Eq:spike2}
&\text{Subject to:}\ \bm{x}\in {\rm S}^{d},\quad \bm{x}\succeq 0,
\end{align}
\end{subequations}
where the samples $\bm{Y}\topdoteq(\bm{\xi}_{0},\cdots,\bm{\xi}_{T-1})^{T}\in \real^{T\times d}$ are generated according to the following model
\begin{align}
\label{Eq:SpikedModel}
\bm{Y}=\sqrt{\textsf{SNR}}\bm{\xi}_{\ast}\bm{\xi}^{T}_{\ast}+\bm{Z}.
\end{align}
Here, $\bm{\xi}_{\ast}\succeq 0, \|\bm{\xi}_{\ast}\|_{2}=1$, and $\bm{Z}=[Z_{ij}]$ is the symmetric noise $Z_{ij}=Z_{ji}$, where $Z_{ij}\sim \mathsf{N}(0,1/T)$ are i.i.d. additive noise, and $Z_{ii}\sim \mathsf{N}(0,2/T)$. 

\begin{remark}
Maximizing the empirical cost function as in eq. \eqref{Eq:emperical} leads to an offline algorithm in which all the samples $\bm{\xi}_{0},\cdots,\bm{\xi}_{T-1}$ are available. In contrast, in Algorithm \ref{alg:1}, the samples $\bm{\xi}_{0},\cdots,\bm{\xi}_{T-1}$ are observed sequentially. That is, the Riemannian primal-dual method we develop in the sequel is an online algorithm.
\end{remark}

\begin{algorithm}[t!]
\caption{Riemannian Primal-Dual Method for Non-negative Online PCA}
\begin{algorithmic}[1]
\State{\textbf{require}:  A sequence of samples $\bm{\xi}_{0},\bm{\xi}_{1},\cdots,$ from eq. \eqref{Eq:spike2}}.
\State{\textbf{initialize}: Choose $\hat{\bm{x}}_{0}\in \mathsf{N}(0,\bm{I}_{d})$ and set $\bm{x}_{0}=\hat{\bm{x}}_{0} /\|\hat{\bm{x}}_{0}\|$. Let $\bm{\lambda}_{0}=0$. A decreasing sequence $\{\eta_{t}\}_{t=0}^{T-1}$ for the step size as in Thm. \ref{Thm:2}. The regularizer free parameter $\alpha\in \real_{+}$.}
\For{$t=0,1,2,\cdots$}
\State{Query the oracle and receive $\bm{\xi}_{t}$.}
\State{Compute $\bm{X}_{t}=2\langle \bm{x}_{t},\bm{\xi}_{t}\rangle \bm{\xi}_{t}+\bm{\lambda}_{t}$.} 
\State{Update variables as follows:}
\State{$\bm{x}_{t+1}\leftarrow \cos(\eta_{t}\|\bm{X}_{t}\|)\bm{x}_{t}+\sin(\eta_{t}\|\bm{X}_{t}\|){\bm{X}_{t}\over\|\bm{X}_{t}\|}$.}
\State{$\bm{\lambda}_{t+1}\leftarrow \Pi_{\real^{d}_{+}}\big((1-\alpha)\bm{\lambda}_{t}-\eta_{t}\bm{x}_{t} \big).$}

\hspace{-7mm}\Return $\hat{\bm{\xi}}_{t+1}=\bm{x}_{t+1}$.
\EndFor
\end{algorithmic}
\label{alg:3}
\end{algorithm}

\subsubsection{Algorithms}
\paragraph{Riemannian Primal-Dual Method}
Since the Riemannian optimization techniques for the classical PCA problem in eqs. \eqref{Eq:emperical}-\eqref{Eq:emperical2} are well-established \cite{absil2009optimization}, here we mostly summarize the results for the primal-dual method in Algorithm \ref{alg:1}. In particular, for each point $\bm{x}\in {\rm S}^{d}\subset \real^{d+1}$, the tangent space $T_{\bm{x}}{\rm S}^{d}$ is given by
\begin{align}
\label{Eq:TangentSphere}
T_{\bm{x}}{\rm S}^{d}=\{\bm{z}\in \real^{d+1}:\langle \bm{x},\bm{z} \rangle=0\}.
\end{align}
Given a vector $\bm{z}\in \real^{d+1}$, the orthogonal projection onto the tangent space $T_{\bm{x}}{\rm S}^{d}$ is $\Pi_{\bm{x}}(\bm{z})=\bm{z}-\bm{x}\langle \bm{x},\bm{z} \rangle$. To verify that this is indeed a projection, notice that $\langle \bm{x},\Pi_{\bm{x}}(\bm{z}) \rangle =\langle \bm{x},\bm{z}\rangle-\langle \bm{x},\bm{x} \rangle \langle \bm{x},\bm{z}\rangle=0$ and thus $\Pi_{\bm{x}}(\bm{z})\in T_{\bm{x}}{\rm S}^{d}$. 

It is a standard exercise to show that the geodesic of ${\rm S}^{d}$ with the initial point $\gamma(0)=\bm{x}\in {\rm S}^{d}$, and the velocity $\dot{\gamma}(0)=\bm{z}\in T_{\bm{x}}{\rm S}^{d}/\{0\}$ is exactly the great circle,
\begin{align}
\label{Eq:geodesic_on_sphere}
\gamma(t)=\cos(t\|\bm{z}\|)\bm{x}+\sin(t\|\bm{z}\|)\dfrac{\bm{z}}{\|\bm{z}\|},
\end{align} 
Hence, the exponential map has the following expression
\begin{align}
\label{Eq:exponential_map}
\Exp_{\bm{x}}({\bm{z}})=\cos(\|\bm{z}\|)\bm{x}+\sin(\|\bm{z}\|)\dfrac{\bm{z}}{\|\bm{z}\|}.
\end{align}
This map can be understood as the rotation of the point $\bm{x}\in {\rm S}^{d}$ in the direction of $\bm{z}$ with the angle $\|\bm{z}\|$. From the exponential map \eqref{Eq:exponential_map}, we verify that the injectivity radius is $i(\bm{x})=\pi$, since given a point on the sphere $\bm{x}\in {\rm S}^{d}$, its antipodal point $-\bm{x}$ can be reached by choosing $\|\bm{z}\|=\pi$. 

Given a sample $\bm{\xi}\in \Xi$, we define the Lagrangian function 
\begin{align*}
\mathfrak{L}(\bm{x},\bm{\lambda};\bm{\xi})=\langle \bm{x},\bm{\xi}\rangle^{2}+\langle\bm{\lambda},\bm{x}\rangle-{\alpha\over 2}\|\bm{\lambda}\|.
\end{align*}
We compute the gradient $\grad_{\bm{x}}\mathfrak{L}(\bm{x},\bm{\lambda};\bm{\xi})\in T_{\bm{x}}{\rm S}^{d}$ as follows 
\begin{align*}
\grad_{\bm{x}}\mathfrak{L}(\bm{x},\bm{\lambda};\bm{\xi})=2\langle \bm{x},\bm{\xi}\rangle \bm{\xi}+\bm{\lambda}.
\end{align*}

The pseudo-code for the problem in eqs. \eqref{Eq:spike1}-\eqref{Eq:spike2} is described in Algorithm \ref{alg:3}. Corresponding to the estimate $\hat{\bm{\xi}}_{T}$, we define per coordinate constraint violation and the overlap performance metrics as follows
\begin{align*}
\Delta^{d}_{\text{CV}}(\hat{\bm{\xi}}_{T})&\topdoteq \|\Pi_{\real^{d}_{-}}(\hat{\bm{\xi}}_{T})\|/\sqrt{d},\\
\Delta_{\text{Overlap}}(\hat{\bm{\xi}}_{T})&\topdoteq |\langle \hat{\bm{\xi}}_{T},\bm{\xi}_{\ast} \rangle|.
\end{align*}

\paragraph{Semi-definite Programming (SDP)}

To evaluate the performance of the Riemannian primal-dual method, we also study the performance of the SDP method for the non-negative PCA problem given in eqs. \eqref{Eq:spike1}-\eqref{Eq:spike2},
\begin{subequations}
\begin{align}
\label{Eq:SDP}
&\min_{\bm{X}} \langle \bm{Y},\bm{X} \rangle\\ 
&\text{Tr}(\bm{X})=1,\\ \label{Eq:SDP1}
&\bm{X}\succeq 0,\quad \bm{X}\geq 0,
\end{align}
\end{subequations}
where $\bm{X}=[X_{ij}]\geq 0$ denotes the element wise inequality $X_{ij}\geq 0$, and $\bm{X}\succeq 0$ is the positive semidefinite constraint.Given the solution $\bm{X}_{*}\in \real^{T\times d}$ of the SDP problem eq. \eqref{Eq:SDP}-\eqref{Eq:SDP1}, we first compute the decomposition $\bm{X}_{*}=\sum_{i=1}^{\min\{T,d\}}\lambda_{i}\bm{u}_{i}\bm{u}_{i}^{T}$, where $\lambda_{1}\geq \lambda_{2}\cdots\geq \lambda_{d}$ are eigenvalues, and $\{\bm{u}_{i}\}_{i=1}^{d}$ are the corresponding eigenvectors. We then take the absolute value of the leading eigenvector $\hat{\bm{\xi}}_{T}=|\bm{u}_{1}|$ as the estimate of the principle component.

\subsubsection{Numerical Results}

We report our numerical results for the mixture model of \cite{montanari2016non}, where the samples are generated according to eq. \eqref{Eq:SpikedModel}. In particular, we let 
\begin{align*}
\xi_{\ast,i}=\begin{cases}
{1\over \sqrt{\delta d}}, &\quad i\in S\\
0, \quad &\hspace{3.4mm} i\not \in S,
\end{cases}
\end{align*}
where $S\subseteq \{1,2,\cdots,d\}$ and $|S|=\delta d$. We consider $T=d=2000$ and choose $\textsf{SNR}\in \{0.05,0.1,\cdots,2\}$ and sparsity level $\delta\in \{0.1,0.5,0.9\}$. For each $\textsf{SNR}$ value and the sparsity level $\delta$, we compute the average value of the overlap $\Delta_{\text{Overlap}}(\hat{\bm{\xi}}_{T})$ over 30 trials. 

\begin{figure}[t!]
  \centering
  \hspace{-4mm}
    \subfigure[]{
    \includegraphics[width=0.5\textwidth]{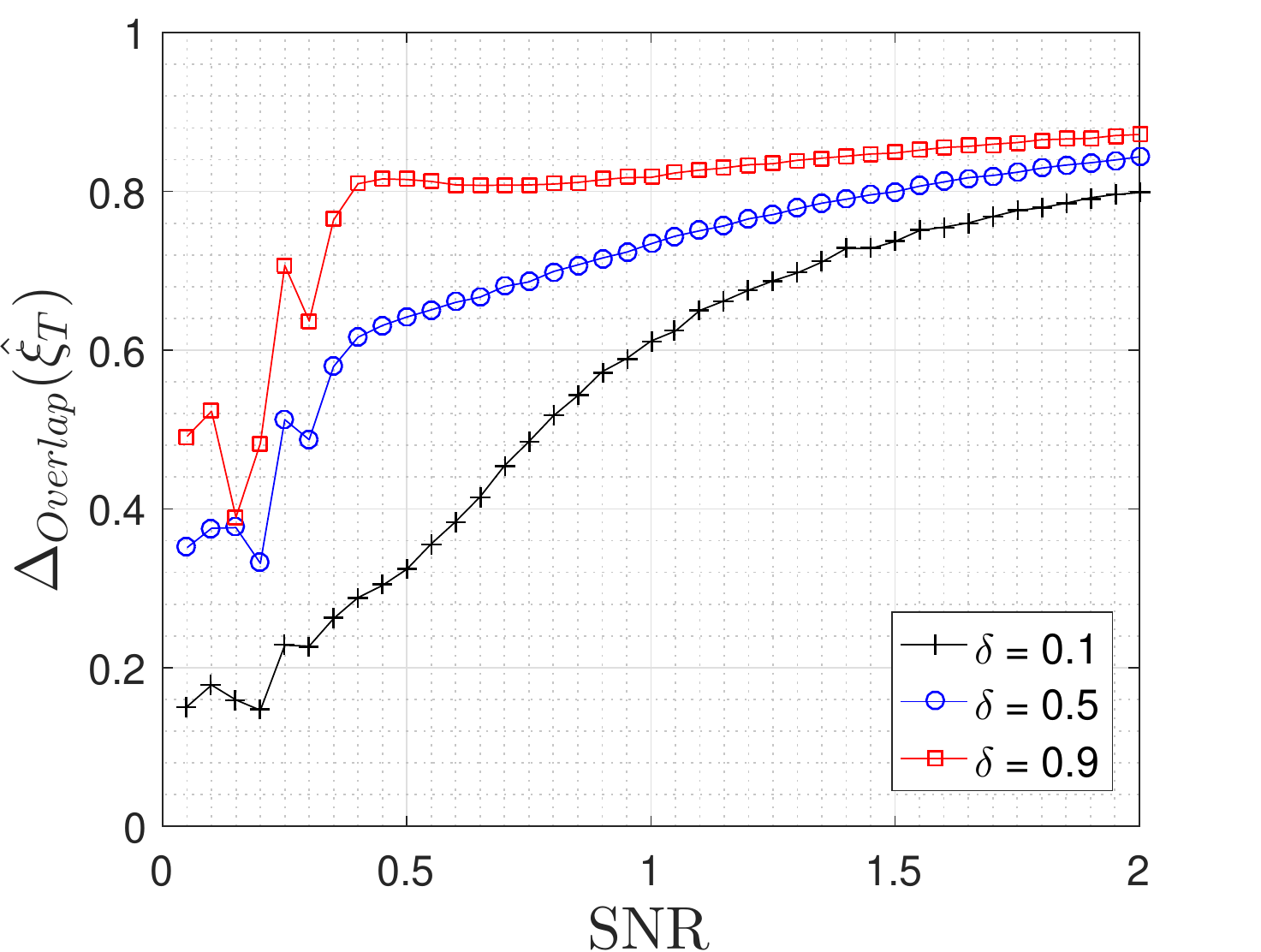}
     }\hspace{-4mm}
     \subfigure[]{
    \includegraphics[width=0.5\textwidth]{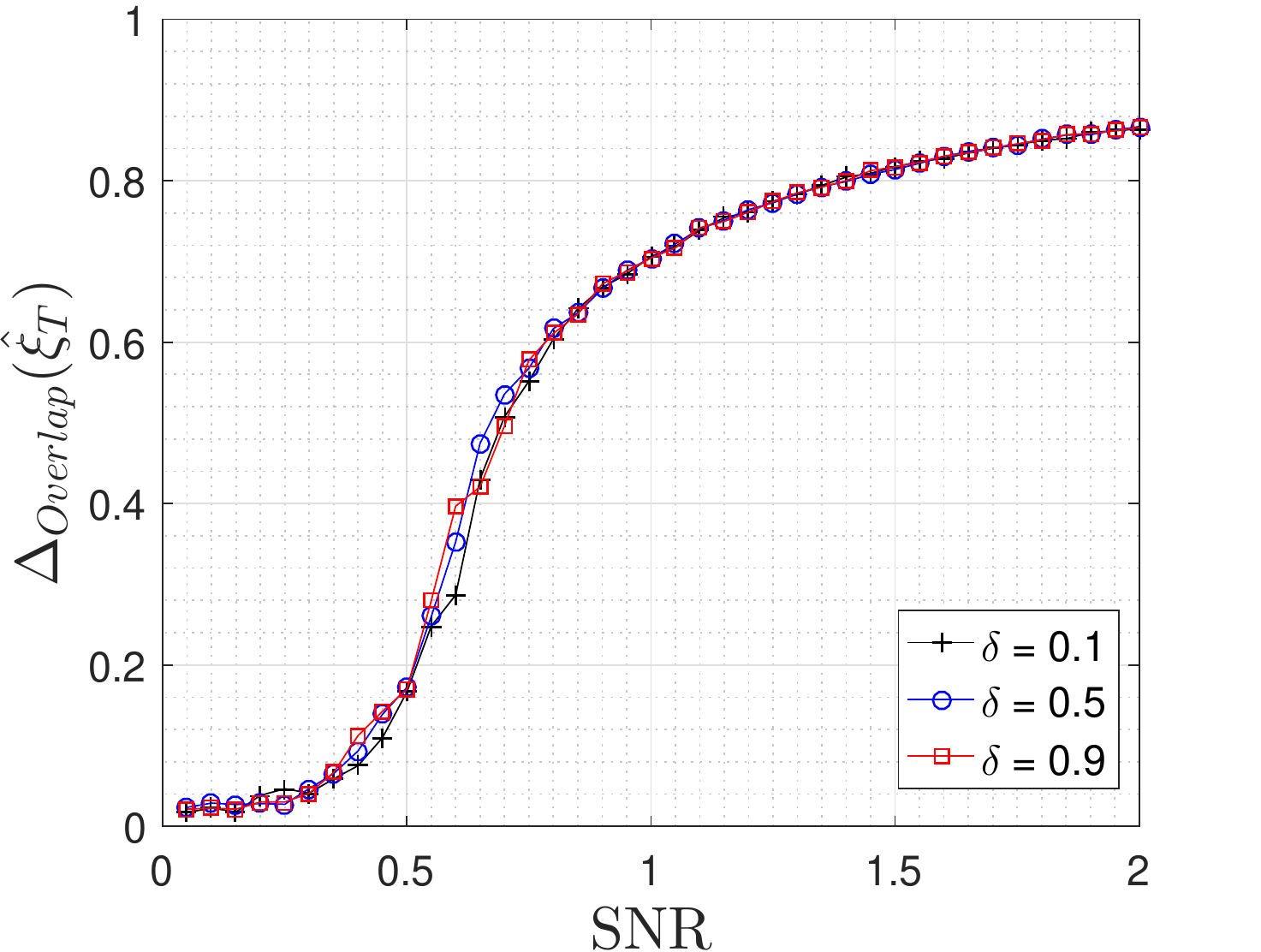}
     }
\caption{Comparison between the average overlap $\Delta_{Overlap}(\hat{\bm{\xi}}_{T})$ of Algorithm \ref{alg:3} for online non-negative PCA, and the spectral method for the vanilla PCA over 30 trials and for the sample size of $T=2000$ and different sparsity levels $\delta$. (a) Riemannian primal-dual method, (b) SVD algorithm.}
\label{Fig:Overlap}

\vspace*{5mm}
\minipage{0.32\textwidth}
  \includegraphics[trim={.6cm 1cm .6cm 1cm},width=\linewidth]{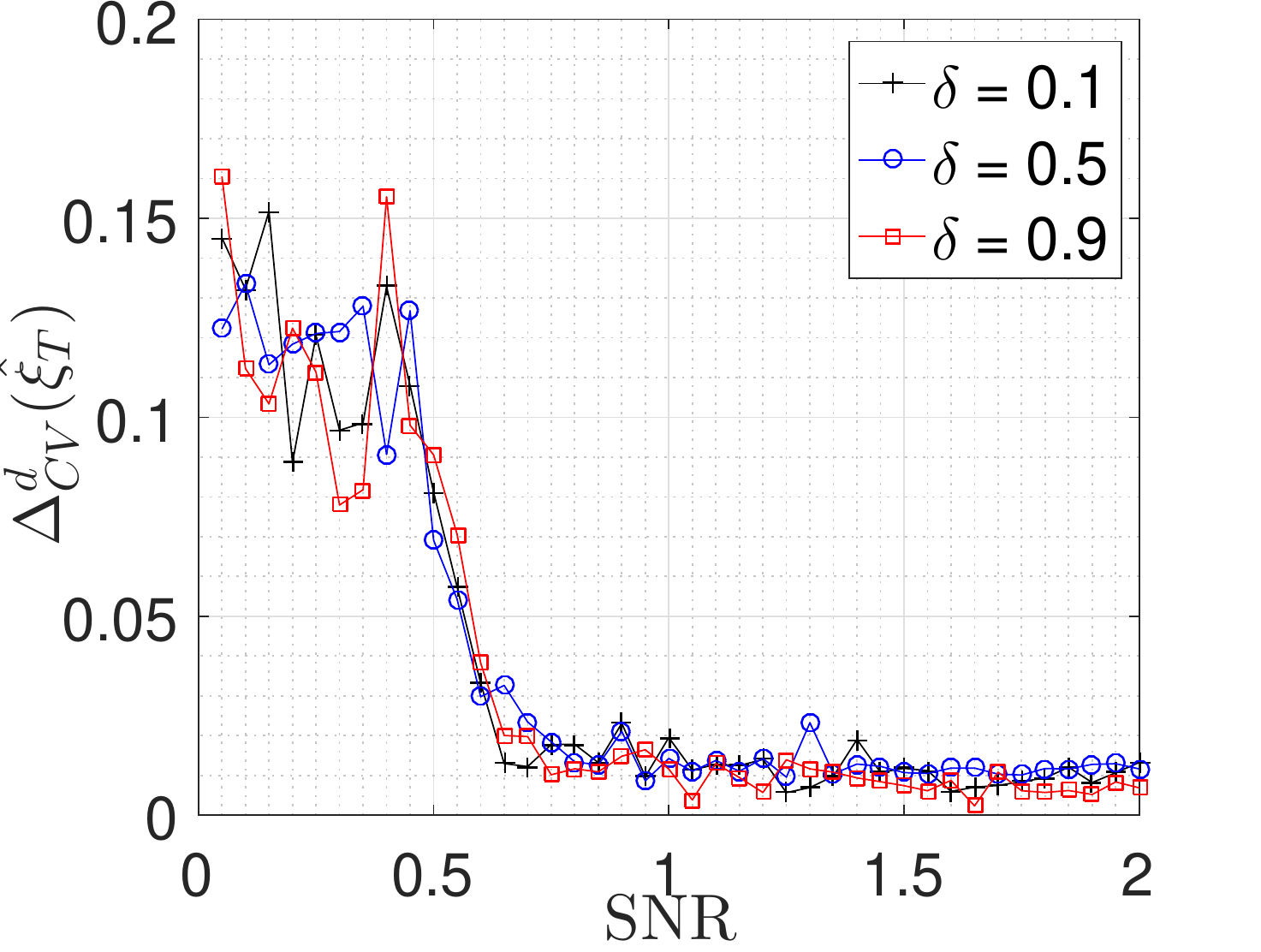}
\endminipage\hfill
\minipage{0.32\textwidth}
  \includegraphics[trim={.6cm 1cm .6cm 1cm},width=\linewidth]{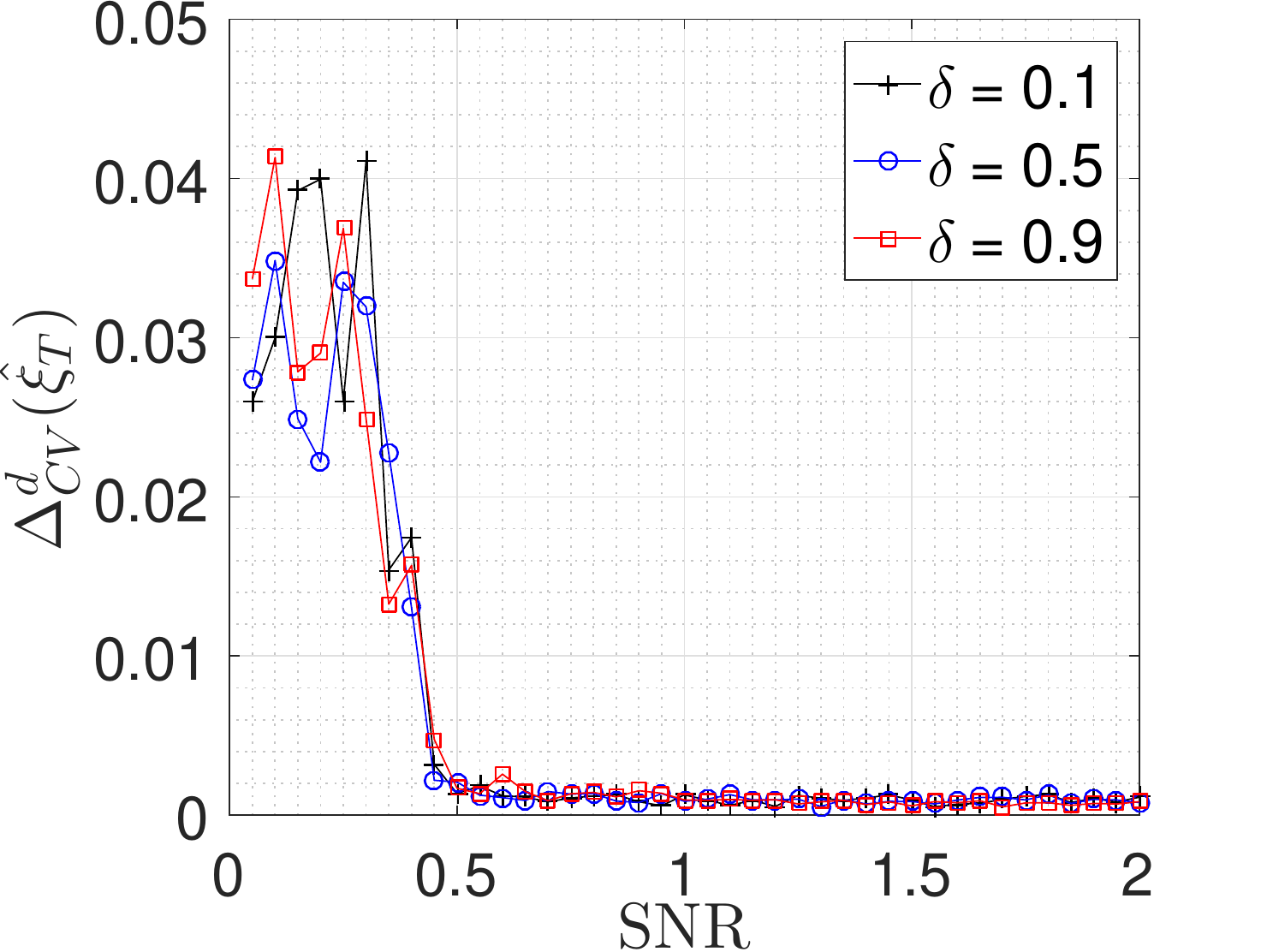}
\endminipage\hfill
\minipage{0.32\textwidth}%
  \includegraphics[trim={.6cm 1cm .6cm 1cm},width=\linewidth]{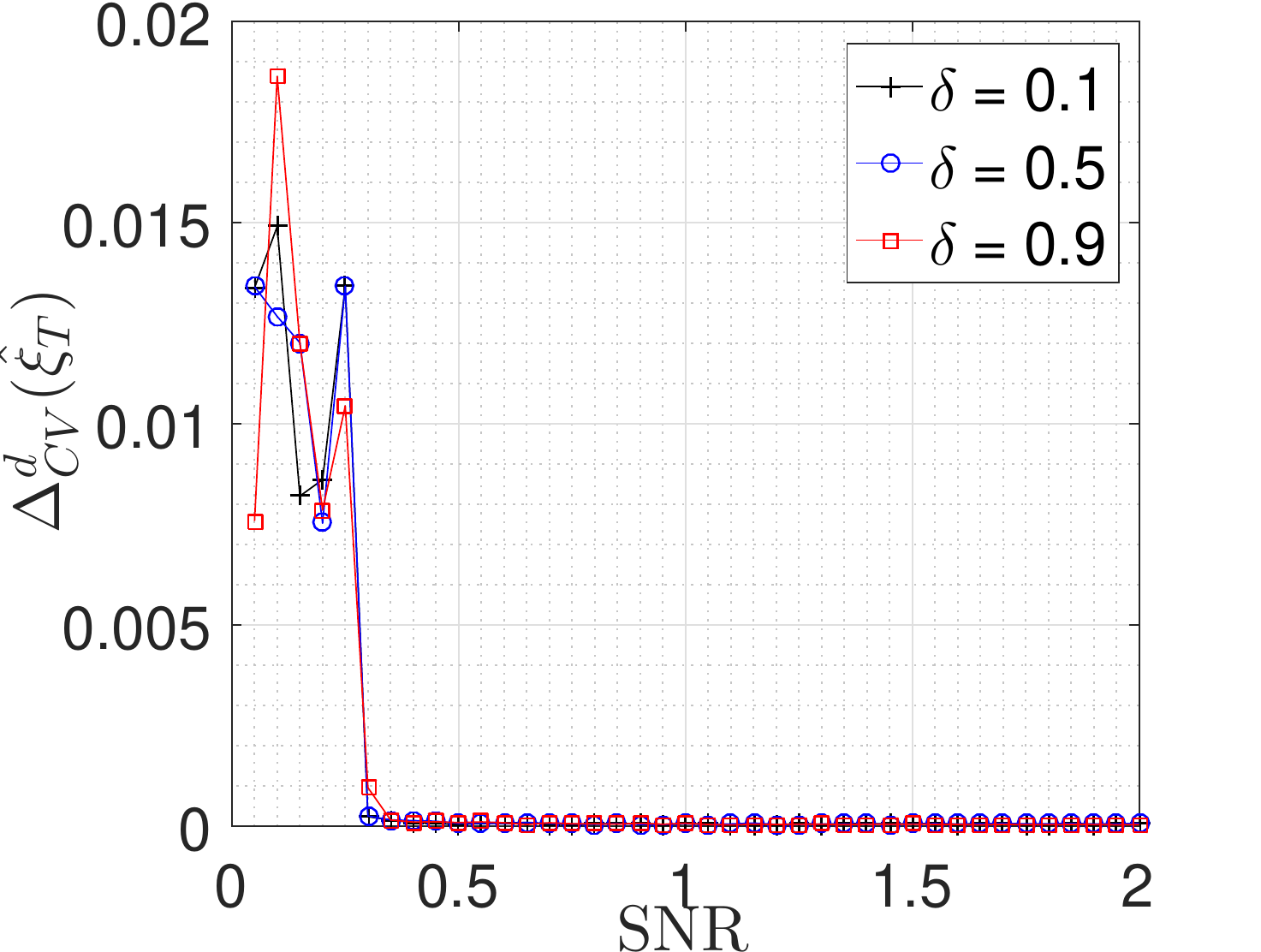}
\endminipage
\caption{Average constraint violation $\Delta^{d}_{\text{CV}}(\hat{\bm{\xi}}_{T})$ of the Riemannian primal-dual solution for different sample sizes over 30 trials. (a) $T=20$, (b) $T=200$, (c) $T=2000$.}
\label{Fig:ConstraintViolation}
\end{figure}

Figure \ref{Fig:Overlap} shows the average overlap $\Delta_{\text{Overlap}}(\bm{\xi}_{T})$ of Algorithm \ref{alg:3} for different $\mathsf{SNR}$ values. It further shows the performance of the spectral method for the vanilla PCA in eqs. \eqref{Eq:population}-\eqref{Eq:population1}, where we first compute the singular value decomposition (SVD) of $\bm{Y}\bm{Y}^{T}/T$ and then choose the leading eigenvector $\hat{\bm{\xi}}(\bm{Y})$ corresponding to the largest eigenvalue. We observe from Figure \ref{Fig:Overlap} that the sparsity level $\delta$ of the principle component $\bm{\xi}_{\ast}$ does not affect the SVD decomposition algorithm. In particular, we obtain the same overlap values from SVD regardless of the sparsity level $\delta$. In contrast, the performance of Algorithm \ref{alg:3} significantly improves in low $\mathsf{SNR}$ regimes as the sparsity level $\delta$ decreases.    

Figure \ref{Fig:ConstraintViolation} shows the average constraint violation $\Delta_{\text{CV}}^{d}(\hat{\bm{\xi}}_{T})$ of the Riemannian primal-dual solution for different sample sizes $T\in \{20,200,2000\}$ over 30 trials. Evidently, the constraint violation is only significant in the low $\mathsf{SNR}$ regimes, where the samples $\bm{\xi}_{0},\cdots,\bm{\xi}_{T-1}$ are nearly random Gaussian. However, by increasing the sample size $T$, the amount of constraint violation of the Riemannian primal-dual solution diminishes. Notice that the amount of constraint violation also depends on the tunable parameter $\alpha$ of the regularizer, which is fixed to $\alpha=0.1$ in all the experiments.

We now turn to the SDP method, and compare its average overlap performance $\Delta_{\text{overlap}}(\hat{\bm{\xi}}_{T})$ over 30 trials with Algorithm \ref{alg:3} in Table \ref{Table:3}. Further, we have computed the algorithm run time in each case and reported  the amount of time (sec) that each trial takes.\footnote{We measured the elapsed time of each algorithm using the $\tt{tic-toc}$ command in MATLAB.} We use \texttt{CVX} \cite{grant2008cvx} to solve the SDP problem in eqs. \eqref{Eq:SDP}-\eqref{Eq:SDP1}, where we choose the best precision settings. We observe that although SDP is a polynomial time algorithm, it scales very poorly in practice and can not be applied to problems with a large sample size $T$. In fact, for the sample size $T=1000$, we were unable to obtain the principle component from SDP due to the physical memory limitations of our setup.\footnote{The experiments are performed in MATLAB on a Windows 7 machine equipped with the dual core 2.10 GHz Intel Xeon processor, and 32GB of RAM.}  In contrast, the Riemannian primal-dual method scales gracefully with the sample size $T$.   

\begin{remark}
We note that the SDP method achieves substantially higher overlaps than Algorithm \ref{alg:3} when the sample size is small in Table \ref{Table:3}. This is due to the fact that the SDP algorithm is an offline algorithm in which all the samples $\bm{\xi}_{0},\cdots,\bm{\xi}_{T-1}$ are known. In comparison, Algorithm \ref{alg:3} is an online algorithm and thus it naturally requires more samples to converge to a good solution.
\end{remark}

\begin{table}[t!]
\label{Table:3}
\caption{Average performance comparison of the online Riemannian primal-dual method and offline SDP for non-negative PCA problem over 30 trials. The reported numbers for the run time is per trial and $\delta=0.9$ and $\textsf{SNR}=1$.} 
\centering 
  \begin{tabular}{l|l l|l l}
    \hline\hline
    \multirow{2}{*}{} &
      \multicolumn{2}{c}{SDP} &
      \multicolumn{2}{c}{Primal-Dual} 
 \\
    Sample Size& $\Delta_{\text{Overlap}}$ & Run-Time(sec)  &$\Delta_{\text{Overlap}}$ &Run-Time(sec) \\
    \hline
    $T=10$ & 0.941 & 0.52 & 0.767 & 0.0029  \\
    \hline
    $T=50$ & 0.824 & 8.12 & 0.807 & 0.0032  \\
    \hline
    $T=100$ &0.799 & 94.49 & 0.809 & 0.0050  \\
    \hline
    $T=200$ &0.841 & 1304.15 & 0.816 & 0.0103  \\
    \hline    
    $T=1000$& NC &  NC    &0.821 & 0.0849  \\    
    \hline  
  \end{tabular}
\end{table}

\subsection{Anchored Synchronization on $\rm{SO}(3)$}
\label{sec:synch}
In general, the synchronization problem is concerned with estimating a set of elements $g_{1},\cdots,g_{n}$ of a group $G$ from noisy observation of their relative positions $y_{ij}=w_{ij}\ast g_{i}\ast g_{j}^{-1}$, where $w_{ij}$ is the noise which is also an element of the group $G$. The synchronization problems on ${\rm Z}_{2}$ and ${\rm U}(1)$ groups have been studied in \cite{javanmard2016phase}. The case of synchronization on ${\rm SO}(3)$ has also received a great deal of attention \cite{boumal2014cramer,singer2011three}. 

The rotation group ${\rm SO(3)}$ is the natural state-space for the orientation of rigid bodies, and thus has many applications in robotics and computer vision for camera registration, see, \textit{e.g.}, \cite{horn1990relative},\cite{gurwitz1989globally}. Therefore, in the sequel we focus on the synchronization problem for recovering the elements of the rotation group $\bm{R}_{1},\bm{R}_{2},\cdots,\bm{R}_{n}\in \rm{SO}(3)$.

\subsubsection{Geometry of the Rotation Group $\rm{SO}(3)$}

The rotation group ${\rm SO(3)}$ can be defined as follows
\begin{align*}
{\rm SO(3)}=\{\bm{R}\in {\rm GL(3,\real)}: \bm{R}^{T}\ast \bm{R}=\bm{I}, \text{det}(\bm{R})=1\},
\end{align*}
where ${\rm GL(3,\real)}$ is the set of $3\times 3$ non-singular matrices, and the group operation $\ast$ is simply the matrix multiplication. Moreover, the metric is defined by $g(\bm{R}_{1},\bm{R}_{2}) = \langle \bm{R}_{1},\bm{R}_{2}\rangle= \text{Tr}(\bm{R}_{1}^{T}\bm{R}_{2})$ for any $\bm{R}_{1},\bm{R}_{2}\in {\rm SO}(3)$, which induces the Frobenius norm denoted by $\|\cdot\|_{F}$.

The Lie algebra $\mathfrak{so}(3)$ associated with the Lie group ${\rm SO(3)}$ is the set of all skew symmetric matrices, \textit{i.e.},
\begin{align*}
\mathfrak{so}(3)=T_{\bm{I}}{\rm SO(3)}=\{\bm{X}\in \real^{3\times 3}:\bm{X}=-\bm{X}^{T}\},
\end{align*}
where the Lie bracket is the commutator $[\bm{X},\bm{Y}]=\bm{X}\bm{Y}-\bm{Y}\bm{X}$. The tangent plan at any other point $\bm{R}\in {\rm SO}(3)$ can be obtained simply by $T_{\bm{R}}{\rm SO(3)}=\bm{R}\mathfrak{so}(3)=\{\bm{R}\bm{X}:\bm{X}\in \mathfrak{so}(3)\}$. The generators of the algebra $\bm{e}_{1},\bm{e}_{2},\bm{e}_{3}$ correspond to the derivative of rotation around each of the standard axis, namely
\begin{align*}
\bm{e}_{1}=\left( {\begin{array}{ccc}0 & 0 & 0 \\ 0 & 0 &-1  \\ 0 &1 &0  \end{array} }\right),  \quad \bm{e}_{2}=\left( {\begin{array}{ccc}0 & 0 & 1 \\ 0 & 0 &0  \\ -1 &0 &0  \end{array} }\right), \quad \bm{e}_{3}=\left( {\begin{array}{ccc}0 & -1 & 0 \\ 1 & 0 &0  \\ 0 &0 &0  \end{array} }\right).
\end{align*}
The commutator of the basis is given by $[\bm{e}_{i},\bm{e}_{j}]=\epsilon_{ijk}\bm{e}_{k}$, for $i,j,k\in \{1,2,3\}$ where $\epsilon_{ijk}$ is the Levi-Civita symbol. 
 Due to the compactness\footnote{Recall that a group is compact, if its point-set topology is compact.}, the rotation group $\rm{SO}(3)$ admits a bi-invariant metric with a non-negative sectional curvatures $K\geq 0$, cf. \cite[Corollary 1.3]{milnor1976curvatures}. In general, consider the $n$-dimensional Lie group $G$ with its associated Lie algebra $\mathfrak{g}$. Given a generator $\bm{e}_{1},\bm{e}_{2},\cdots,\bm{e}_{n}$ of the algebra $\mathfrak{g}$, the \textit{structure constants} $\alpha_{ijk}$ are defined as follows
\begin{align*}
[\bm{e}_{i},\bm{e}_{j}]=\sum_{k=1}^{n}\alpha_{ijk}\bm{e}_{k},  \quad \ \alpha_{ijk}\topdoteq\langle [\bm{e}_{i},\bm{e}_{j}],\bm{e}_{k}\rangle.
\end{align*}
Note that the structure constants depend on the choice of the metric. For two orthonormal vectors $\bm{X},\bm{Y}\in \mathfrak{g}$, a generator can always be constructed such that $\bm{e}_{1}=\bm{X}$ and $\bm{e}_{2}=\bm{Y}$. As shown by Milnor \cite{milnor1976curvatures}, the sectional curvature for a bi-invariant metric is given by 
\begin{align}
\label{Eq:Lie_group_Milnor}
K(\bm{X},\bm{Y})=\dfrac{1}{4}\|[\bm{X},\bm{Y}]\|^{2} ={1\over 4}\sum_{j=1}^{n}(\alpha_{2j1})^{2}\geq 0.
\end{align}
In the special case of the rotation group $\rm{SO}(3)$, the identity of eq. \eqref{Eq:Lie_group_Milnor} gives us $K(\bm{X},\bm{Y})=1/4$ for all $\bm{X},\bm{Y}\in \mathfrak{so}(3)$. Hence, the condition of Assumption \ref{Assumption:4} is satisfied. 

The exponential map $\Exp:\mathfrak{so}(3)\rightarrow {\rm SO}(3)$ for the matrix Lie groups is simply the matrix exponential, \textit{i.e.}, $\Exp(\bm{X})=\sum_{k=0}^{\infty}\bm{X}^{k}/k!$ for $\bm{X}\in \mathfrak{so}(3)$. However, dealing with the infinite series of the matrix powers is challenging for numerical reasons. We thus alternatively use \textit{Rodrigues' rotation formula} \cite[p. 28]{murray1994mathematical},
\begin{align}
\label{Eq:exponential_map_SO}
\Exp(\bm{X})=\bm{I}+\dfrac{\sin(\|\bm{X}\|_{F})}{\|\bm{X}\|_{F}}\bm{X}+\dfrac{1-\cos(\|\bm{X}\|_{F})}{\|\bm{X}\|_{F}^{2}}\bm{X}^{2}, \quad \bm{X}\in \mathfrak{so}(3), \bm{X}\not= \bm{0}.
\end{align}
The logarithmic map $\log_{\bm{I}}:{\rm SO}(3)\rightarrow \mathfrak{so}(3)$ is also given by
\begin{align}
\label{Eq:Logarithmic_map_SO}
\log_{\bm{I}}(\bm{R})=\left\{
    \begin{array}{ll}
        0 & \phi=0 \\
        \dfrac{\phi}{2\sin \phi}(\bm{R}-\bm{R}^{T}) & \phi\in (0,\pi)
    \end{array}
\right.
\end{align}
where $\phi$ satisfies $\cos(\phi)=(\text{Tr}(\bm{R})-1)/2$ and $\|\log(\bm{R})\|_{F}=\phi$. We note that when $\text{Tr}(\bm{R})=-1$ two solutions for the $\log$ function exist, namely, if $\bm{u}=(u_{1},u_{2},u_{3})$ is the unit eigenvector of $\bm{R}$ associated with the eigenvalue $1$, then $\log(\bm{R})=\pm \pi (u_{1}\bm{e}_{1}+u_{2}\bm{e}_{2}+u_{3}\bm{e}_{3})$, see \cite{park1997smooth}. Hence, the injectivity radius is $\text{inj}({\rm SO(3)})=\pi$.

From the logarithmic map, the geodesic distance between the two rotation matrices $\bm{R}_{1},\bm{R}_{2}\in {\rm SO}(3)$ can be determined as follows
\begin{align}
\label{Eq:geodesic_distance}
d(\bm{R}_{1},\bm{R}_{2})&=\|\log_{\bm{I}}(\bm{R}^{T}_{1}\bm{R}_{2})\|_{F}=\|\log_{\bm{I}}(\bm{R}^{T}_{2}\bm{R}_{1})\|_{F}\\ \label{Eq:geodesic_distance_1}
&=\arccos\left({(\text{Tr}(\bm{R}_{1}^{T}\bm{R}_{2})-1)\over 2}\right)=\arccos\left({(\langle \bm{R}_{1},\bm{R}_{2}\rangle -1) \over 2}  \right).
\end{align} 

\begin{remark}
Based on the analysis of Section \ref{subsection:produtct_manifold}, the geodesics and the exponential map of the product manifold $\mathcal{M}={\rm SO}(3)\otimes \cdots \otimes {\rm SO(3)}$ can also be defined coordinate wise, under the extension defined in eq. \eqref{Eq:construction}. 
\end{remark}

\subsubsection{Measurement Model}
Consider the undirected, connected graph $\mathcal{G}=(V,E)$ with vertices $V\topdoteq \{1,2,\cdots,n\}$ and edges $E\subseteq V\times V$. The set of measurements of the relative rotations on the undirected graph $\mathcal{G}$ can be described as follows,
\begin{align}
\label{Eq:noise_model}
\bm{Y}_{ij}=\bm{W}_{ij}\bm{R}^{0}_{i}(\bm{R}^{0}_{j})^{T}, \quad (i,j)\in E,
\end{align} 
where here $\bm{W}_{ij}\in {\rm SO(3)}$ are i.i.d. multiplicative noise that are sampled from a distribution with the probability density function $g:\rm{SO(3)}\rightarrow \real_{+}$ with respect to the Haar measure $\mu$ (see Appendix \ref{sec:Lie Groups}). In particular,
\begin{align*}
\int_{\rm{SO}(3)}g(\bm{R})\mu(d\bm{R})=1.
\end{align*}
We assume that the multiplicative noise is unbiased, \textit{i.e.}, $\expect_{g}[\bm{W}_{ij}]=\bm{I}$. Notice that the measurement $(i,j)$ and $(j,i)$ on the undirected graph are the same, and hence we have $\bm{Y}_{ij}^{T}=(\bm{W}_{ij}\bm{R}^{0}_{i}(\bm{R}^{0}_{j})^{T})^{T}=\bm{W}_{ji}\bm{R}^{0}_{j}(\bm{R}^{0}_{i})^{T}=\bm{Y}_{ji}$. Thus necessarily, $\bm{W}_{ji}=\bm{R}^{0}_{j}(\bm{R}^{0}_{i})^{T}\bm{W}_{ij}^{T}\bm{R}_{i}(\bm{R}_{j}^{0})^{T}$. The symmetry on the observations also enforces certain restrictions on an admissible density function $g$; see \cite{boumal2014cramer}.

\subsection{Algorithm}
To apply the Riemannian primal-dual framework to the synchronization problem, we begin with the problem of minimizing the mean squared error (MMSE) on the product manifold $\mathcal{M}\topdoteq {\rm SO(3)}\otimes {\rm SO(3)}\otimes \cdots \otimes {\rm SO(3)}$,
\begin{align}
\label{Eq:min_squared_error}
\text{minimize}_{{\bm{R}}\in \mathcal{M}} \expect_{\bm{Y}}\bigg[\sum_{(i,j)\in E}d^{2}(\bm{I},\bm{Y}_{ij}\bm{R}_{j}\bm{R}_{i}^{T})\bigg],
\end{align}
where $\bm{R}\topdoteq (\bm{R}_{1},\cdots,\bm{R}_{n})$, $\bm{Y}\topdoteq (\bm{Y}_{ij})_{(i,j)\in E}$, and $\bm{W}_{ij}=\bm{Y}_{ij}\bm{R}_{j}\bm{R}_{i}^{T}$. Notice that the mean squared error (MSE) in eq. \eqref{Eq:min_squared_error} is invariant under the right transformation $\bm{R}_{i}\mapsto \bm{R}_{i}\bm{A}$ for $\bm{A}\in {\rm SO(3)}$ and all $i=1,2,\cdots,n$, \textit{i.e.}, $d^{2}(\bm{I},\bm{Y}_{ij}\bm{R}_{j}\bm{R}_{i}^{T})=d^{2}(\bm{I},\bm{Y}_{ij}\bm{R}_{j}\bm{A}\bm{A}^{T}\bm{R}_{i}^{T})$. A symmetry breaking strategy is to consider some reference nodes (anchors), that is $\bm{R}_{i}=\bm{R}^{0}_{i}$ for $i\in \widetilde{V}\subset V$, see \cite{boumal2014cramer}. Employing anchors to eliminate redundancies and pin down the underlying state configuration is also a common approach in the localization of sensor networks, where anchors consist of those sensors that are able to estimate their locations based on a standalone localization mechanism such as GPS \cite{romer2005time}.

Including the anchors in the MMSE problem \eqref{Eq:min_squared_error} yields
\begin{subequations}
\begin{align}
\label{Eq:Max_Likelihood_1}
&\text{minimize}_{\bm{R}\in \mathcal{M}} {1\over 2}\expect_{\bm{Y}}\bigg[\sum_{(i,j)\in E}\|\bm{I}-\bm{Y}_{ij}\bm{R}_{j}\bm{R}_{i}^{T} \|_{F}^{2}\bigg],\\ \label{Eq:Max_Likelihood_2}
&\text{subject to}: \|\bm{I}-\bm{R}^{T}_{i}\bm{R}^{0}_{i}\|_{F}^{2}\leq 0, \quad i\in \widetilde{V}\subseteq V.
\end{align}
\end{subequations}
Two remarks on the reformulation in eqs. \eqref{Eq:Max_Likelihood_1}-\eqref{Eq:Max_Likelihood_2} are in order. 

First, note that the geodesic distance in eq. \eqref{Eq:min_squared_error} is replaced with the Frebenous distance $d_{F}(\bm{R}_{1},\bm{R}_{1})\topdoteq\|\bm{I}-\bm{R}_{1}^{T}\bm{R}_{2}\|_{F},\bm{R}_{1},\bm{R}_{2}\in {\rm SO}(3)$. This relaxation is permissible since the geodesic distance and the Frebenous distance are boundedly equivalent with respect to topology of ${\rm SO}(3)$ \cite{huynh2009metrics}, \textit{i.e.}, there are positive constants $a,b>0$ such that
\begin{align*}
a\cdot d_{F}(\bm{R}_{1},\bm{R}_{2})\leq d(\bm{R}_{1},\bm{R}_{2})\leq b\cdot d_{F}(\bm{R}_{1},\bm{R}_{2}).
\end{align*}

Second, notice that the equality constraints $\bm{R}_{i}=\bm{R}^{0}_{i}$ are written in the inequality forms. An advantage of the current formulation  is that the inequality constraints can be easily modified to impose soft constraints $\|\bm{I}-\bm{R}_{i}^{T}\bm{R}^{0}_{i}\|_{F}^{2}\leq \epsilon$ on the rotation states for some $\epsilon>0$. This is particularly useful when side information $\bm{R}^{0}_{i},i\in \widetilde{V}$ are themselves some estimates of the true states.

We now apply the Riemannian primal-dual method to the problem in eqs. \eqref{Eq:Max_Likelihood_1}-\eqref{Eq:Max_Likelihood_2}. We define the Lagrangian function,
\begin{align*}
&L:(\real^{3\times 3})^{n}\times \real_{+}^{n} \rightarrow \real_{+} \\
&(\bm{R},\bm{\lambda};\bm{Y})\mapsto {1\over 2}\sum_{(i,j)\in E}\|\bm{I}-\bm{Y}_{ij}\bm{R}_{j}\bm{R}_{i}^{T}\|_{F}^{2}+{1\over 2}\sum_{i\in \widetilde{V}}\lambda_{i} \|\bm{I}-\bm{R}_{i}\bm{R}^{0}_{i}\|^{2}_{F}-\dfrac{\alpha}{2}\|\bm{\lambda}\|_{2}^{2}.
\end{align*}
The elements of the gradient vector $\grad L(\bm{R},\bm{\lambda};\bm{Y})=(\bm{Z}_{1},\cdots,\bm{Z}_{n})$ are given by
\begin{align}
\nonumber
\bm{Z}_{i}=&-\sum_{j\in N(i)}(\bm{I}-\bm{R}_{i}\bm{R}^{T}_{j}\bm{Y}_{ji})\bm{Y}_{ij}\bm{R}_{j}- \sum_{j\in N(i)}\bm{Y}_{ij}(\bm{I}-\bm{Y}_{ji}\bm{R}_{i}\bm{R}_{j}^{T})\bm{R}_{j} \\ \label{Eq:gradient_of_SO(3)} &- \lambda_{i}\bm{R}_{i}^{0}(\bm{I}-(\bm{R}_{i}^{0})^{T}\bm{R}_{i})\mathds{1}_{i\in \tilde{V}},
\end{align}
where $\mathds{1}_{\{\cdot\}}$ is the indicator function, and we recall $N(i)\topdoteq \{j\in V:(i,j)\in E\}$.

Given the restriction $\mathfrak{L}(\bm{R},\bm{\lambda};\bm{Y})=L(\bm{R},\bm{\lambda};\bm{Y})\vert_{\mathcal{M}}$, we obtain the $i$-th element of the gradient $\grad \mathfrak{L}(\bm{R},\bm{\lambda};\bm{Y})$ by projecting $\bm{Z}_{i}$ onto $T_{\bm{R}_{i}}{\rm SO(3)}=\bm{R}_{i}\mathfrak{so}(3)$,
\begin{align}
\nonumber
(\grad \mathfrak{L}(\bm{R},\bm{\lambda};\bm{Y}))_{i}&=\Pi_{\bm{R}_{i}\mathfrak{so}(3)}(\bm{Z}_{i}) =\bm{R}_{i}\text{skew}(\bm{R}_{i}^{T}\bm{Z}_{i})\\  
&={1\over 2} \bm{R}_{i}(\bm{R}_{i}^{T}\bm{Z}_{i}-\bm{Z}_{i}^{T}\bm{R}_{i}), 
\end{align}
where $\text{skew}(\bm{X})\topdoteq(\bm{X}-\bm{X}^{T})/2$ is the skew symmetric matrix based on $\bm{X}$. Notice that $(\grad \mathfrak{L}(\bm{R},\bm{\lambda};\bm{Y}))_{i}\in \bm{R}_{i}\mathfrak{so}(3)$. Since the exponential map is a map from $T_{\bm{I}}{\rm SO(3)}=\mathfrak{so}(3)$ to ${\rm SO(3)}$, we subsequently use the left translation to compute an element of $\bm{s}_{i}\in \mathfrak{so}(3)$ as follows
\begin{align}
\label{Eq:SO(3)_Compute}
\bm{s}_{i}=\bm{R}_{i}^{-1}(\grad \mathfrak{L}(\bm{R},\bm{\lambda};\bm{Y}))_{i}={1\over 2}(\bm{R}_{i}^{T}\bm{Z}_{i}-\bm{Z}_{i}^{T}\bm{R}_{i}), 
\end{align}
which is clearly a skew symmetric matrix. 

We now are in position to describe the structure of the Riemannian primal-dual algorithm for solving \eqref{Eq:Max_Likelihood_1}-\eqref{Eq:Max_Likelihood_2}, see the pseudocode of Algorithm \ref{alg:4}.

\begin{algorithm}[t!]
\caption{Riemannian Primal-Dual Method for Synchronization on ${\rm SO}(3)$}
\label{alg:4}
\begin{algorithmic}[1]
\State{\textbf{require}: A connected synchronization graph $E\subseteq V\times V$ on the vertices $V$.}
\State{\textbf{initialize}: $\bm{R}^{0}_{i}\in {\rm SO}(3)$, $\lambda^{0}_{i}=0$ for all $i\in V$.}
\For{$t=0,1,\cdots,T-1$}
\State{Let $\bm{W}_{ij}^{t}\sim g(\bm{X})$ be i.i.d. noise terms for all edges $(i,j)\in E$.}
\State{Observe the noisy samples $\bm{Y}_{ij}^{t}=\bm{W}_{ij}^{t}\bm{R}^{0}_{i}(\bm{R}^{0}_{j})^{T}$ of all edges $(i,j)\in E$.}
\State{Compute $\grad \mathfrak{L}(\bm{R}^{t},\bm{\lambda}^{t})=(\bm{X}_{1}^{t},\bm{X}_{2}^{t},\cdots,\bm{X}_{n}^{t})$ from eqs. \eqref{Eq:gradient_of_SO(3)}-\eqref{Eq:SO(3)_Compute}.}
\State{$\bm{R}^{t+1}_{i}\leftarrow \bm{R}^{t}_{i}\left(\bm{I}-\dfrac{\sin(\eta_{t}\|\bm{X}_{i}^{t}\|_{F})}{\|\bm{X}_{i}^{t}\|_{F}}\bm{X}_{i}^{t}+\dfrac{1-\cos(\eta_{t}\|\bm{X}_{i}^{t}\|_{F})}{\|\bm{X}_{i}^{t}\|_{F}^{2}}(\bm{X}_{i}^{t})^{2}\right)$.}
\State{$\lambda^{t+1}_{i}\leftarrow \max\{0,(1-\alpha)\lambda^{t}_{i}+\eta_{t}\|\bm{R}_{i}-\bm{R}_{i}^{0}\|^{2}_{F}\}$. }
\EndFor

\hspace{-7mm}\Return $(\bm{R}_{1}^{T},\bm{R}_{2}^{T},\cdots,\bm{R}_{n}^{T})$.
\end{algorithmic}
\end{algorithm}

\subsubsection{Numerical Results}

We now evaluate the performance of Algorithm \ref{alg:4} on the connected Erd\"{o}s-R\'{e}yni random graph $G(n,p)$ we described in Section \ref{sec:MAX-CUT}. Further, we use the isotropic Langevin distribution as an analogue of the Gaussian distribution for the rotation group in the following simulations \cite{boumal2014cramer}, 
\begin{align}
\label{Eq:F_function}
g(\bm{R})={1\over Z}\cdot \exp(\beta \langle \bm{R},\bm{N}\rangle), \quad \bm{R}\in {\rm SO}(3),
\end{align}
where $\bm{N}\in {\rm SO(3)}$ denotes the mean value, $\beta>0$ is the concentration parameter, and 
$Z$ is the normalization factor,
\begin{align}
Z&=\int_{{\rm SO}(3)}\exp(\beta \langle \bm{R},\bm{N}\rangle)d\mu(\bm{X})\\
&\stackrel{{\rm (a)}}{=}\int_{{\rm SO}(3)}\exp(\beta \langle \bm{R},\bm{I}\rangle)d\mu(\bm{R})\\
&\stackrel{{\rm (b)}}{=}\exp(\beta)(I_{0}(2\beta)-I_{1}(2\beta)),
\end{align} 
where ${\rm (a)}$ is due to the bi-invariance of the Haar measure (cf. Appendix \ref{sec:Lie Groups}), and in ${\rm (b)}$, $I_{0}(2\beta),I_{1}(2\beta)$ are the Bessel functions of the first kind.

\begin{figure}[t!]
  \centering
    \subfigure[]{
    \includegraphics[width=0.45\textwidth]{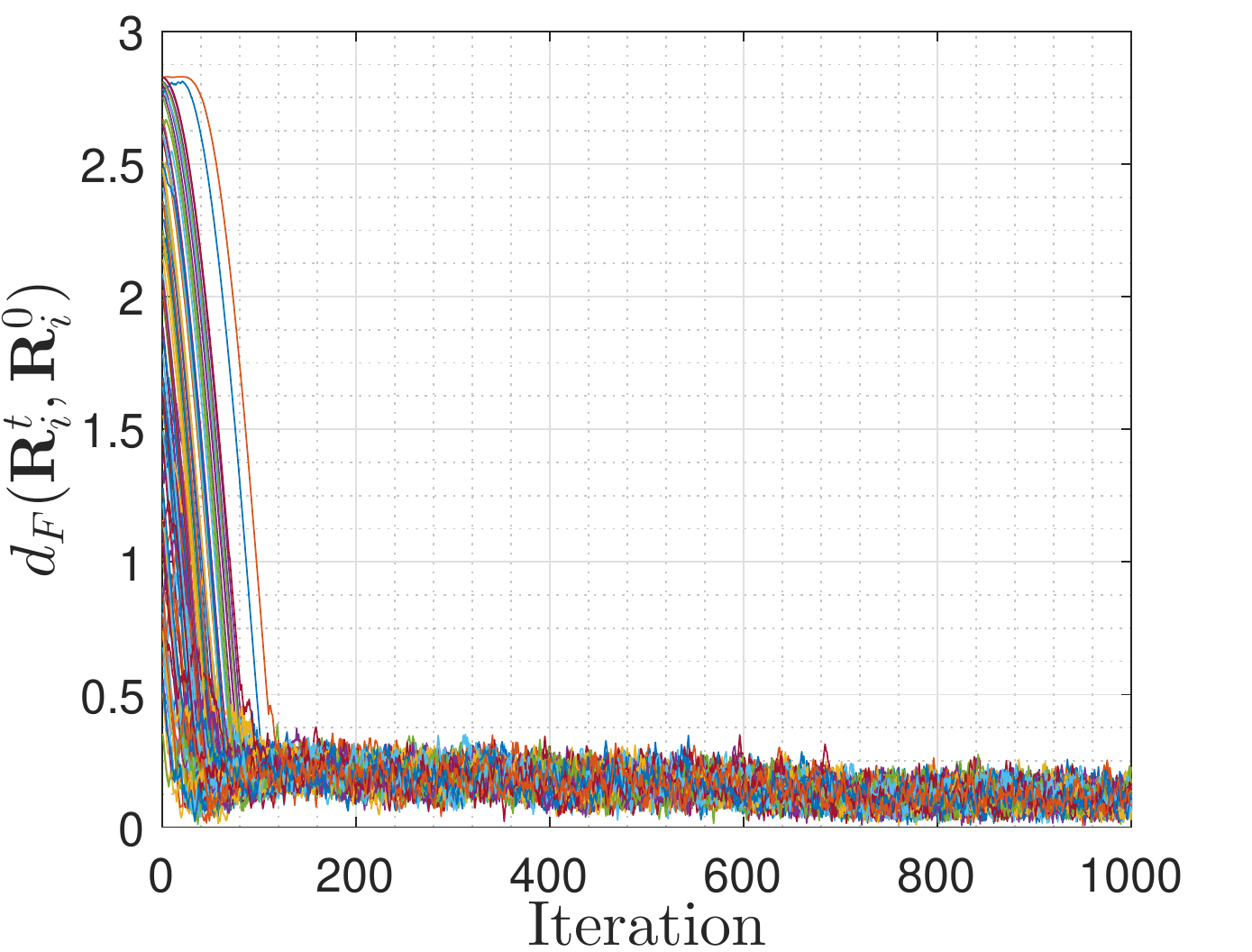}
     }
    \subfigure[]{
    \includegraphics[width=0.45\textwidth]{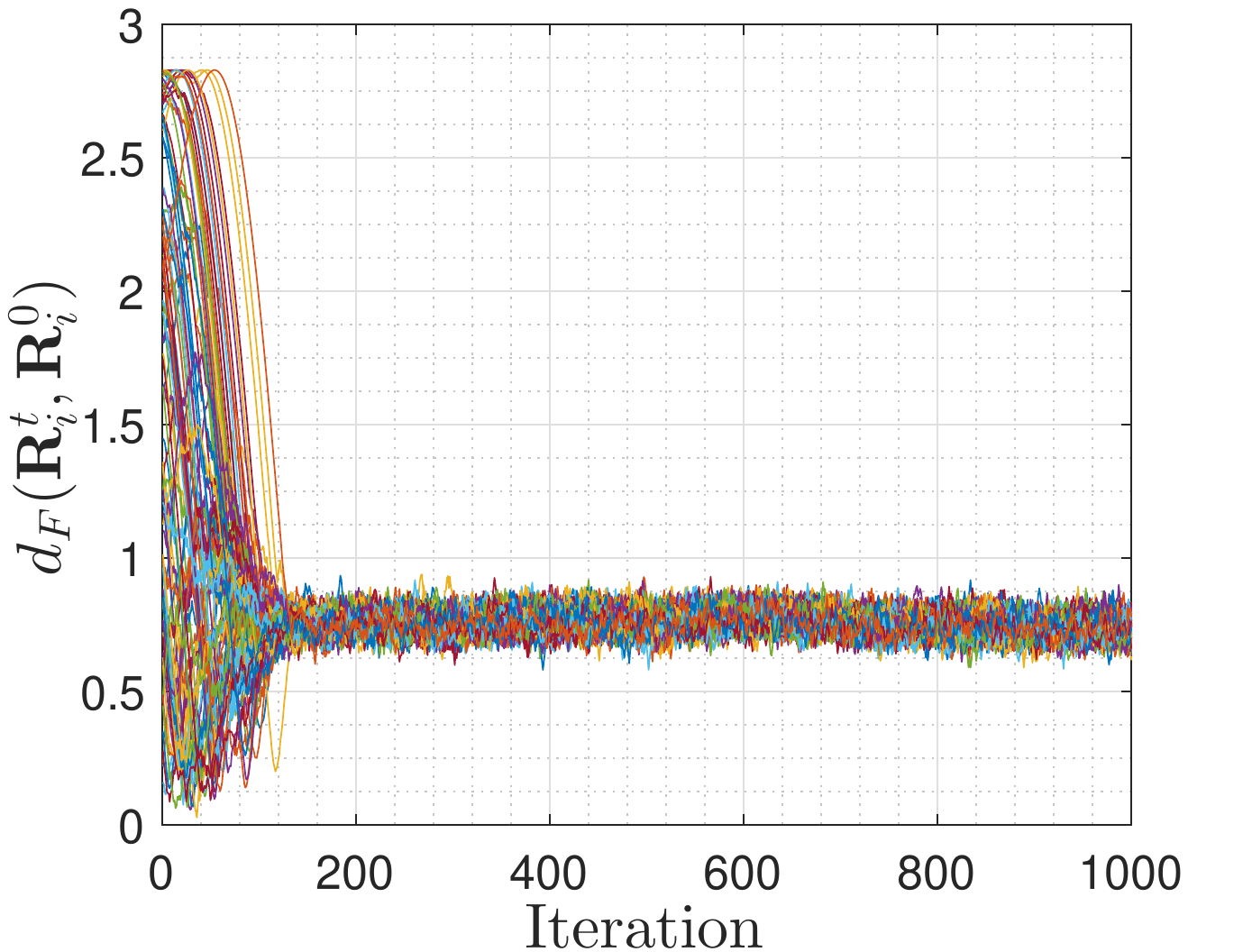}
     }
\caption{Synchronization on the Erd\"{o}s-R\'{e}yni random graph $G(n,p)$ with $n=100$, $p=0.05$, $\eta_{t}=\eta=0.05$ and $\alpha=10^{-3}$. Panel (a): anchored synchronization with a single anchor $\tilde{V}=\{i\}$. Panel (b): anchor-free synchronization.}
\label{Fig:SO3_1}
\vspace{4mm}
\minipage{0.33\textwidth}
  \includegraphics[trim={1cm 0cm 1cm 1cm},width=\linewidth]{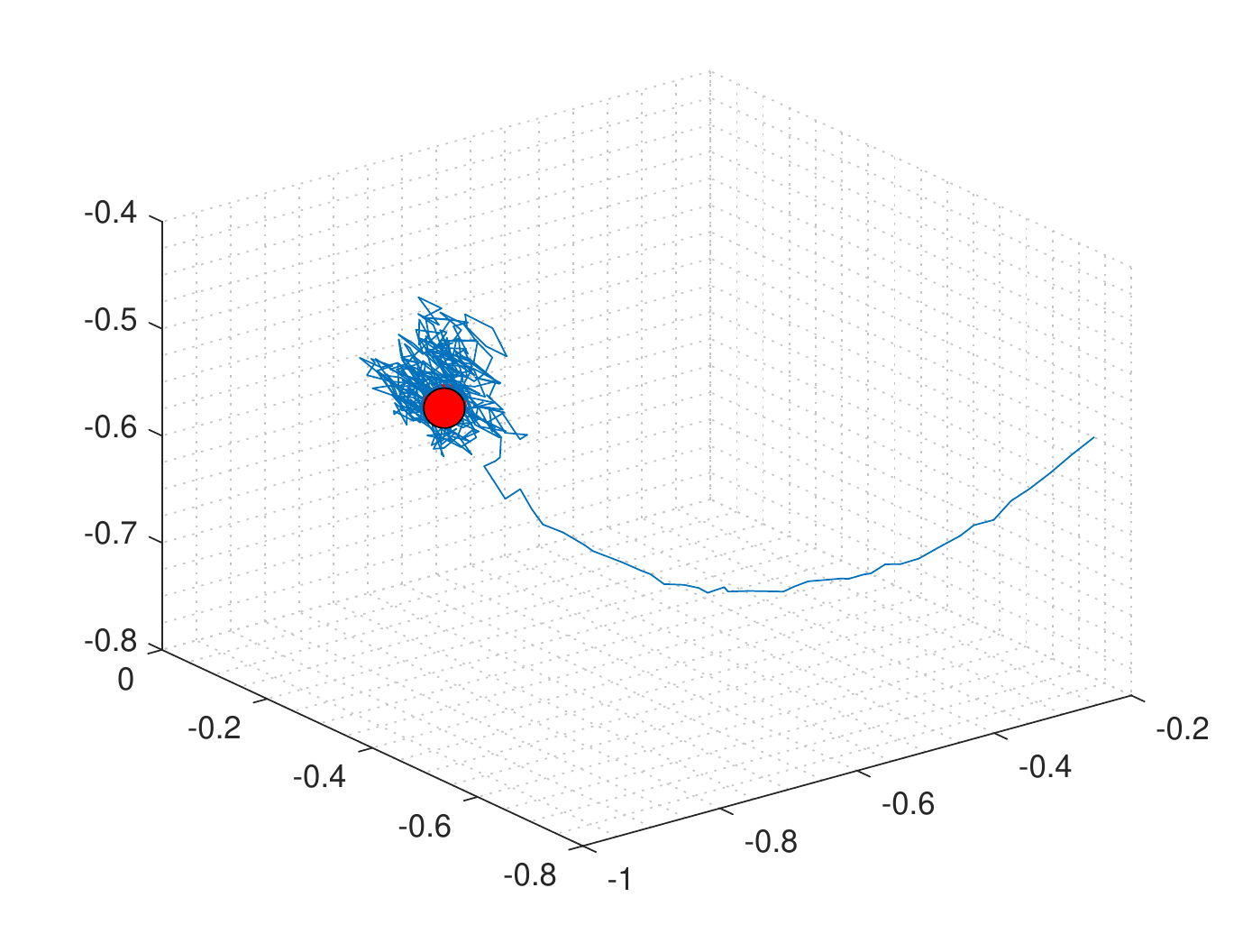}
\endminipage\hfill
\minipage{0.33\textwidth}
  \includegraphics[trim={1cm 0cm 1cm 1cm},width=\linewidth]{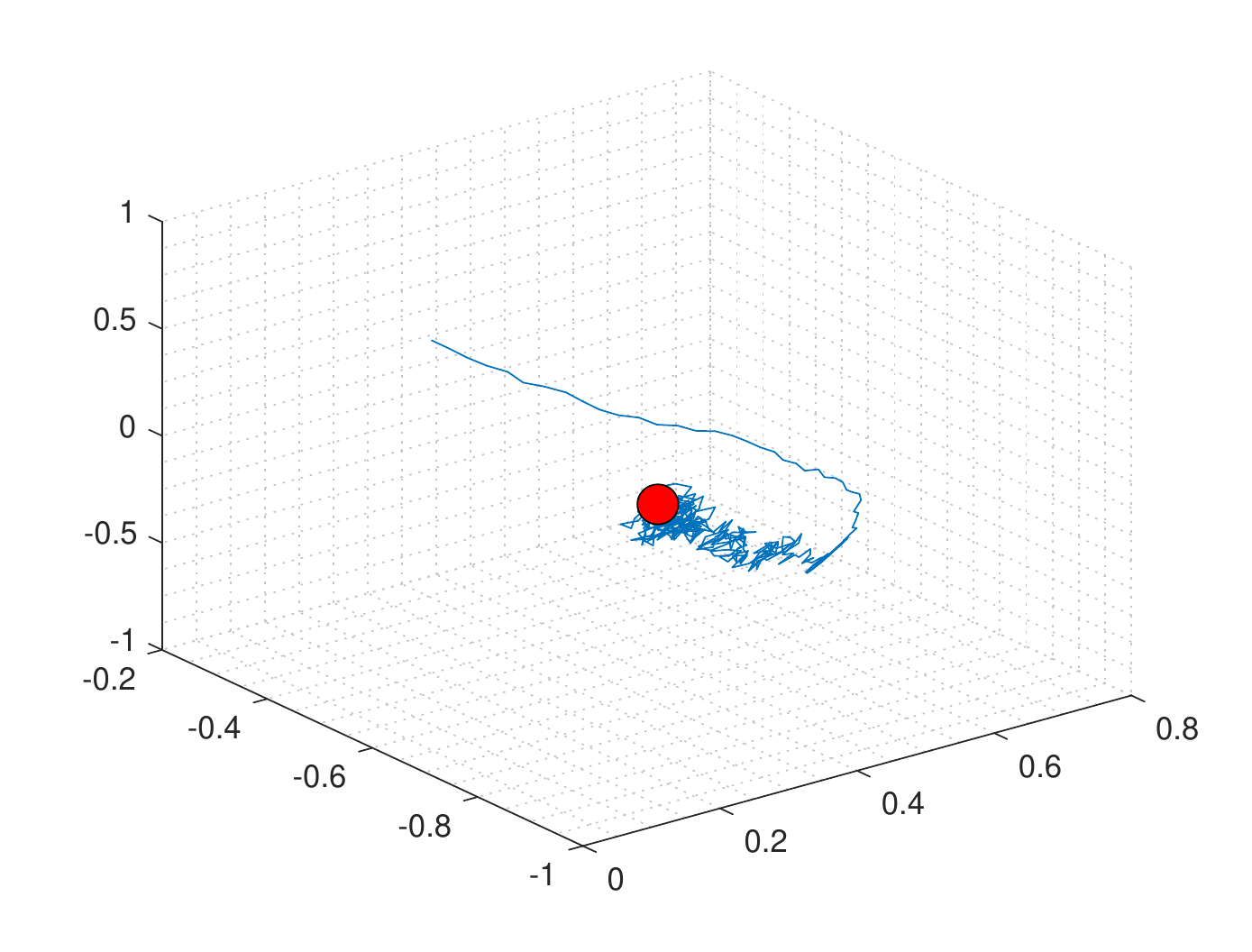}
\endminipage\hfill
\minipage{0.33\textwidth}%
  \includegraphics[trim={1cm 0cm 1cm 1cm},width=\linewidth]{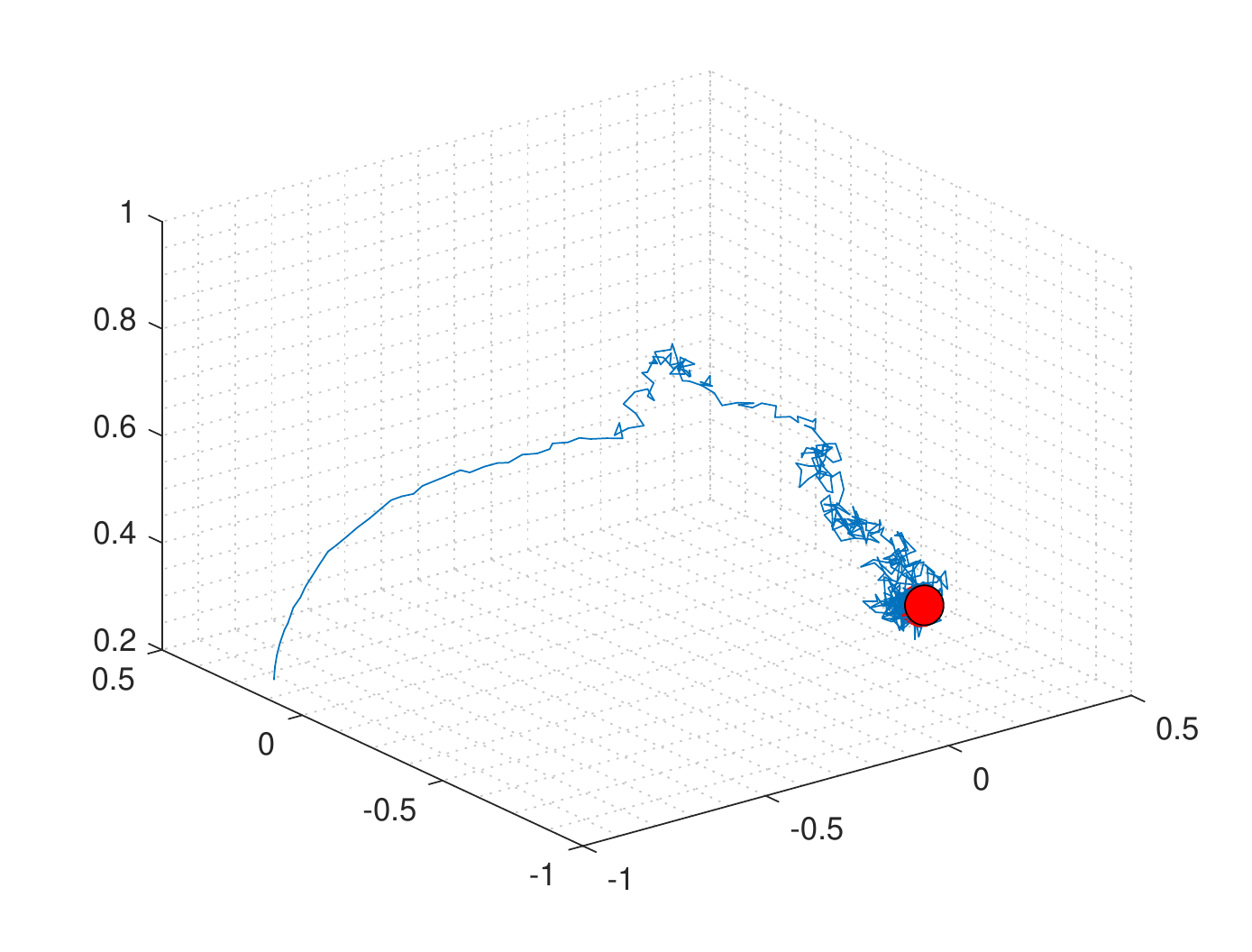}
\endminipage
\\ \vspace{10mm}
\minipage{0.33\textwidth}
  \includegraphics[trim={1cm 1cm 1cm .5cm},width=\linewidth]{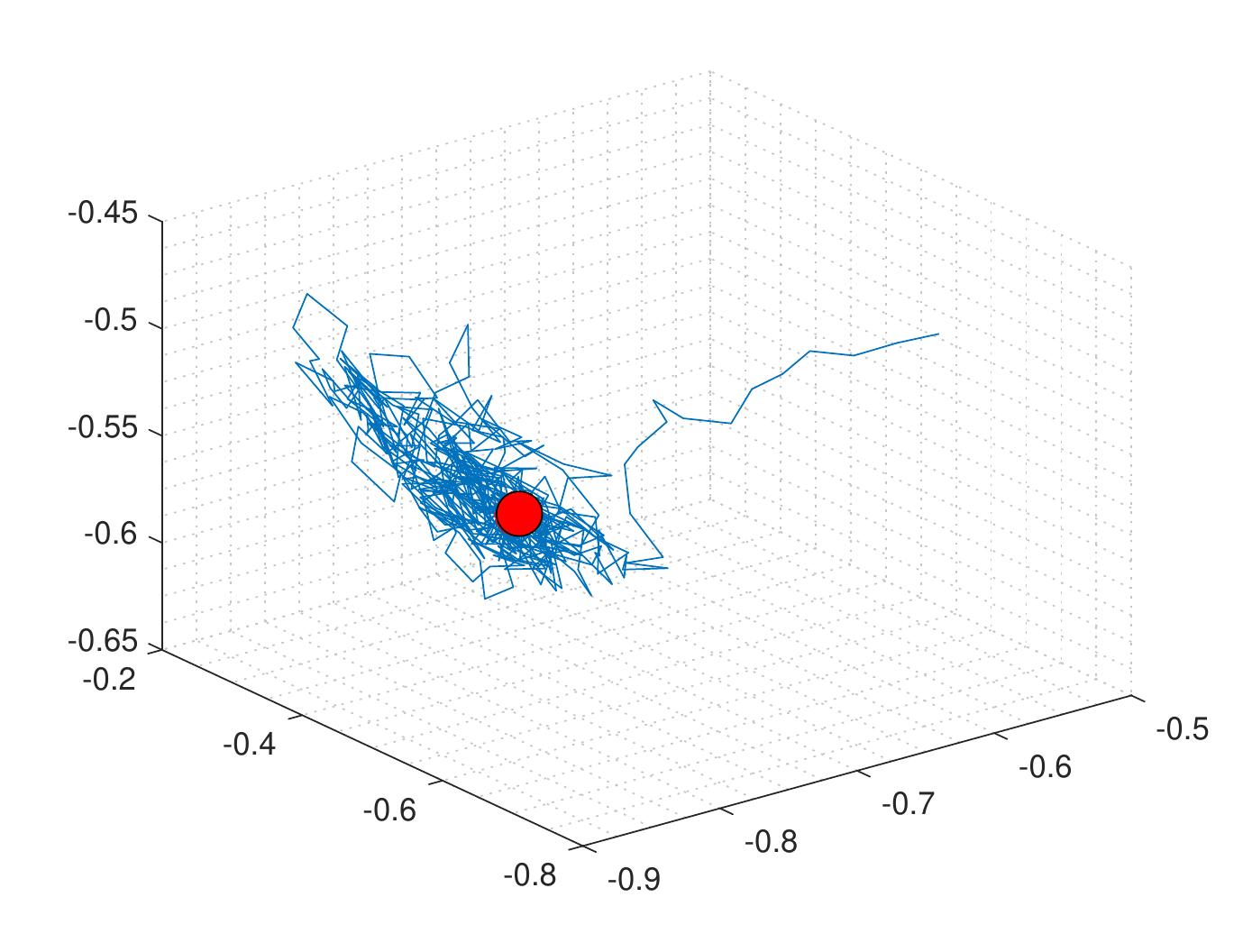}
\endminipage\hfill
\minipage{0.33\textwidth}
  \includegraphics[trim={1cm 1cm 1cm .5cm},width=\linewidth]{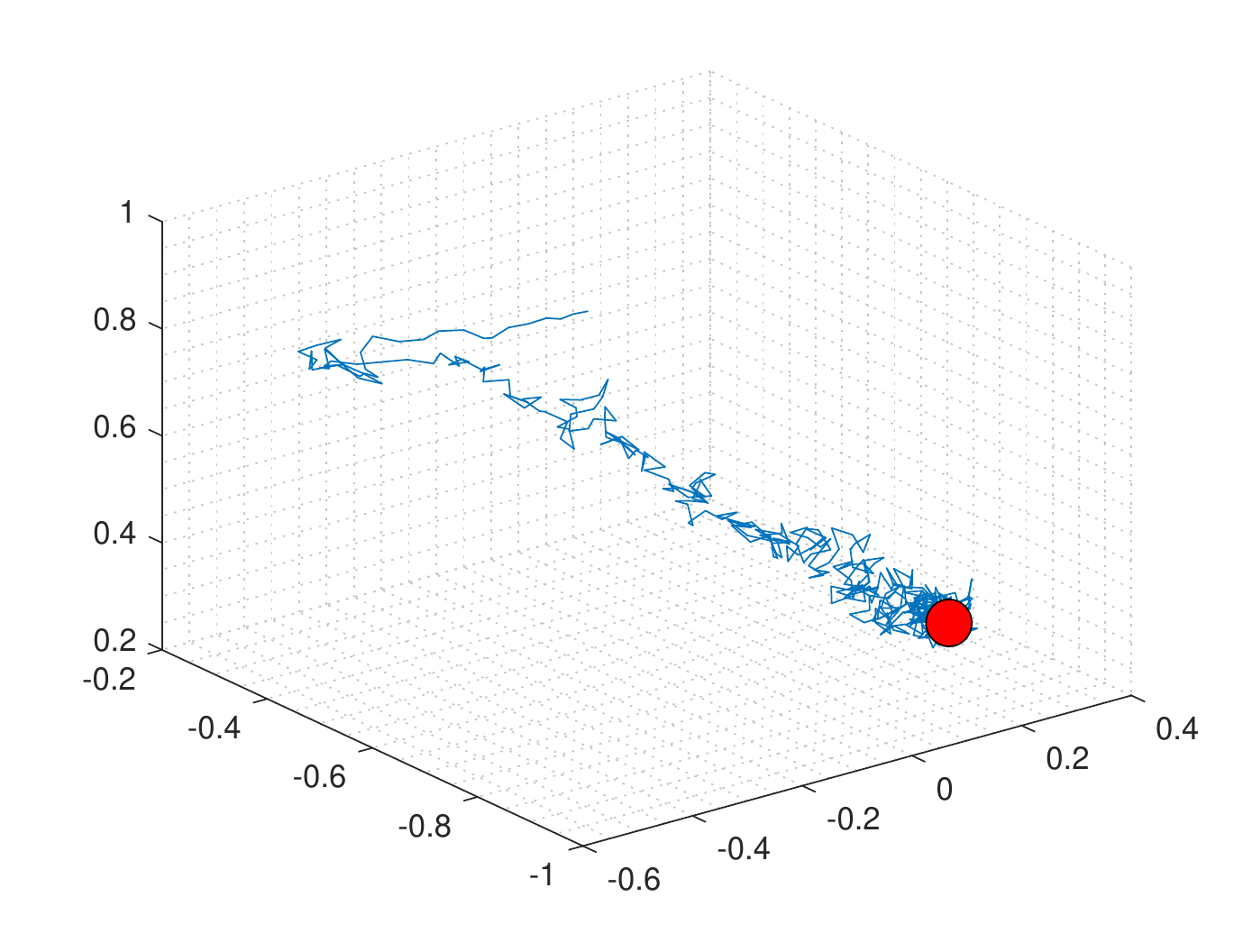}
\endminipage\hfill
\minipage{0.33\textwidth}%
  \includegraphics[trim={1cm 1cm 1cm .5cm},width=\linewidth]{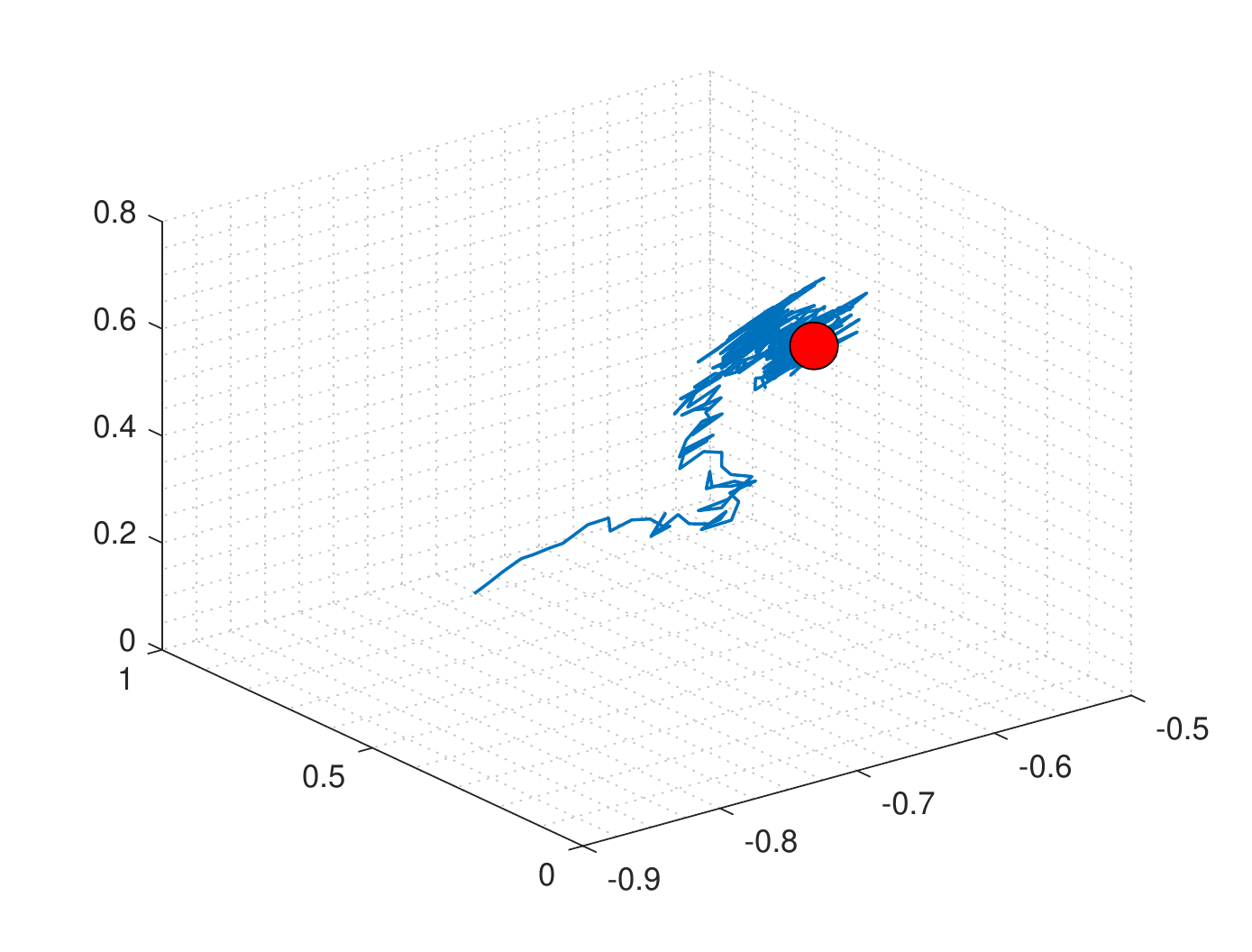}
\endminipage
\caption{The trajectories (blue curves) of column vectors of the matrix $\bm{R}_{1}^{t}$ in $\real^{3}$ at different iterations $t=1,2,\cdots,T$ and for two different trials. The column vectors of the true state $\bm{R}_{1}^{0}$ are shown by the red filled circles, where $i=1\not\in \tilde{V}$ is \textit{not} the anchor node.  From left to right: first column, second column, and third column.}
\label{Fig:SO3_2}
\end{figure}

Let $\beta=10$, and consider the Erd\"{o}s-R\'{e}yni random graph $G(n,p)$ with $n=100$ and $p=0.05$. Further, let $\eta_{t}=\eta=0.05$ and $\alpha=10^{-3}$. Figure \ref{Fig:SO3_1} shows the Frobenius distance $d_{F}(\bm{R}^{t},\bm{R}_{i}^{0})=\|\bm{R}^{t}_{i}-\bm{R}_{i}^{0}\|_{F}$ for all $i\in V$ and at different iterations $t=1,2,\cdots,1000$, where both anchored and anchor-free synchronizations are considered. In our simulations, we used the rejection sampling method to sample from the density function \eqref{Eq:F_function}. 

In the case of anchored synchronization, we consider only one anchor $\tilde{V}=\{i\}$, where $i\in V$ is chosen randomly. We observe from Figure \ref{Fig:SO3_1} that in the absence of any anchors, the algorithm does not converge to the true rotation states $\bm{R}_{i}^{0},i\in V$ due to the invariance of the MSE under rotation. Nevertheless, by including an anchor, the symmetry of the MSE is broken and the algorithm correctly detects the true underlying states.  

Figure \ref{Fig:SO3_2} shows the trajectories of the column vectors of the first node in the graph $\bm{R}_{1}^{t}=[\bm{R}_{11}^{t}|\bm{R}_{12}^{t}|\bm{R}_{13}^{t}]$ with a blue curve at different iterations $t=1,2,\cdots,T$ and for two different trials of Algorithm \ref{alg:4}.  Further, the column vectors of the true matrix $\bm{R}_{1}^{0}=[\bm{R}_{11}^{0}|\bm{R}_{12}^{0}|\bm{R}_{13}^{0}]$ are shown by the red filled circles. Clearly, the trajectory of each column vector $\bm{R}_{1i},i=1,2,3$ clusters around its true state $\bm{R}_{1i}^{0}$. Note that the trajectories oscillate around their true states due to using a constant step size $\eta_{t}=\eta=0.05$. However, these oscillations can be suppressed if a smaller step size is used. 

\subsection{Weighted MAX-CUT Problem}
\label{sec:MAX-CUT}

\subsubsection{Problem Description}
\label{sec:MAX-CUT_Problem_Description}
Consider a weighted graph $G=(V,E,\omega)$ with vertices $V\topdoteq\{1,2,\cdots,n\}$, and edges $E\subseteq V\times V$, where $\omega:E\rightarrow \real_{+}$ assigns a non-negative weight $\omega_{ij}$ to each edge $(i,j)\in E$ and $\omega_{ij}=0$ when $(i,j)\not\in E$. Further, $\omega_{ii}=0,i\in V$ and $\max_{i,j\in [n]}\omega_{ij}<\infty$. The non-negativity condition of the weight coefficients is  inessential and is only imposed to simplify our analysis in the sequel. The weighted MAX-CUT problem is concerned with finding a set of vertices $S$ (not necessarily unique) that maximizes the total weights of the edges in the cut $(S,V/S)$, that is $w(S)=\sum_{i\in S,j\in V/S}\omega_{ij}$. 

The MAX-CUT problem has diverse applications which include statistical physics, where the combinatorial optimization methods are useful for finding the ground states \cite{d1997optimal}. Another application of the weighted MAX-CUT problem is in layout design of very large scale integrated (VLSI) circuits, where the goal is to minimize the number of locations that wires change layers to minimize the risk of breaking the chip. 

To find the maximum cut $w(S)$, the following integer quadratic program can be formulated
\begin{subequations}
\begin{align}
\label{Eq:integerprogram-1}
&\max_{x_{1},\cdots,x_{n}}\ \dfrac{1}{2n}\sum_{i<j}\omega_{ij}(1-x_{i}x_{j})\\ \label{Eq:integerprogram-2}
&\text{subject to}: x_{i}\in \{-1,1\}, \quad \text{for all}\ i\in V,
\end{align}
\end{subequations}
where $x_{i}=1$ if $i\in S$, and $x_{i}=-1$ if $i\in V/S$. Since the problem in eqs. \eqref{Eq:integerprogram-1}-\eqref{Eq:integerprogram-2} is NP-complete, an approximate algorithm based on the following relaxation of eqs. \eqref{Eq:integerprogram-1}-\eqref{Eq:integerprogram-2} is proposed in \cite{goemans1995improved},
\begin{subequations}
\begin{align}
\label{Eq:relaxation-1}
&\max_{\bm{\sigma}_{1},\cdots,\bm{\sigma}_{n}}\ \dfrac{1}{2n}\sum_{i<j}\omega_{ij}(1-\langle \bm{\bm{\sigma}}_{i},\bm{\bm{\sigma}}_{j} \rangle)\\ \label{Eq:relaxation-2}
&\text{subject to}: \bm{\bm{\sigma}}_{i}\in {\rm S}^{d}, \quad \text{for all}\ i\in V,
\end{align}
\end{subequations}

A different relaxation of the integer programming problem in \eqref{Eq:integerprogram-1}-\eqref{Eq:integerprogram-2} is also established in \cite{goemans1995improved},
\begin{subequations}
\begin{align}
\label{Eq:relaxation_m-1}
&\max_{\bm{\sigma}_{1},\cdots,\bm{\sigma}_{n}} {1\over n}\sum_{i<j}\omega_{ij}\dfrac{\arccos(\langle\bm{\sigma}_{i},\bm{\sigma}_{j}\rangle)}{\pi}\\
\label{Eq:relaxation_m-2}
&\text{subject to}: \bm{\sigma}_{i}\in {\rm S}^{d}, \quad \text{for all}\ i\in V,
\end{align}
\end{subequations}

The relaxation in eqs. \eqref{Eq:relaxation_m-1}-\eqref{Eq:relaxation_m-2} has received less attention compared to the formulation in eqs. \eqref{Eq:relaxation-1}-\eqref{Eq:relaxation-2}, since it does not admit a semi-definite programming solution as in eqs. \eqref{Eq:relaxation-1}-\eqref{Eq:relaxation-2}. However, we show that the optimization methods on manifolds can be employed effectively to search for good local minima of \eqref{Eq:relaxation_m-1}-\eqref{Eq:relaxation_m-2}. 

It is also noteworthy that from a pure geometric perspective, the relaxation in \eqref{Eq:relaxation_m-1}-\eqref{Eq:relaxation_m-2} is more intriguing than \eqref{Eq:relaxation-1}-\eqref{Eq:relaxation-2} since it has a nice geometric interpretation in terms of maximizing the mutual distance of $n$ points on the sphere as we discuss later. 

\subsubsection{Semi-definite Relaxation \cite{goemans1995improved}}
\label{subseq:semi-definite}
The problem in eqs. \eqref{Eq:relaxation-1}-\eqref{Eq:relaxation-2} is non-convex, due to a non-convex feasible set. A semi-definite programming (SDP) relaxation for this problem can be characterized by defining the matrix $\bm{Y}\topdoteq[Y_{ij}]$ with elements $Y_{ij}\topdoteq\langle \bm{\bm{\sigma}}_{i},\bm{\bm{\sigma}}_{j} \rangle$, as well as $\bm{W}\topdoteq[\omega_{ij}]$. The optimization in \eqref{Eq:relaxation-1}-\eqref{Eq:relaxation-2} can now be rewritten as a trace maximization problem over elliptope,
\begin{subequations}
\begin{align}
\label{Eq:SDP-1}
&\max_{\bm{Y}}\ \dfrac{1}{4}\text{Tr}(\bm{W}(\bm{1}-\bm{Y}))\\ \label{Eq:SDP-2}
&\text{subject to}: \text{diag}(\bm{Y})=\bm{I},\quad \bm{Y}\succeq 0,
\end{align}
\end{subequations}
where $\text{diag}(\bm{Y})$ is a diagonal matrix, where diagonal elements are $\{Y_{ii}\}_{i=1}^{n}$ and the off-diagonal elements are all zeros. Further, $\bm{1}$ and $\bm{I}$ denote the all ones and identity matrices, respectively. Given the solution $\bm{Y}_{\ast}$ of \eqref{Eq:SDP-1}-\eqref{Eq:SDP-2}, we compute the eigenvalue decomposition $\bm{Y}_{\ast}=\bm{U}\bm{\Sigma} \bm{U}^{T}$, where $\bm{\Sigma}=\text{diag}(\lambda_{1},\cdots,\lambda_{n})$ is a diagonal matrix with the eigenvalues $\lambda_{1}\geq \lambda_{2}\geq \cdots \geq \lambda_{n}$ as the diagonal elements, and $\bm{U}=(\bm{U}_{1},\cdots,\bm{U}_{n})$ is a column matrix composed of the eigenvectors $\bm{U}_{1},\cdots,\bm{U}_{n}$. We then choose the principle component $\bm{U}_{1}=(U_{11},U_{12},\cdots,U_{1n})$ and select the following cut
\begin{align*}
S=\{i\in V: \text{sign}(U_{1i})>0\}.
\end{align*}

Note that the spectral rounding method we used here is due to \cite{javanmard2016phase} and is different from the hyperplane rounding technique of \cite{goemans1995improved}. 

\subsubsection{Primal-Dual Method on the Product Manifold}
\label{subsection:produtct_manifold}
As we have already observed in the PCA problem, in practice the SDP method hardly scales beyond $n$ of a few hundreds, although it yields the best known integrality gap for the MAX-CUT problem.\footnote{In fact, the unique game conjecture \cite{khot2007optimal}, which is a statement about the optimal inapproximability of certain class of computer science problems, implies that the SDP method is optimal for the MAX-CUT problem in the sense that it yields the conjectured optimal integrality gap of approximately $0.878$.} Recently, Javanmard, Montanari, and Ricci-Tersenghi \cite{javanmard2016phase} have shown empirically that a heauristic gradient ascent method for the closely related community detection problem converges to good local optima that are undistinguishable from the SDP optimum. Montanari \cite{montanari2016grothendieck} has further proved that the local optima of eqs. \eqref{Eq:relaxation-1}-\eqref{Eq:relaxation-2} are indeed within a gap from the SDP optimum in \eqref{Eq:SDP-1}-\eqref{Eq:SDP-2}, where the gap is proportional to $n/\sqrt{d}$. Motivated by those observations and with the understanding that the objective value of eqs. \eqref{Eq:relaxation-1}-\eqref{Eq:relaxation-2} does not have spurious local minima, we apply the Riemannian optimization techniques to the MAX-CUT problem. 

In the sequel, we characterize two different algorithms for the MAX-CUT problem using the two different MAX-CUT relaxations we discussed in Section \ref{sec:MAX-CUT_Problem_Description}. In both algorithms, we allow for certain restrictions on the desired cut on the graph. In particular, given the sub-graph $\tilde{E}\subset E$, we find a cut $(S,V/S)$ that maximizes the objective function under the restriction that it includes edges belonging to $\tilde{E}$, \textit{i.e.}, either $(k,l)\in S\times V/S$ or $(l,k)\in S\times V/S$ for $(l,k)\in \tilde{E}$. We remark that the MAX-CUT problem under graph constraints has also been investigated in \cite{lee2016max} using a linear programming relaxation based on the Sherali-Adams hierarchy.
\begin{itemize}[leftmargin=*]

\item \textbf{Randomized Hyper-Plane Rounding Method}: We consider an optimization problem on the manifold of direct product of $d$-spheres $\mathcal{M}={\rm S}^{d}\otimes {\rm S}^{d}\otimes \cdots \otimes{\rm S}^{d}\cong {\rm S}^{nd} $ equipped with the product metric $g_{\mathcal{M}}=g\oplus g \oplus \cdots\oplus g$ and the product topology, where $g$ is the metric induced by the embedding of each sphere in $\real^{d+1}$. We now study a relaxation of eqs. \eqref{Eq:integerprogram-1}-\eqref{Eq:integerprogram-2} as follows
\begin{subequations}
\begin{align}
\label{Eq:doublecone-1}
&\text{Maximize}_{\bm{\sigma}_{1}\times \cdots\times \bm{\sigma}_{n}\in \mathcal{M}}\ \dfrac{1}{2n}\sum_{i<j}\omega_{ij}(1-\langle \bm{\bm{\sigma}}_{i},\bm{\sigma}_{j} \rangle)\\ \label{Eq:doublecone-2}
&\text{subject to}: \langle \bm{\sigma}_{k},\bm{\sigma}_{l} \rangle=-1, \quad \text{for all}\ (k,l)\in \tilde{E}. 
\end{align}
\end{subequations}
Notice that replacing the equality constraint \eqref{Eq:doublecone-2} with the inequality constraint $\langle \bm{\sigma}_{k},\bm{\sigma}_{l} \rangle\leq -1$ does not change the problem. Given a locally optimal solution $\bm{\sigma}^{\ast}=(\bm{\sigma}_{1}^{\ast},\cdots,\bm{\sigma}_{n}^{\ast})$ of eqs. \eqref{Eq:doublecone-1}-\eqref{Eq:doublecone-2}, we next generate a sequence of random vectors $\bm{u}_{1},\cdots,\bm{u}_{N}$ $\in {\rm S}^{d}$ and find the cut $S_{k}=\{i\in V:\langle \bm{\sigma}^{\ast}_{i},\bm{u}_{k} \rangle\geq 0\},k=1,2,\cdots,N$. Notice that corresponding to each random vector $\bm{u}_{k}$, a separating hyperplanes is defined in $\real^{n+1}$. Further, notice that the constraints in \eqref{Eq:doublecone-2} ensures that, irrespective of a chosen hyperplane, the designated edges of $\tilde{E}$ are included in $S_{k}\times V/S_{k}$ for all $k=1,2,\cdots,N$.
 
Now, we characterize the geometry of the constraint manifold $\mathcal{M}$. For each point $\bm{p}\in {\rm S}^{d}\subset \real^{d+1}$, the tangent space $T_{\bm{p}}{\rm S}^{d}$ is defined in eq. \eqref{Eq:TangentSphere}.  Given $\bm{p}=(p_{1},\cdots,p_{n})\in \mathcal{M}={\rm S}^{d}\otimes
 \cdots \otimes {\rm S}^{d}$, the tangent space of the product manifold takes the following form
\begin{align} 
 T_{\bm{p}}\mathcal{M}&=T_{(p_{1},\cdots,p_{n})}{\rm S}^{d}\otimes
 \cdots \otimes {\rm S}^{d}\\ \label{Eq:isomorphism}
 &\cong T_{p_{1}}{\rm S}^{d} \oplus T_{p_{2}}{\rm S}^{d}\oplus \cdots \oplus T_{p_{n}}{\rm S}^{d},
\end{align}
where the isomorphism is induced by the differential of the projection onto each sphere with the differential of the inclusion map as its inverse. 

Let $\bm{\lambda}\topdoteq(\lambda_{1},\cdots,\lambda_{n})$. We view $\mathcal{M}={\rm S}^{d}\otimes\cdots\otimes {\rm S}^{d}\cong {\rm S}^{n\times d}$ as a sub-manifold of $\real^{n\times(d+1)}$ and define the following Lagrangian function
\begin{align*}
&L(\bm{\sigma},\bm{\lambda}):\real^{n\times(d+1)}\times \real_{+}^{n}\rightarrow \real,\\
&(\bm{\sigma},\bm{\lambda})\mapsto {1\over 2n}\sum_{i<j}\omega_{ij}(1-\langle\bm{\sigma}_{i},\bm{\sigma}_{j}\rangle)+\sum_{(k,l)\in \tilde{E}}\lambda_{kl}(1+\langle\bm{\sigma}_{k},\bm{\sigma}_{l}\rangle)-\dfrac{\alpha}{2}\|\bm{\lambda}\|^{2}_{2}.
\end{align*}
The directional derivative of $L$ in the direction of $\bm{z}\topdoteq(\bm{z}_{1},\cdots,\bm{z}_{n})^{T}\in \real^{n\times(d+1)}$ can be computed as below
\begin{align*}
dL(\bm{\sigma},\bm{\lambda})(\bm{z})&\topdoteq\lim_{\varepsilon \rightarrow 0}\dfrac{L(\bm{\sigma}+\varepsilon \bm{z},\bm{\lambda})-L(\bm{\sigma},\bm{\lambda})}{\varepsilon}\\
&\stackrel{\rm{(a)}}{=}-{1\over 2n}\sum_{i<j}\omega_{ij}\langle \bm{z}_{i},\bm{\sigma}_{j}\rangle -{1\over 2n}\sum_{i<j}\omega_{ij}\langle \bm{\sigma}_{i},\bm{z}_{j}\rangle\\ &\hspace{4mm}+\sum_{(k,l)\in \tilde{E}}\lambda_{kl}\langle \bm{z}_{k},\bm{\sigma}_{l}\rangle+\sum_{(k,l)\in \tilde{E}}\lambda_{kl}\langle \bm{\sigma}_{k},\bm{z}_{l}\rangle,\\
&\stackrel{\rm{(b)}}{=}-{1\over 2n}\sum_{i,j=1}^{n}\omega_{ij}\langle \bm{z}_{i},\bm{\sigma}_{j}\rangle+2\sum_{(k,l)\in \tilde{E}}\lambda_{kl}\langle \bm{z}_{k},\bm{\sigma}_{l}\rangle,
\end{align*}
where to derive $\rm{(b)}$, we exchanged the role of indices $i,j$ in the second term of $\rm{(a)}$ and used the fact that $\omega_{ij}=\omega_{ji}$, and $\omega_{ii}=0$ for all $i,j\in [n]$. Let $N(i)\topdoteq\{j\in V:(i,j)\in \tilde{E}\}$. We further made a similar observation for the terms with the dual variables. Based on Definition \ref{Definition:gradient}, we compute the gradient as an element of the tangent space $T_{\bm{\sigma}}\real^{n\times(d+1)}\cong \real^{n\times(d+1)}$ as follows
\begin{align*}
dL(\bm{\sigma},\bm{\lambda})(\bm{z})=\text{Tr}(\grad L(\bm{\sigma},\bm{\lambda})\cdot \bm{z}^{T})\ \dot{=} \ \langle \grad L(\bm{\sigma},\bm{\lambda}),\bm{z}\rangle,
\end{align*}
where the $i$-th element of the gradient vector $\grad L(\bm{\sigma},\bm{\lambda})$ is 
\begin{align*}
(\grad L(\bm{\sigma},\bm{\lambda}))_{i}=-{1\over 2n}\sum_{j=1}^{n}{\omega_{ij}\bm{\sigma}_{j}}+2\sum_{j\in N(i)}\lambda_{ij}\bm{\sigma}_{j}.  
\end{align*}
Define $\mathfrak{L}(\bm{\sigma},\bm{\lambda})$ as the restriction of $L(\bm{\sigma},\bm{\lambda})$ to the product manifold $\mathcal{M}$, \textit{i.e.}, $\mathfrak{L}(\bm{\sigma},\bm{\lambda})$ $\dot{=} L(\bm{\sigma},\bm{\lambda})\vert_{\mathcal{M}}$. We now compute $\grad \mathfrak{L}(\bm{\sigma},\bm{\lambda})$. Due to the isomorphism in eq. \eqref{Eq:isomorphism}, it suffices to compute the projection of each coordinate $i\in [n]$ of the vector $\grad L(\bm{\sigma},\bm{\lambda})$ onto its corresponding tangent plane $T_{\bm{\sigma}_{i}}{\rm S}^{d}$. In particular, we have
\begin{align}
\label{Eq:Gradient11}
{\grad \mathfrak{L}}(\bm{\sigma},\bm{\lambda})=\left(\bm{X}_{1},\cdots,\bm{X}_{n}\right),
\end{align}
where 
\begin{align}
\label{Eq:Gradient22}
\bm{X}_{i}&=\Pi_{T_{\bm{\sigma}_{i}}{\rm S}^{d}}\left((\grad L(\bm{\sigma},\bm{\lambda}))_{i}\right)\\ \label{Eq:Gradient33}
&=\dfrac{1}{2n}\sum_{j=1}^{n}\omega_{ij}(\bm{\sigma}_{i}\langle \bm{\sigma}_{i},\bm{\sigma}_{j}\rangle -\bm{\sigma}_{j})+2\sum_{j\in N(i)}\lambda_{ij}(\bm{\sigma}_{i}\langle \bm{\sigma}_{i},\bm{\sigma}_{j}\rangle -\bm{\sigma}_{j}).
\end{align}

To describe the algorithm, we need to specify the stopping criteria. We define the constraint violation associated with the solution vector $\bm{\sigma}=(\bm{\sigma}_{1},\cdots,\bm{\sigma}_{m})\in {\rm S}^{d}\otimes \cdots \otimes {\rm S}^{d}$ as follows
\begin{align}
\label{Eq:stopping_criterion1}
\Delta_{1}(\bm{\sigma})\topdoteq \dfrac{1}{\sqrt{|\tilde{E}|}}\left \|\Pi_{\real_{+}^{|\tilde{E}|}}\big(h_{ij}(\bm{\sigma}_{i},\bm{\sigma}_{j})\big)_{(i,j)\in \tilde{E}}\right\|,
\end{align}
where $h_{ij}(\bm{\sigma}_{i},\bm{\sigma}_{j})\topdoteq 1+\langle \bm{\sigma}_{i},\bm{\sigma}_{j} \rangle$.

We also define the norm of the gradient associated with the primal-dual pair $(\bm{\sigma},\bm{\lambda})$ as follows
\begin{align}
\label{Eq:stopping_criterion}
\Delta_{2}(\bm{\sigma},\bm{\lambda})\topdoteq\dfrac{1}{\sqrt{n}} \left\|\grad \mathfrak{L}(\bm{\sigma},\bm{\lambda})\right\|_{F},
\end{align}
where $\|\cdot\|_{F}$ is the Frobenius norm, and we recall $\bm{\lambda}\topdoteq(\lambda_{1},\cdots,\lambda_{m})$. Based on these stopping criteria, the pseudocode for the MAX-CUT problem is described in Algorithm \ref{alg:2}. 

\begin{remark}
Note that when the step size is small $\eta_{t}\ll 1$, we have $\cos(\eta_{t}\|\bm{s}_{i}^{t}\|)\approx 1$ and $\sin(\eta_{t}\|\bm{s}_{i}^{t}\|)\approx \eta_{t}\|\bm{s}_{i}^{t}\|$ and Step 6 of Algorithm \ref{alg:3} approximately is
\begin{align*}
\bm{\sigma}_{i}^{t+1}\approx \bm{\sigma}_{i}^{t}+\eta_{t}\bm{s}_{i}^{t},
\end{align*}
which is similar to the primal ascent step of the Euclidean primal-dual method.  
\end{remark}

\begin{algorithm}[t!]
\caption{Randomized Hyper-Plane Rounding Method for MAX-CUT}
\label{alg:2}
\begin{algorithmic}[1]
\State{\textbf{require}:  Tolerances $\tt{tol}_{1}$, $\tt{tol}_{2}$. A graph $E\subseteq V\times V$ with weighted edges $\omega_{ij},(i,j)\in E$. A sub-graph $\tilde{E}\subset E$.}
\State{\textbf{initialize}: Choose $\bm{u}_{1},\cdots,\bm{u}_{N},\bm{\sigma}^{0}_{i}\in {\rm S}^{d}$, $\lambda^{0}_{i}=0$ for all $i\in V$, and $\Delta_{1}>\tt{tol}_{1}$ and $\Delta_{2}>\tt{tol}_{2}$. A decreasing sequence $\{\eta_{t}\}_{t=0}^{\infty}$ for the step size as in Thm. \ref{Thm:2}. The regularizer free parameter $\alpha\in \real_{+}$ and the dimension $d>1$.}
\While{$\Delta_{1}>\tt{tol}_{1}$ and $\Delta_{2}>\tt{tol}_{2}$}
\State{Compute $\grad \mathfrak{L}(\bm{\sigma}^{t},\bm{\lambda}^{t})=(\bm{X}^{t}_{1},\cdots,\bm{X}^{t}_{n})$ from eqs. \eqref{Eq:Gradient11},\eqref{Eq:Gradient33}.} 
\State{Update variables as follows:}
\State{$\bm{\sigma}_{i}^{t+1}\leftarrow \cos(\eta_{t}\|\bm{X}_{i}^{t}\|)\bm{\bm{\sigma}}_{i}^{t}+\sin(\eta_{t}\|\bm{X}_{i}^{t}\|){\bm{X}_{i}^{t}\over\|\bm{X}_{i}^{t}\|},\quad i=1,\cdots,n.$}
\State{$\lambda_{ij}^{t+1}\leftarrow \max\big\{0,(1-\alpha)\lambda_{ij}^{t}+\eta_{t}(1+\langle \bm{\sigma}^{t}_{i},\bm{\sigma}^{t}_{j}\rangle))\big\},\quad	(i,j)\in \tilde{E}$}.
\State{Compute $\Delta_{1}=\Delta_{1}(\bm{\sigma}^{t})$  and $\Delta_{2}=\Delta_{2}(\bm{\sigma}^{t},\bm{\lambda}^{t})$.}
\EndWhile

\hspace{-7mm}\Return $S_{k}=\{i\in \mathcal{V}:\langle \bm{\sigma}^{t+1}_{i},\bm{u}_{k} \rangle\geq 0\},k=1,2,\cdots,N$.
\end{algorithmic}
\end{algorithm}
\item \textbf{Randomized Distance Rounding Method}: Here, we consider a different approach for the MAX-CUT problem based on the relaxation \eqref{Eq:relaxation_m-1}-\eqref{Eq:relaxation_m-2}. To apply geometric methods to this relaxation, notice that the intrinsic distance between two points on the sphere $\bm{p},\bm{q}\in {\rm S}^{d}$ is defined by
\begin{align*}
d(\bm{p},\bm{q})=\arccos(\langle \bm{p},\bm{q}\rangle).
\end{align*}
Consequently, the problem \eqref{Eq:relaxation_m-1}-\eqref{Eq:relaxation_m-2} can be regarded as maximizing the mutual distance between $n$ points on the sphere. In particular, the relaxation in \eqref{Eq:relaxation_m-1}-\eqref{Eq:relaxation_m-2} can now be rewritten as 
\begin{subequations}
\begin{align}
\label{Eq:distance_m-1}
&\text{Maximize}_{(\bm{\sigma}_{1},\cdots,\bm{\sigma}_{n})\in {\rm S}^{d}\otimes \cdots \otimes {\rm S}^{d}} {1\over n}\sum_{i<j}\omega_{ij}\dfrac{d(\bm{\sigma}_{i},\bm{\sigma}_{j})}{\pi}\\
\label{Eq:distance_m-2}
&\text{subject to}: d(\bm{\sigma}_{l},\bm{\sigma}_{k})=\pi, \quad (k,l)\in \tilde{E},
\end{align}
\end{subequations}
where we included the graph based constraints on the pairs $(\bm{\sigma}_{l},\bm{\sigma}_{k}),(k,l)\in \tilde{E}$. Notice that the constant factor $\text{diam}({\rm S}^{d})=\pi$ in the denominator is the diameter of the unit sphere, where $\text{diam}(\mathcal{M})\topdoteq \max_{\bm{p},\bm{q}\in \mathcal{M}}d(\bm{p},\bm{q})$. As a result, the problem in \eqref{Eq:distance_m-1}-\eqref{Eq:distance_m-2} can be extended to an optimization on the general smooth manifold $\mathcal{M}$,
\begin{subequations}
\begin{align}
\label{Eq:distance_abstract-1}
&\text{Maximize}_{(\bm{\sigma}_{1},\cdots,\bm{\sigma}_{n})\in \mathcal{M}\otimes \cdots \otimes \mathcal{M}} {1\over n}\sum_{i<j}\omega_{ij}\dfrac{d(\bm{\sigma}_{i},\bm{\sigma}_{j})}{\text{diam}(\mathcal{M})}\\
\label{Eq:distance_abstract-2}
&\text{subject to}: d(\bm{\sigma_{l}},\bm{\sigma}_{k})=\text{diam}(\mathcal{M}), \quad (k,l)\in \tilde{E}.
\end{align}
\end{subequations}
The only restriction on an admissible manifold $\mathcal{M}$ in eqs. \eqref{Eq:distance_abstract-1}-\eqref{Eq:distance_abstract-2} comes from the fact that the diameter must be bounded $\text{diam}(\mathcal{M})<+\infty$. This rules out the hyperbolic manifolds, unless $\mathcal{M}$ is a sub-manifold of a hyperbolic manifold with a bounded diameter. 

Given a locally optimal solution $\bm{\sigma}^{*}\topdoteq (\bm{\sigma}^{\ast}_{1},\cdots,\bm{\sigma}^{\ast}_{n})$ of the problem \eqref{Eq:distance_abstract-1}-\eqref{Eq:distance_abstract-2}, we next generate a random sequence of points $\bm{u}_{1},\cdots,\bm{u}_{N}\in \mathcal{M}$. We then compute $N$ different cuts $S_{1},\cdots,S_{N}\subset V$, where
\begin{align*}
S_{k}=\{i\in V: d(\bm{\sigma}_{i},\bm{u}_{k})\leq \text{diam}(\mathcal{M})/2\}.
\end{align*}

Now, define the Lagrangian function as follows
\begin{align*}
\mathfrak{L}(\bm{\sigma},\bm{\lambda})\topdoteq  {1\over n\cdot \text{diam}(\mathcal{M})}\sum_{i<j}\omega_{ij}{d(\bm{\sigma}_{i},\bm{\sigma}_{j})}+\hspace{-2mm}\sum_{(k,l)\in \tilde{E}}\lambda_{kl}(\text{diam}(\mathcal{M})-d(\bm{\sigma}_{l},\bm{\sigma}_{k}))-\dfrac{\alpha}{2}\|\bm{\lambda}\|_{2}^{2}.
\end{align*}

Since the gradient of the distance function is the logarithmic map $\grad_{\bm{\sigma}_{i}}d(\bm{\sigma}_{i},\bm{\sigma}_{j})=\log_{\bm{\sigma}_{i}}(\bm{\sigma}_{j})$ for $\bm{\sigma}_{i}\not =-\bm{\sigma}_{j}$, the $i$-th element of the gradient of the Lagrangian function is
\begin{align}
\label{Eq:explicit_form}
(\grad_{\bm{\sigma}_{i}}\mathfrak{L}(\bm{\sigma},\bm{\lambda}))_{i}={1\over n\cdot \text{diam}(\mathcal{M})}\sum_{j=1}^{n}\omega_{ij}\log_{\bm{\sigma}_{i}}(\bm{\sigma}_{j})-\sum_{j\in N(i)}\lambda_{ij}\log_{\bm{\sigma}_{i}}(\bm{\sigma}_{j}).
\end{align}

To derive a more explicit form for the gradient $\grad_{\bm{\sigma}_{i}}\mathfrak{L}(\bm{\sigma},\bm{\lambda})$, in the sequel we focus on two types of manifolds, namely \textit{i}) the unit sphere $\mathcal{M}={\rm S}^{d}$, and  \textit{ii}) the rotation group $\mathcal{M}={\rm SO}(d)$. In both cases, the diameter is given by $\text{diam}(\mathcal{M})=\pi$.

In Case $(\textit{i})$, the Riemannian logarithmic map $\log_{\bm{\sigma}_{i}}:T_{\bm{\sigma}_{i}}{\rm S}^{d}\rightarrow {\rm S}^{d}$ takes the following form
\begin{align}
\label{Eq:Case_I_Sphere}
\log_{\bm{\sigma}_{i}}(\bm{\sigma}_{j})= \left\{
    \begin{array}{ll}
        \dfrac{\arccos(\langle \bm{\sigma}_{i},\bm{\sigma}_{j}\rangle)}{\sqrt{1-\langle\bm{\sigma}_{i},\bm{\sigma}_{j} \rangle^{2}}}(\bm{I}-\bm{\sigma}_{i}\bm{\sigma}_{i}^{T})\bm{\sigma}_{j} & \bm{\sigma}_{j}\not\in \{\bm{\sigma}_{i},-\bm{\sigma}_{i}\}  \\
        0 & \bm{\sigma}_{j}=\bm{\sigma}_{i}.
    \end{array}
\right.
\end{align}

In Case $(\textit{ii})$, the Riemannian logarithmic map $\log_{\bm{\sigma}_{i}}:T_{\bm{\sigma}_{i}}\mathfrak{so}(d)\rightarrow {\rm SO}(d)$ can be written from \eqref{Eq:Logarithmic_map_SO} as follows
\begin{align}
\label{Eq:Case_II_SO}
\log_{\bm{\sigma}_{i}}(\bm{\sigma}_{j})= \left\{
    \begin{array}{ll}
        \dfrac{\arccos\left((\text{Tr}(\bm{\sigma}_{i}^{T}\bm{\sigma}_{j})-1)/2\right)}{2\sqrt{1-\left((\text{Tr}(\bm{\sigma}_{i}^{T}\bm{\sigma}_{j})-1)/2\right)^2}}(\bm{\sigma}_{i}^{T}\bm{\sigma}_{j}-\bm{\sigma}_{j}^{T}\bm{\sigma}_{i}) & \bm{\sigma}_{j}\not\in \{\bm{\sigma}_{i},-\bm{\sigma}_{i}\}  \\
        0 & \bm{\sigma}_{j}=\bm{\sigma}_{i}.
    \end{array}
\right.
\end{align}

\begin{algorithm}[t!]
\caption{Randomized Distance Rounding Method for MAX-CUT ($\mathcal{M}={\rm S}^{d}$).}
\label{alg:5}
\begin{algorithmic}[1]
\State{\textbf{require}:  Tolerances $\tt{tol}_{1}$, $\tt{tol}_{2}$. A graph $E\subseteq V\times V$ with weighted edges $\omega_{ij},(i,j)\in E$. A subgraph $\tilde{E}\subset E$.}
\State{\textbf{initialize}: Choose $\bm{\sigma}^{0}_{i}\in {\rm S}^{d}$, $\lambda^{0}_{i}=0$ for all $i\in V$, and $\Delta_{1}>\tt{tol}_{1}$ and $\Delta_{2}>\tt{tol}_{2}$. A decreasing sequence $\{\eta_{t}\}_{t=0}^{\infty}$ for the step size as in Thm. \ref{Thm:2}. The regularizer free parameter $\alpha\in \real_{+}$ and the dimension $d>1$.}
\While{$\Delta_{1}>\tt{tol}_{1}$ and $\Delta_{2}>\tt{tol}_{2}$}
\State{Compute $\grad \mathfrak{L}(\bm{\sigma}^{t},\bm{\lambda}^{t})=(\bm{X}^{t}_{1},\cdots,\bm{X}^{t}_{n})$ from eqs. \eqref{Eq:explicit_form} and \eqref{Eq:Case_I_Sphere}.} 
\State{Update variables as follows:}
\State{$\bm{\sigma}_{i}^{t+1}\leftarrow \cos(\eta_{t}\|\bm{X}_{i}^{t}\|)\bm{\sigma}_{i}^{t}+\sin(\eta_{t}\|\bm{X}_{i}^{t}\|){\bm{X}_{i}^{t}\over\|\bm{X}_{i}^{t}\|},\quad i=1,\cdots,n.$}
\State{$\lambda_{ij}^{t+1}\leftarrow \max\big\{0,(1-\alpha)\lambda_{ij}^{t}+\eta_{t}(\pi-d(\bm{\sigma}^{t}_{i},\bm{\sigma}^{t}_{j}))\big\},\quad	(i,j)\in \tilde{E}$}.
\State{Compute $\Delta_{1}=\Delta_{1}(\bm{\sigma}^{t})$  and $\Delta_{2}=\Delta_{2}(\bm{\sigma}^{t},\bm{\lambda}^{t})$.}
\EndWhile

\hspace{-7mm}\Return  $S_{k}=\{i\in V: d(\bm{\sigma}^{t+1}_{i},\bm{u}_{k})\leq \pi/2 \},k=1,2,\cdots,N$.
\end{algorithmic}
\end{algorithm}

\begin{algorithm}[t!]
\caption{Randomized Distance Rounding Method for MAX-CUT ($\mathcal{M}={\rm SO}(d)$).}
\label{alg:6}
\begin{algorithmic}[1]
\State{\textbf{require}:  Tolerances $\tt{tol}_{1}$, $\tt{tol}_{2}$. A graph $E\subseteq V\times V$ with weighted edges $\omega_{ij},(i,j)\in E$. A subgraph $\tilde{E}\subset E$.}
\State{\textbf{initialize}: Choose $\bm{\sigma}^{0}_{i}\in {\rm S}^{d}$, $\lambda^{0}_{i}=0$ for all $i\in V$, and $\Delta_{1}>\tt{tol}_{1}$ and $\Delta_{2}>\tt{tol}_{2}$. A decreasing sequence $\{\eta_{t}\}_{t=0}^{\infty}$ for the step size as in Thm. \ref{Thm:2}. The regularizer free parameter $\alpha\in \real_{+}$ and the dimension $d>1$.}
\While{$\Delta_{1}>\tt{tol}_{1}$ and $\Delta_{2}>\tt{tol}_{2}$}
\State{Compute $\grad \mathfrak{L}(\bm{\sigma}^{t},\bm{\lambda}^{t})=(\bm{X}^{t}_{1},\cdots,\bm{X}^{t}_{n})$ from eqs. \eqref{Eq:explicit_form} and \eqref{Eq:Case_II_SO}.} 
\State{Update variables as follows:}
\State{$\bm{\sigma}^{t+1}_{i}$ $\hspace*{-1.1mm}\leftarrow \hspace*{-1.1mm} \bm{\sigma}^{t}_{i}\left(\bm{I}-\dfrac{\sin(\eta_{t}\|\bm{X}_{i}^{t}\|_{F})}{\|\bm{X}_{i}^{t}\|_{F}}\bm{X}_{i}^{t}+\dfrac{1-\cos(\eta_{t}\|\bm{X}_{i}^{t}\|_{F})}{\|\bm{X}_{i}^{t}\|_{F}^{2}}(\bm{X}_{i}^{t})^{2}\right), i=1,\cdots,n$.}
\State{$\lambda_{ij}^{t+1}\leftarrow \max\big\{0,(1-\alpha)\lambda_{ij}^{t}+\eta_{t}(\pi- d(\bm{\sigma}^{t}_{i},\bm{\sigma}^{t}_{j}))\big\},\quad	(i,j)\in \tilde{E}$}.
\State{Compute $\Delta_{1}=\Delta_{1}(\bm{\sigma}^{t})$  and $\Delta_{2}=\Delta_{2}(\bm{\sigma}^{t},\bm{\lambda}^{t})$.}
\EndWhile

\hspace{-7mm}\Return $S_{k}=\{i\in V: d(\bm{\sigma}^{t+1}_{i},\bm{u}_{k})\leq \pi/2 \},k=1,2,\cdots,N$.
\end{algorithmic}
\end{algorithm}

Now, Algorithms \ref{alg:5} and \ref{alg:6} describe the pseudo-codes of two randomized primal-dual algorithms for the MAX-CUT problem on the unit sphere ${\rm S}^{d}$ and the rotation group ${\rm SO}(d)$, respectively. 

\end{itemize}
\subsubsection{Numerical Results}

Here, we only present the experimental results for the performance of Algorithm \ref{alg:4} as well as the semi-definite programming of Section \ref{subseq:semi-definite}. In particular, we apply Algorithm \ref{alg:4} to obtain the maximum cut in the Erd\"{o}s-R\'{e}nyi random graphs $G(n,p)$, where each edge is connected with the probability of $p\in [0,1]$ independent from every other edge, and $n$ is the number of vertices. We only consider the connected realizations of $G(n,p)$ in the simulations, so that each cut-set $S$ determines a unique partition of the vertices. The weights $\omega_{ij}$ corresponding to the edges $(i,j)\in E$ are i.i.d. and $\omega_{ij}=|\omega|, \omega\sim \mathsf{N}(0,1)$, \textit{i.e.}, $\omega_{ij}$ is distributed according to the folded Gaussian distribution.

\begin{figure}[t!]
\centering
\subfigure[]{
  \includegraphics[trim={.4cm 0cm .4cm .4cm},width=6cm]{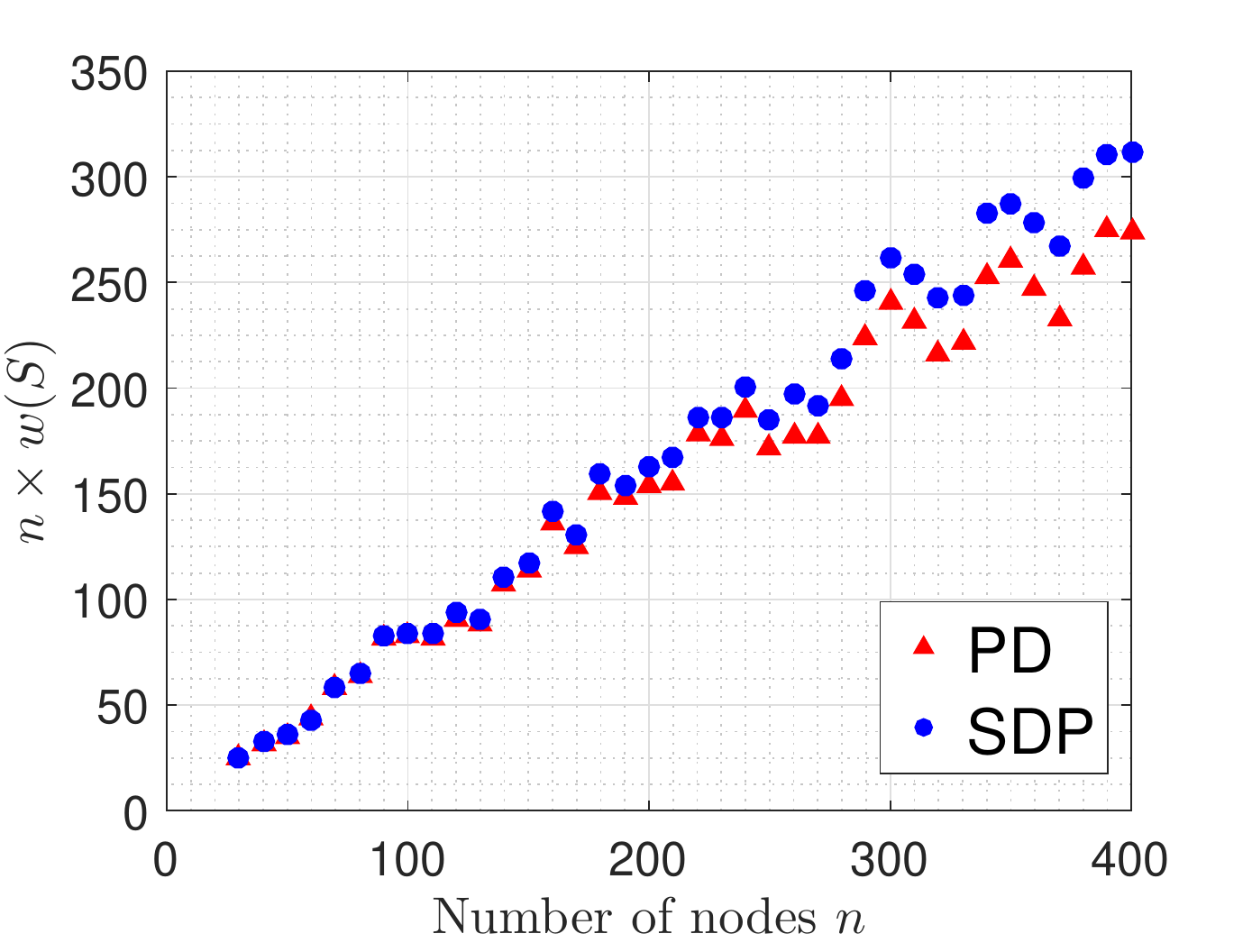}}
\subfigure[]{
  \includegraphics[trim={.4cm 0cm .4cm .4cm},width=6cm]{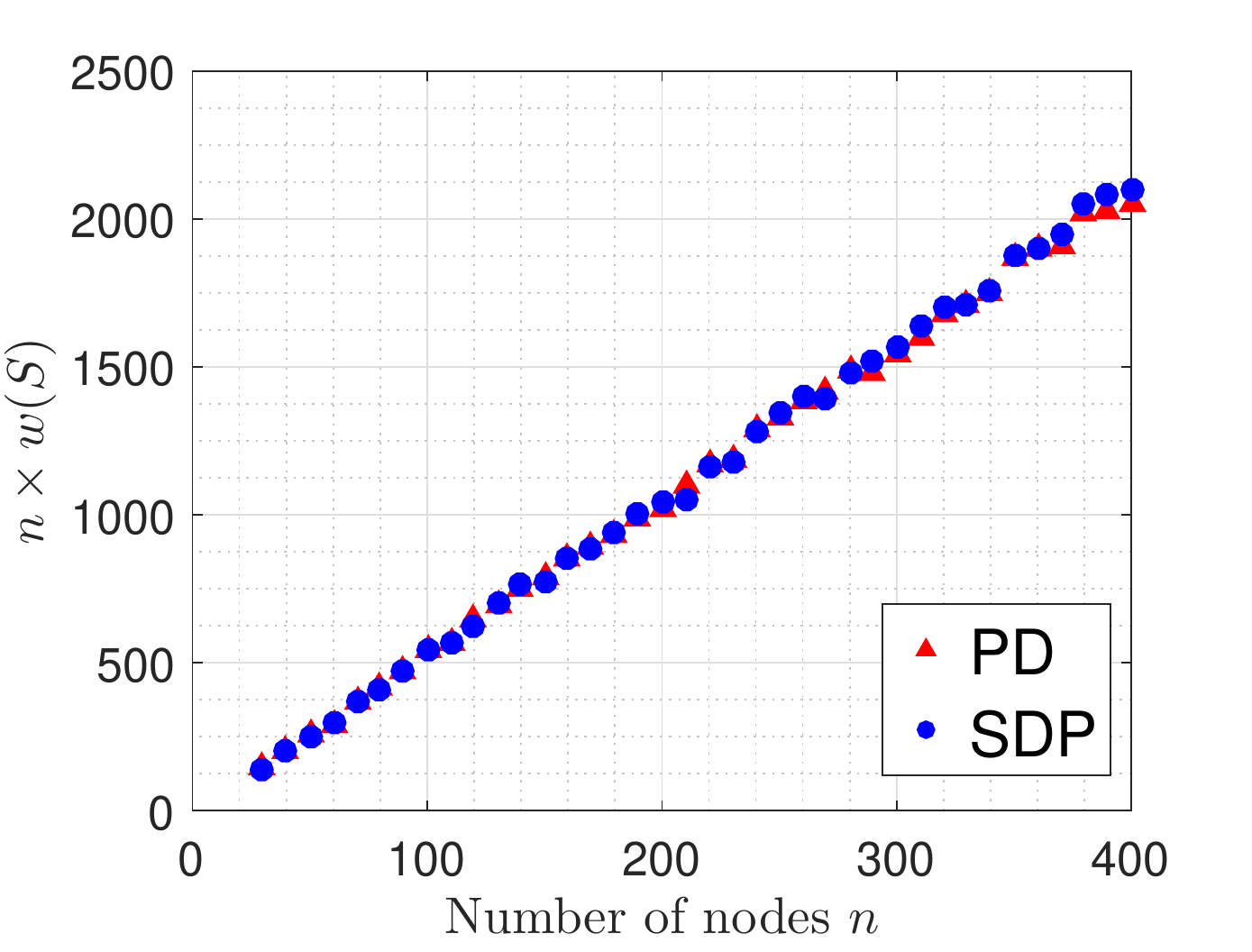}}
\caption{Denormalized objective value $n\times w(S)$ of the Riemannian optimization method for the unconstrained problem $\tilde{E}=\emptyset$ and the SDP algorithm on $G(n,p)$ with different number of nodes $n$. (a) sparse graph $p=1/100$, (b) dense graph $p=1/10$.}
\label{Fig:ObjectiveValues}
\end{figure}

First, we consider the unconstrained settings, where $\tilde{E}=\emptyset$. Since the SDP method cannot be applied to large problems beyond $n$ of a few hundreds, we restrict our analysis to the small graph sizes of $n\in \{30,40,50,\cdots,400\}$ and we consider a fixed choice of the lifting dimension $d=3$. We use \texttt{CVX} \cite{grant2008cvx} with the best precision to solve the SDP problem in eqs. \eqref{Eq:SDP-1}-\eqref{Eq:SDP-2}. Figure \ref{Fig:ObjectiveValues} shows the objective values $w(S)$ attained by the primal-dual algorithm and the SDP method, where $\Delta_{1}=0$ (no constraint), and $\Delta_{2}=10^{-3}$. We also put a limit of $5000$ iterations, after which the algorithm generates the output even if the threshold $\Delta_{2}=10^{-3}$ is not reached.

In our simulations, we let $\alpha=0$, and $\eta_{t}=\eta=1$ for $t\leq 1000$, and $\eta_{t}=\eta=1/100$ for $t>1000$. Further, we generate $N=1000$ samples of random vectors on the sphere $\bm{u}_{1},\cdots,\bm{u}_{N}$ and compute different cuts $S_{1},\cdots,S_{N}$. We then choose the best cut among $S_{1},\cdots,S_{N}$ that maximizes the objective function in \eqref{Eq:integerprogram-1}.     
 
Interestingly, we observe from Figure \ref{Fig:ObjectiveValues} that when the random graph is dense ($p=1/10$), the objective values from the SDP and Riemannian methods is indistinguishable, even when the graph size $n$ grows. Evidently, for sparse graphs $(p=1/100)$, the SDP generally provides better solutions, although the scalability issue of the SDP prohibits its use for larger graphs.

We remark that in all of our simulations with Algorithm \ref{alg:3}, we obtained good solutions, irrespective of the initialization point. This observation is an empirical evidence for our earlier suggestion that all the local minima of the MAX-CUT problem are close to the global minimum or minima. In fact, this observation has already been rigorously proved in \cite{montanari2016grothendieck}.

\begin{figure}[t!]
  \centering
    \subfigure[]{
    \includegraphics[width=0.45\textwidth]{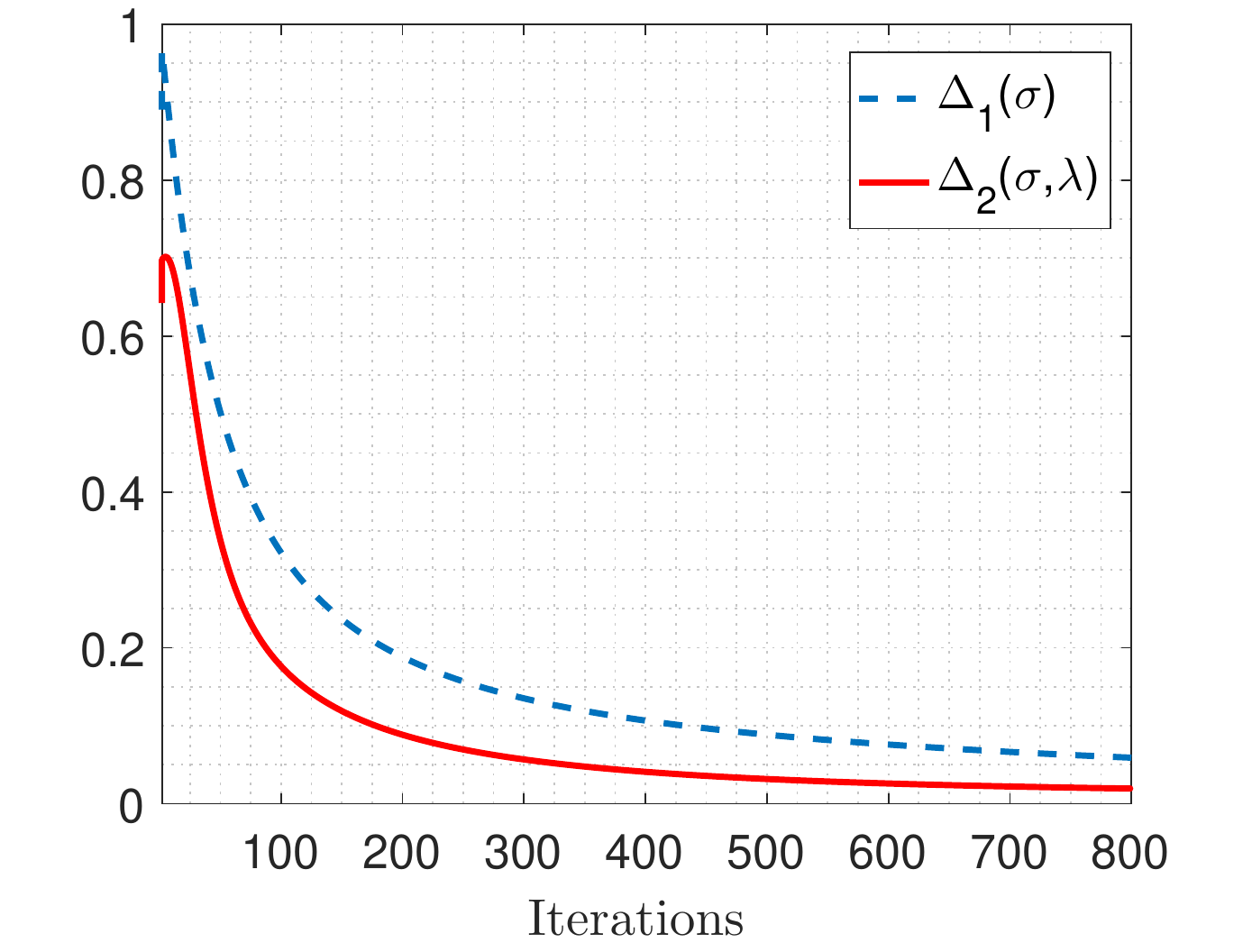}
     }
    \subfigure[]{
    \includegraphics[width=0.45\textwidth]{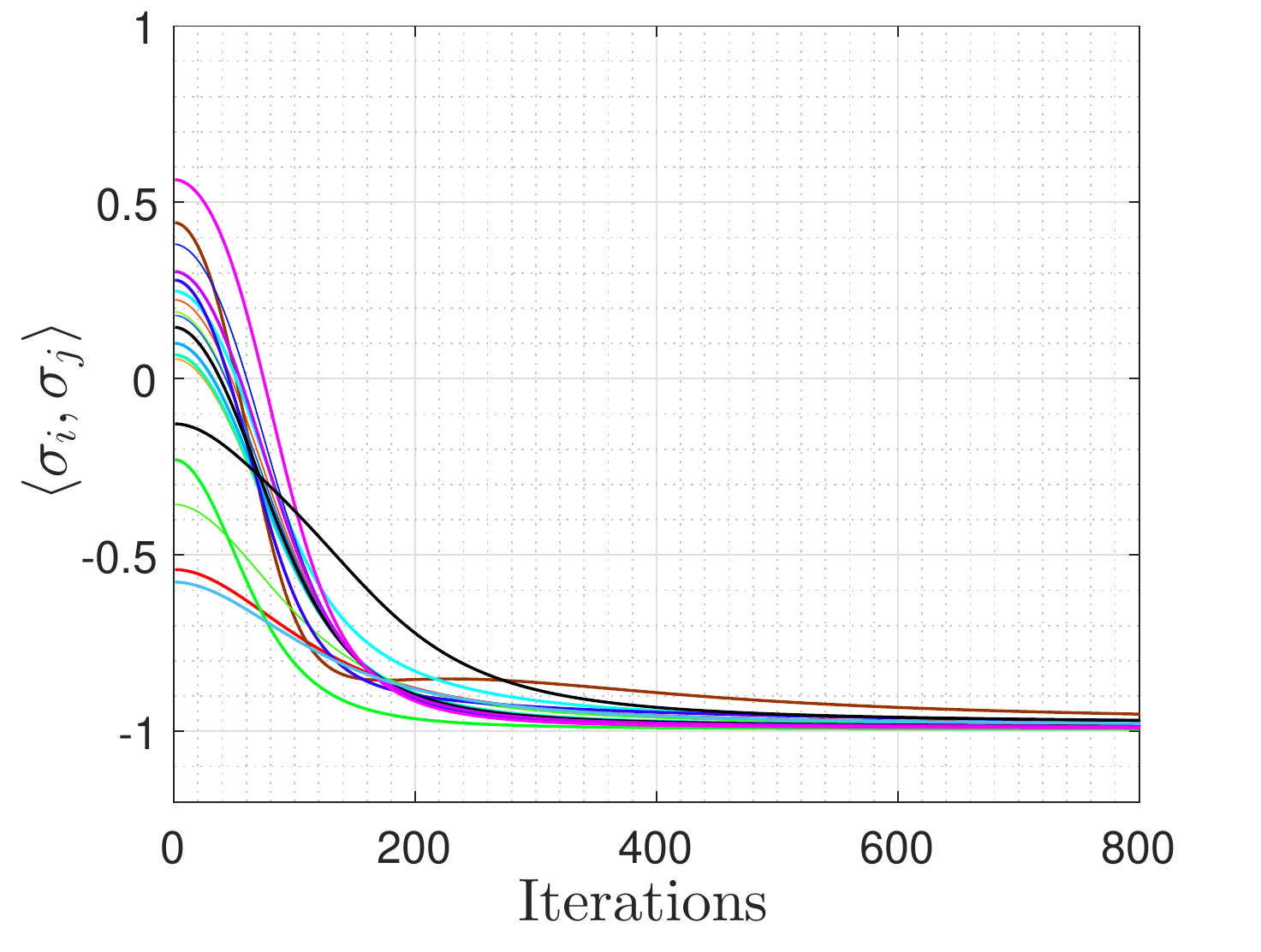}
     }
\caption{Panel (a): The decay of errors $\Delta_{1}(\bm{\sigma})$ and $\Delta_{2}(\bm{\sigma},\bm{\lambda})$ as defined in eqs. \eqref{Eq:stopping_criterion1} and \eqref{Eq:stopping_criterion}, resp. Panel (b): The inner product $\langle\bm{\sigma}_{i},\bm{\sigma}_{j} \rangle,(i,j)\in \tilde{E}$ for $\tilde{E}\subset E,|\tilde{E}|=20$ at different iterations with $\alpha=0.01$ and $\eta_{t}=\eta=1/100$.}
\label{Fig:MAXCUT}
\end{figure}

Next, we consider the constrained settings, where the elements of $E\subset V\times V$ are sampled randomly to obtain a subgraph $\tilde{E}\subset E$ with the cardinality $|\tilde{E}|=20$. Figure \ref{Fig:MAXCUT} denotes the the inner product $\langle\bm{\sigma}_{i},\bm{\sigma}_{j}\rangle$, for all $(i,j)\in \tilde{E}$. Evidently, the inner products $\langle\bm{\sigma}_{i},\bm{\sigma}_{j} \rangle,(i,j)\in \tilde{E}$ satisfy the constraints \eqref{Eq:doublecone-2} asymptotically. We also demonstrated the errors in the constraint violation and the norm of the gradient as defined in eqs. \eqref{Eq:stopping_criterion1} and \eqref{Eq:stopping_criterion}.

\section{Proofs} 
\label{sec:Proofs}
In this section, we provide the proofs of main theorems and their corollaries, while deferring the more technical proofs to the appendices. In particular, we first prove the convergence results for Algorithm \ref{alg:1} when the sectional curvature is non-positive $\kappa \leq K(X_{p},Y_{p})\leq 0, X_{p},Y_{p}\in T_{p}\mathcal{M}$ for all $p\in \mathcal{M}$. We then prove a different convergence results for positive curvature case $K(X_{p},Y_{p})>\kappa\geq 0$ and asymptotically non-negative case $K(X_{p},Y_{p})\geq -\kappa/d^{2}(p,q)$ for a fixed $q\in \mathcal{M}$. We also prove a high probability convergence bound under certain technical assumptions on the norm of the estimation error $\|e_{t}\|$.

\subsection{Proof of Theorem \ref{Thm:1} (Hyperbolic Manifolds)}
\label{app:Proof_of_Theorem}
We start by writing the Taylor expansion of $\ell\circ \Exp_{p}$ at the point $p\in\mathcal{M}$ for some smooth generic function $\ell\in C^{2}(\mathcal{M},\real)$. Note that $\ell\circ \Exp_{p}$ is a function defined on the tanget space $T_{p}\mathcal{M}\cong \real^{n}$. Furthermore, we note that $\dfrac{d}{d\epsilon}\Big \vert_{\epsilon=0}\ell(\Exp_{p}(\epsilon X))= \nabla_{X}\ell(p)=\langle \grad_{p}\ell,X\rangle$ for any $X_{p}\in T_{p}\mathcal{M}$. Thus, based on the first order Taylor expansion with the integral remainder we derive (cf. \cite{smith1994optimization})
\begin{align}
\label{Eq:expansion}
\ell(\Exp_{p}(\epsilon X_{p}))=\ell(p)+\epsilon\langle \grad_{p}\ell,X_{p}\rangle+\epsilon^{2}\int_{0}^{1}(1-s)\nabla_{X_{p}}^{2}\ell(\Exp_{p}s\epsilon X_{p})\ \D s,
\end{align}
where $\nabla_{X_{p}}^{2}\ell\topdoteq (\nabla\ell(X_{p},X_{p}))_{p}$.

Now, we identify the parameters $p=x_{t}$, $X_{p}=X_{t}=-\grad_{x_{t}}\mathfrak{L}(x_{t},\lambda_{t};\xi_{t})\in T_{x_{t}}\mathcal{M}$, $\epsilon=\eta_{t}$, and $\ell(x)={1\over 2} d^{2}(x,x_{\ast})$. Recall from (A2) in Section \ref{subsec:Hyperbolic}, the squared distance function is differentiable, and $\grad_{x} {1\over 2} d^{2}(x,x_{\ast})=-\log_{x}(x_{\ast})$, where the Riemannian logarithm exists due to (A3). Therefore, from the expansion in eq. \eqref{Eq:expansion} we compute
\begin{align}
\nonumber
{1\over 2} d^{2}(\Exp_{x_{t}}(\eta_{t} X_{t}),x_{\ast})=&{1\over 2} d^{2}(x_{t},x_{\ast})-\eta_{t}\langle \log_{x_t}(x_{\ast}) ,X_{t} \rangle\\
\label{Eq:Taylor_Rules}
&+{1\over 2} \eta_{t}^{2}\int_{0}^{1}(1-s)\nabla_{X_{t}}^{2}d^{2}(\Exp_{x_{t}}(s\eta_{t}X_{t}),x_{\ast})\ \D s.
\end{align}
Due to the update rule in Algorithm \ref{alg:1} we have $x_{t+1}=\Exp_{x_{t}}(\eta_{t} X_{t})$. Therefore,
\begin{align}
\nonumber
{1\over 2} d^{2}(x_{t+1},x_{\ast})=&{1\over 2} d^{2}(x_{t},x_{\ast})-\eta_{t}\langle \log_{x_{t}}(x_{\ast}) ,X_{t} \rangle \\
\label{Eq:last_term}
&+{1\over 2}\eta_{t}^{2}\int_{0}^{1}(1-s)\nabla_{X_{t}}^{2}d^{2}(\Exp_{x_t}(s\eta_{t}X_{t}),x_{\ast})\ \D s.
\end{align}
Now, to bound the last term of the expression in eq. \eqref{Eq:last_term}, in the following we consider the special case of the Hadamard manifolds, where the sectional curvature is negative everywhere on $\mathcal{M}$. 

In this case, the maximum eigenvalue $\lambda_{\max}$ of the Hessian of the squared distance function $\ell(x)={1\over 2}d^{2}(x,x_{\ast})$ is upper bounded by (cf. \cite{bonnabel2013stochastic},\cite{cordero2001riemannian})
\begin{align}
\label{Eq:inequality0}
\lambda_{\max}\left(\nabla^{2}{d^{2}(x,x_{\ast})\over 2}\right)&\leq \dfrac{\sqrt{|\kappa|{1\over 2}d^{2}(x,x_{\ast})}}{\tanh\left(\sqrt{|\kappa|{1\over 2}d^{2}(x,x_{\ast})}\right)}\\ \label{Eq:inequality1}
&\leq d(x,x_{\ast})\sqrt{{1\over 2}|\kappa|}+1,
\end{align}
where we recall that $\kappa<0$ is a lower bound on the sectional curvature. Therefore,
\begin{align}
\nabla_{X_{t}}^{2}d^{2}(\Exp_{x_t}(s\eta_{t}X_{t}),x_{\ast})\leq (d(\Exp_{x_t}(s\eta_{t}X_{t}),x_{\ast})\sqrt{|\kappa|/2}+1)\|X_{t}\|^{2}.
\end{align}
From eq. \eqref{Eq:last_term} and the inequality \eqref{Eq:inequality1} we derive
\begin{align}
\nonumber
{1\over 2} d^{2}(x_{t+1},x_{\ast})\leq &{1\over 2} d^{2}(x_{t},x_{\ast})-\eta_{t}\langle \log_{x_{t}}(x_{\ast}) ,X_{t} \rangle +{1\over 4}{\eta_{t}^{2}}\|X_{t}\|^{2} \\ \label{Eq:Comb1}
& +\dfrac{\sqrt{|\kappa|}\eta_{t}^{2}}{\sqrt{2}}\|X_{t}\|^{2} \int_{0}^{1}(1-s) d(\Exp_{x_{t}}(s\eta_{t}X_{t}),x_{\ast})\ \D s.
\end{align}
Applying the method of integration by parts to the last term gives us
\begin{align}
\nonumber
\int_{0}^{1}(1-s) d(\Exp_{x_{t}}(s\eta_{t}X_{t}),&x_{\ast})\ \D s =\left(s-{s^{2}\over 2}\right) d(\Exp_{x_{t}}(s\eta_{t}X_{t}),x_{\ast})\vert_{s=0}^{1}\\ \nonumber &-\int_{0}^{1}\left(s-{s^{2}\over 2}\right){\D \over \D s}  d(\Exp_{x_{t}}(s\eta_{t}X_{t}),x_{\ast})\ \D s\\ \nonumber
&\stackrel{\rm{(a)}}{=}\dfrac{1}{2}d(x_{t+1},x_{\ast})\\ \nonumber &-\eta_{t}\int_{0}^{1}\left(s-{s^{2}\over 2}\right)\dfrac{\langle P_{x_{t}\rightarrow \gamma(s)} X_{t},\log_{\Exp_{x_{t}}(s\eta_{t}X_{t})}(x_{\ast})\rangle}{d(\Exp_{x_{t}}(s\eta_{t}X_{t}),x_{\ast})}\ \D s\\  \nonumber
&\stackrel{\rm{(b)}}{\leq} \dfrac{1}{2}d(x_{t+1},x_{\ast})+\eta_{t}\|P_{x_{t}\rightarrow \gamma(s)}X_{t}\|\int_{0}^{1}\left(s-{s^{2}\over 2}\right)\ \D s,\\ \label{Eq:Comb2}
&\stackrel{\rm{(c)}}{=}\dfrac{1}{2}d(x_{t+1},x_{\ast})+\dfrac{1}{3}\eta_{t}\|X_{t}\|,
\end{align}
where to derive $\rm {(a)}$ we used the fact that $\gamma(s)\topdoteq \Exp_{x_{t}}(s\eta_{t}X_{t})$ is a geodesic with the initial point $\gamma(0)=x_{t}$, and the tangent $\dot{\gamma}(0)=\eta_{t}X_{t}$. Therefore, 
\begin{align*}
{\D \over \D s}\Exp_{x_{t}}(s\eta_{t}X_{t})=\eta_{t}P_{x_{t}\rightarrow \gamma(s)}X_{t},
\end{align*}
where $P_{x_{t}\rightarrow \gamma(s)}X_{t}$ is the parallel-transport of the vector $X_{t}$ from the tangent plan $T_{x_{t}}\mathcal{M}$ to the tangent plan $T_{\gamma(s)}\mathcal{M}$. To derive the step $\rm{(b)}$, we used the following inequality for the inner product in the integrand, 
\begin{align*}
\left\vert\dfrac{\langle P_{x_{t}\rightarrow \gamma(s)} X_{t},\log_{\Exp_{x_{t}}(s\eta_{t}X_{t})}(x_{\ast})\rangle}{\|P_{x_{t}\rightarrow \gamma(s)} X_{t}\|d(\Exp_{x_{t}}(s\eta_{t}X_{t}),x_{\ast})}\right|\leq 1.
\end{align*} 
Lastly, in step $\rm{(c)}$, we used the fact that a geodesic $\gamma(s)$ has a constant velocity and thus $\|P_{x_{t}\rightarrow \gamma(s)}X_{t}\|=\|X_{t}\|$.

Putting together the inequalities \eqref{Eq:Comb1} and \eqref{Eq:Comb2} yields
\begin{align}
\nonumber
\eta_{t}\langle \log_{x_{t}}(x_{\ast}) ,X_{t} \rangle\leq &{1\over 2} d^{2}(x_{t},x_{\ast})-{1\over 2} d^{2}(x_{t+1},x_{\ast})\\ 
\label{Eq:Before_Triangle} &+{\eta_{t}^{2}}\|X_{t}\|^{2}\left({1\over 4}+{\sqrt{|\kappa|}\over 2\sqrt{2}}d(x_{t+1},x_{\ast})\right)+{\sqrt{|\kappa|}\over 3\sqrt{2}}\eta_{t}^{3}\|X_{t}\|^{3}.
\end{align}
Since $(x_{t},\lambda_{t})\in \ball_{R}(x_{\ast})\times \real_{+}$ for all $t\in [T]$, we have $d(x_{t+1},x_{\ast})<R$ and from eq. \eqref{Eq:Before_Triangle} we derive
\begin{align}
\nonumber
\eta_{t}\langle \log_{x_{t}}(x_{\ast}) ,X_{t} \rangle\leq &{1\over 2} d^{2}(x_{t},x_{\ast})-{1\over 2} d^{2}(x_{t+1},x_{\ast})\\  \label{Eq:the_same_inequality}
 &+{\eta_{t}^{2}}\|X_{t}\|^{2}\left({1\over 4}+{\sqrt{|\kappa|}R \over 2\sqrt{2}}\right)+{\sqrt{|\kappa|}\over 3\sqrt{2}}\eta_{t}^{3}\|X_{t}\|^{3}.
\end{align}
By the geodesic convexity condition in Assumption \ref{Assumption:g-convex}, we also have the following inequality 
\begin{align}
\label{Eq:geodesic_convex_Ref}
\mathfrak{L}(x_{t},\lambda_{t})&-\mathfrak{L}(x_{\ast},\lambda_{t})\leq -\langle \log_{x_{t}}(x_{\ast}),\grad_{x_{t}}\mathfrak{L}(x_{t},\lambda_{t})\rangle.
\end{align}
Thus, from \eqref{Eq:the_same_inequality} and using the fact that $X_{t}=-\mathfrak{L}(x_{t},\lambda_{t};\xi_{t})=e_{t}-\mathfrak{L}(x_{t},\lambda_{t})$, we have after summation over $t=0,1,\cdots,T-1$ that
\begin{align}
\nonumber
\sum_{t=0}^{T-1}\eta_{t}(\mathfrak{L}(x_{t},\lambda_{t})-\mathfrak{L}(x_{\ast},\lambda_{t}))
&\leq {1\over 2} d^{2}(x_{0},x_{\ast})+\left({1\over 4}+{\sqrt{|\kappa|}R \over 2\sqrt{2}}\right)\sum_{t=0}^{T-1}{\eta_{t}^{2}}\|X_{t}\|^{2}\\ \label{Eq:similar_result} &+{\sqrt{|\kappa|}\over 3\sqrt{2}}\sum_{t=0}^{T-1}\eta_{t}^{3}\|X_{t}\|^{3}+\sum_{t=0}^{T-1}\eta_{t}\langle \log_{x_{t}}(x_{\ast}),e_{t} \rangle.
\end{align}

We now obtain a similar inequality for dual variables. For all $\lambda\in \real_{+}^{m}$, we have
\begin{align}
\nonumber
\|\lambda_{t+1}-\lambda\|^{2}&=\|\Pi_{\real^{m}_{+}}(\lambda_{t}+\eta_{t} \grad_{\lambda}\mathfrak{L}(x_{t},\lambda_{t})-\lambda\|^{2}
\\ \nonumber &\leq \|\lambda_{t}+\eta_{t} \grad_{\lambda_{t}}\mathfrak{L}(x_{t},\lambda_{t})-\lambda\|^{2}
\\ \label{Eq:Earlier_Result}  &\leq \|\lambda_{t}-\lambda\|^{2}+2\eta_{t}\langle \grad_{\lambda_{t}}\mathfrak{L}(x_{t},\lambda_{t}) ,\lambda_{t}-\lambda \rangle  +\eta_{t}^{2}\|\grad_{\lambda_{t}}\mathfrak{L}(x_{t},\lambda_{t})\|^{2}.
\end{align}
Summation over $t=0,1,\cdots,T-1$ and using the telescoping sum series gives us
\begin{align}
\label{Eq:Dual_Ineq} 
\|\lambda_{T+1}-\lambda\|^{2} \leq \|\lambda\|^{2}+\sum_{t=0}^{T-1}2\eta_{t}\langle \grad_{\lambda_{t}}\mathfrak{L}(x_{t},\lambda_{t}) ,\lambda_{t}-\lambda \rangle  +\sum_{t=0}^{T-1}\eta_{t}^{2}\|\grad_{\lambda_{t}}\mathfrak{L}(x_{t},\lambda_{t})\|^{2},
\end{align}
where we used the fact that $\lambda_{0}=0$. Since the left hand side of \eqref{Eq:Dual_Ineq} is positive, we can write
\begin{align}
\sum_{t=0}^{T-1}\eta_{t}\langle \grad_{\lambda_{t}}\mathfrak{L}(x_{t},\lambda_{t}) ,\lambda-\lambda_{t} \rangle  \leq {1\over 2}\|\lambda\|^{2}+{1\over 2}\sum_{t=0}^{T-1}\eta_{t}^{2}\|\grad_{\lambda_{t}}\mathfrak{L}(x_{t},\lambda_{t})\|^{2},
\end{align}
Due to the concavity of $\mathfrak{L}(x_{t},\cdot)$ we have the following inequality
\begin{align*}
\mathfrak{L}(x_{t},\lambda)-\mathfrak{L}(x_{t},\lambda_{t})\leq \langle \grad_{\lambda_{t}}\mathfrak{L}(x_{t},\lambda_{t}) ,\lambda-\lambda_{t} \rangle.
\end{align*}
Consequently,
\begin{align}
\label{Eq:sum2}
\sum_{t=0}^{T-1}\eta_{t}\left(\mathfrak{L}(x_{t},\lambda)-\mathfrak{L}(x_{t},\lambda_{t})\right)\leq {1\over 2} \|\lambda\|^{2}+{1\over 2} \sum_{t=0}^{T-1}\eta_{t}^{2}\|\grad_{\lambda}\mathfrak{L}(x_{t},\lambda_{t})\|^{2}.
\end{align}
Combining eqs. \eqref{Eq:similar_result} and \eqref{Eq:sum2} gives
\begin{align}
\nonumber
&\sum_{t=0}^{T-1}\eta_{t}(\mathfrak{L}(x_{t},\lambda)-\mathfrak{L}(x_{\ast},\lambda_{t}))\leq {1\over 2} \|\lambda\|^{2}+{1\over 2}d^{2}(x_{0},x_{\ast})+\sum_{t=0}^{T-1}\eta_{t}\langle \log_{x_{t}}(x_{\ast}),e_{t} \rangle\\ \label{Eq:Before_Using_Bounds}
&+{1\over 2} \sum_{t=0}^{T-1}\eta_{t}^{2}\|\grad_{\lambda}\mathfrak{L}(x_{t},\lambda_{t})\|^{2}+\left({1\over 4}+{\sqrt{|\kappa|}R \over 2\sqrt{2}}\right)\sum_{t=0}^{T-1}{\eta_{t}^{2}}\|X_{t}\|^{2}+{\sqrt{|\kappa|}\over 3\sqrt{2}}\sum_{t=0}^{T-1}\eta_{t}^{3}\|X_{t}\|^{3}.
\end{align}
To further proceed, we state a lemma:
\begin{lemma} 
\label{Lemma:3}
Consider Assumption \ref{Assumption:1}, and let $M\topdoteq\max\{M_{f},M_{g}\}$ and $\eta_{t}\leq 1/\alpha$ for all $t\in [T]$. Then, for all realizations of the stochastic process $\xi_{0},\cdots,\xi_{T-1}\in \Xi$, the following inequalities hold
\begin{align*}
\|\grad_{x_{t}}\mathfrak{L}(x_{t},\lambda_{t})\|&\leq M\left(1+\dfrac{\sqrt{m}G}{\alpha} \right),\\
\|\grad_{\lambda_{t}}\mathfrak{L}(x_{t},\lambda_{t})\|&\leq 2\sqrt{m}G,
\end{align*}
for all $t\in [T]$.
\end{lemma}

Now, recall that $X_{t}=e_{t}-\grad_{x_{t}}\mathfrak{L}(x_{t},\lambda_{t})$ and consider the following inequality
\begin{align}
\nonumber
\|X_{t}\|^{n}&\leq 2^{n-1}\|e_{t}\|^{n}+2^{n-1}\|\grad_{x_{t}}\mathfrak{L}(x_{t},\lambda_{t})\|^{n}\\ \label{Eq:comeback1}
&\leq 2^{n-1}\|e_{t}\|^{n}+2^{n-1} M\left(1+\dfrac{\sqrt{m}G}{\alpha} \right)^{n}, \quad n=2,3,
\end{align}
where the last inequality is due to the upper bound in Lemma \ref{Lemma:3}. Moreover,
\begin{align}
\label{Eq:comeback2}
\|\grad_{\lambda_{t}}\mathfrak{L}(x_{t},\lambda_{t})\|&\leq 2\sqrt{m}G.
\end{align}
Based on these two upper bounds, we obtain from \eqref{Eq:Before_Using_Bounds} that,
\begin{align}
\nonumber
&\sum_{t=0}^{T-1}\eta_{t}(\mathfrak{L}(x_{t},\lambda)-\mathfrak{L}(x_{\ast},\lambda_{t}))\leq {1\over 2}\|\lambda\|^{2}+
{1\over 2} d^{2}(x_{0},x_{\ast})+\sum_{t=0}^{T-1}\eta_{t}\langle \log_{x_{t}}(x_{\ast}),e_{t} \rangle
\\ \label{Eq:right_hand_side}
&+{1\over 2}\sum_{t=0}^{T-1}\Big(A+(1+{2\sqrt{|\kappa|/2}R})\|e_{t}\|^{2}\Big){\eta_{t}^{2}}+{4\sqrt{|\kappa|}\over 3\sqrt{2}}\sum_{t=0}^{T-1}(B+\|e_{t}\|^{3})\eta_{t}^{3},
\end{align}
where $A$ and $B$ are two non-negative constants defined in eqs. \eqref{Eq:Def_A} and \eqref{Eq:Def_B}, respectively.
We now expand the left hand side to obtain
\begin{align*}
&\sum_{t=0}^{T-1}\eta_{t}\left(f(x_{t})-f(x_{\ast})+\langle \lambda, h(x_{t}) \rangle -\langle \lambda_{t},h(x_{\ast}) \rangle+ \dfrac{\alpha}{2}(\|\lambda_{t}\|^{2}-\|\lambda\|^{2})\right)\\ &\hspace{70mm} \leq \text{r.h.s. of eq. \eqref{Eq:right_hand_side}}.
\end{align*}
Since $h(x_{\ast})\preceq 0$ and $\lambda_{t}\succeq  0$, and due to positivity of $\|\lambda_{t}\|^{2}$, we can remove these terms from the left hand side which then leaves us
\begin{align}
 \nonumber
&\sum_{t=0}^{T-1}\eta_{t}\left(f(x_{t})-f(x_{\ast})\right)+\left(\langle \lambda,\sum_{t=0}^{T-1}\eta_{t} h(x_{t}) \rangle -\dfrac{(\sum_{t=0}^{T-1}\eta_{t})\alpha+1}{2}\|\lambda\|^{2}\right)\\ \nonumber
\leq &{1\over 2} d^{2}(x_{0},x_{\ast})+\sum_{t=0}^{T-1}\eta_{t}\langle \log_{x_{t}}(x_{\ast}),e_{t} \rangle
+{1\over 2}\sum_{t=0}^{T-1}\Big(A+(1+{2\sqrt{|\kappa|/2}R})\|e_{t}\|^{2}\Big){\eta_{t}^{2}}\\ \label{Eq:right_hand_side_2} &+{4\sqrt{|\kappa|}\over 3\sqrt{2}}\sum_{t=0}^{T-1}(B+\|e_{t}\|^{3})\eta_{t}^{3}.
\end{align}
Maximizing the terms inside the parenthesis in the l.h.s. with respect to the vector $\lambda \in \real_{+}^{m}$ gives us
\begin{align}
\label{Eq:invoke}
&\sum_{t=0}^{T-1}\eta_{t}\left(f(x_{t})-f(x_{\ast})\right)+\dfrac{2}{(\sum_{t=0}^{T-1}\eta_{t})\alpha+1}\Big\| \Big[\sum_{t=0}^{T-1}\eta_{t}h(x_{t}) \Big]_{+}\Big\|^{2}
\leq \text{r.h.s. of eq. \eqref{Eq:right_hand_side_2}.}
\end{align}
We eliminate the second term on the l.h.s. since it is positive. We then take the expectation from both sides,
\begin{align}
\nonumber
\sum_{t=0}^{T-1}\eta_{t}\left(\expect[f(x_{t})]-f(x_{\ast})\right)\leq &\dfrac{1}{2}R^{2}+{1\over 2}\sum_{t=0}^{T-1}\Big(A+(1+{2\sqrt{|\kappa|/2}R})\expect[\|e_{t}\|^{2}]\Big){\eta_{t}^{2}}\\ \label{Eq:combine} &+{4 \sqrt{|\kappa|}\over 3\sqrt{2}}\sum_{t=0}^{T-1}(B+\expect[\|e_{t}\|^{3}])\eta_{t}^{3}.
\end{align}
To derive this inequality, we used the fact that $x_{t}$ is $\mathfrak{F}_{t}$-measurable since it is computed from all the random variables $\xi_{0},\xi_{1},\cdots,\xi_{t-1}$. Hence,
\begin{align}
\nonumber
\expect[\langle \log_{x_{t}}(x_{\ast}),e_{t} \rangle]&=\expect[\expect[\langle \log_{x_{t}}(x_{\ast}), e_{t} \rangle|\mathfrak{F}_{t}]]\\ &=\expect[\langle\log_{x_{t}}(x_{\ast}),\expect[\mathfrak{L}(x_{t},\lambda_{t})-\mathfrak{L}(x_{t},\lambda_{t};\xi_{t})|\mathfrak{F}_{t}] \rangle]\\ \label{Eq:Inequality_Martingale}
&=0.
\end{align}
Combining eq. \eqref{Eq:combine} with the inequality
\begin{align*}
\Big(\min_{t\in [T]}\expect[f(x_{t})]-f(x_{\ast})\Big)\sum_{t=0}^{T-1}\eta_{t} \leq \sum_{t=0}^{T-1}\eta_{t}(\expect[f(x_{t})]-f(x_{\ast})),
\end{align*}
yields the inequality of Theorem \ref{Thm:1}.

\subsubsection{Proof of Corollary \ref{Cor:2}}

Given the step size $\eta_{t}=1/\sqrt{t+1}$ and based on the Riemann sum approximation of integrals, we compute the following upper bounds 
\begin{subequations}
\begin{align}
\label{Eq:Stepsize_1}
& \sum_{t=0}^{T-1}\eta_{t}= \sum_{t=0}^{T-1}\dfrac{1}{\sqrt{t+1}}\leq 1+ \int_{0}^{T-1}\hspace{-2mm}\dfrac{1}{\sqrt{t+1}}dt= 2\sqrt{T}-1,\\
&\sum_{t=0}^{T-1}\eta_{t}^{2}= \sum_{t=0}^{T-1}\dfrac{1}{t+1}\leq 1+\int_{0}^{T-1}\hspace{-2mm}\dfrac{1}{t+1}dt=1+\log(T),\\
\label{Eq:Stepsize_4}
&\sum_{t=0}^{T-1}\eta_{t}^{3}= \sum_{t=0}^{T-1}\dfrac{1}{(t+1)^{3\over 2}}\leq 1+\int_{0}^{T-1}\hspace{-2mm}\dfrac{1}{(t+1)^{3\over 2}}dt=3-{2\over \sqrt{T}}.
\end{align}
\end{subequations}
Using the same technique, we can also prove a lower bound $\sum_{t=0}^{T-1}\eta_{t}\geq 2(\sqrt{T}-1)$. From these upper and lower bounds and based on the inequality of Theorem \ref{Thm:1}, we derive
\begin{align}
\min_{t\in [T]}\expect[f(x_{t})]-f(x_{\ast})=\mathcal{O}\left(\dfrac{\log(T)}{\alpha^{2}(\sqrt{T}-1)}+\dfrac{1}{\alpha^{3}(\sqrt{T}-1)} \right).
\end{align}
Now, we recall our previous results in eq. \eqref{Eq:invoke}, 
\begin{align}
\nonumber
\min_{t\in [T]}\left(f(x_{t})-f(x_{\ast})\right)&+\dfrac{2}{(\sum_{t=0}^{T-1}\eta_{t})\alpha+1}\Big\| \Big[\sum_{t=0}^{T-1}\eta_{t}h(x_{t}) \Big]_{+}\Big\|_{2}^{2}
\\ &\hspace{-6mm} \leq \dfrac{1}{2\sum_{t=0}^{T-1}\eta_{t}} \left(d^{2}(x_{0},x_{\ast})+ A\sum_{t=0}^{T-1}\eta_{t}^{2}+\dfrac{8 \sqrt{|\kappa|}}{3}B\sum_{t=0}^{T-1}\eta_{t}^{3}\right),
\end{align}
where we set $e_{t}=0$ for simplicity. Let the constant $\mathcal{F}\geq 0$ be such that 
\begin{align*}
\min_{t\in [T]}\left(f(x_{t})-f(x_{\ast})\right)=-\mathcal{F}\leq 0. 
\end{align*}
From the bounds \eqref{Eq:Stepsize_1}-
\eqref{Eq:Stepsize_4} we derive
\begin{align}
\nonumber
&\Big\| \Big[\sum_{t=0}^{T-1}\eta_{t}h(x_{t}) \Big]_{+}\Big\|_{2}^{2}\leq \dfrac{(2\sqrt{T}-1)\alpha+1}{2}\mathcal{F}\\ \nonumber
&+ \dfrac{(2\sqrt{T}-1)\alpha+1}{8(\sqrt{T}-1)}\left(d^{2}(x_{0},x_{\ast})+ A(1+\log(T))+\dfrac{8 \sqrt{|\kappa|}}{3}B(3-2/\sqrt{T})\right).
\end{align}
Now, we divide both sides by $1/\sum_{t=0}^{T-1}\eta_{t}\leq 1/(2\sqrt{T}-2)$ and use the definition $\widehat{\eta}_{t}\topdoteq \eta_{t}/\sum_{t=0}^{T-1}\eta_{t}$ to derive
\begin{align}
\label{Eq:Therefore}
&\Big\| \Big[{\sum_{t=0}^{T-1}\widehat{\eta}_{t}h(x_{t})} \Big]_{+}\Big\|_{2}^{2}\leq \dfrac{(2\sqrt{T}-1)\alpha+1}{4(\sqrt{T}-1)}\mathcal{F}\\ \nonumber  &+\dfrac{(2\sqrt{T}-1)\alpha+1}{8(\sqrt{T}-1)^{2}}\left(R^{2}+ A(1+\log(T))+\dfrac{8\sqrt{|\kappa|}}{3}B(3-2/\sqrt{T})\right),
\end{align}
where we used the fact that $d^{2}(x_{0},x_{\ast})<R^{2}$, and also used Jensen's inequality for $2$-norm, \textit{i.e.},
\begin{align*}
\Big\| \Big[{\sum_{t=0}^{T-1}\widehat{\eta}_{t}h(x_{t})}\Big]_{+}\Big\|_{2}^{2}\leq \dfrac{1}{\sum_{t=0}^{T-1}\eta_{t}}\Big\|\Big[\sum_{t=0}^{T-1}\eta_{t}h(x_{t}) \Big]_{+}\Big\|_{2}^{2}.
\end{align*}
Therefore, from \eqref{Eq:Therefore}, we obtain the following asymptotic result
\begin{align*}
\Big\| \Big[{\sum_{t=0}^{T-1}\widehat{\eta}_{t}h(x_{t})} \Big]_{+}\Big\|_{2}^{2}=\mathcal{O}\left(\alpha\right).
\end{align*}

\subsection{Proof of Theorem \ref{Thm:2} (Elliptic Manifold)}
\label{subsection:proof_of_Thm2}

To prove the theorem, we first state a proposition: 

\begin{proposition}
\label{Prop:Triangle Inequality on Elliptic Manifolds}
\textsc{\small{(Law of Cosine on Elliptic Manifolds)}}
Let $\{x_{t},\lambda_{t}\}_{t=1}^{T}$ be the sequence of the primal-dual points generated by Algorithm \ref{alg:1} with the step size $\eta_{t}\leq  \min\{1/\alpha, i(\mathcal{M})/\|X_{t}\|\}$. Then, for all $t=1,2,\cdots,T$, we have
\begin{align*}
\nonumber
\cos(\sqrt{\kappa}d(x_{t+1},x_{\ast}))\geq &\cos(\sqrt{\kappa}d(x_{t},x_{\ast}))\cos(\sqrt{\kappa}d(x_{t},x_{t+1}))\\
&+\sin(\sqrt{\kappa}d(x_{t},x_{\ast})) \sin(\sqrt{\kappa}d(x_{t},x_{t+1}))\cos\vartheta_{t},
\end{align*}
where $\vartheta_{t}$ is the angle between the arms $(x_{t},x_{\ast})$ and $(x_{t},x_{t+1})$.
\end{proposition}

The proof of Proposition \ref{Prop:Triangle Inequality on Elliptic Manifolds} is based on the comparative geometry techniques and is deferred to Appendix \ref{App:Proof_of_Prop_Triangle}.

Now, recall that $x_{t}\not\in \mathcal{C}(x_{\ast})$, \textit{i.e.}, $x_{t}$ is neither a conjugate or a cut point for $x_{\ast}$. Hence, $\log_{x_{t}}(x_{\ast})$ exists and the cosine similarity is well-defined, 
\begin{align}
\label{Eq:Cautions}
\cos \vartheta_{t}=\dfrac{\langle X_{t},\log_{x_{t}}(x_{\ast})\rangle}{\|X_{t}\|\cdot d(x_{t},x_{\ast})},
\end{align}
where we recall the definition $X_{t}=-\grad_{x_{t}} \mathfrak{L}(x_{t},\lambda_{t};\xi_{t})$. From Proposition \ref{Prop:Triangle Inequality on Elliptic Manifolds}, we thus obtain,
\begin{align*}
\cos(\sqrt{\kappa}d(x_{t+1},x_{\ast}))\geq &\cos(\sqrt{\kappa}d(x_{t},x_{\ast}))\cos(\eta_{t}\sqrt{\kappa}\|X_{t}\|)\\ &+\sin(\sqrt{\kappa}d(x_{t},x_{\ast}))\sin(\eta_{t}\sqrt{\kappa}\|X_{t}\|)\dfrac{\langle X_{t},\log_{x_{t}}(x_{\ast})\rangle}{\|X_{t}\|\cdot d(x_{t},x_{\ast})}.
\end{align*}
In the next step, we write the first order Taylor's expansion with the remainder term,
\begin{align*}
\cos(\eta_{t}\sqrt{\kappa}\|X_{t}\|)&=1-\dfrac{1}{2}\eta^{2}_{t}\kappa\|X_{t}\|^{2}\cos(\sigma)  \\
\sin(\eta_{t}\sqrt{\kappa}\|X_{t}\|)&=\eta_{t}\sqrt{\kappa}\|X_{t}\|-\dfrac{1}{2}\eta^{2}_{t}\kappa\|X_{t}\|^{2}\sin(\tilde{\sigma}), 
\end{align*}
where $\sigma,\tilde{\sigma}\in (0,\eta_{t}\sqrt{\kappa}\|X_{t}\|)$. Hence,
\begin{align*}
\cos(\sqrt{\kappa}d(x_{t+1},x_{\ast}))\geq  &\cos(\sqrt{\kappa}d(x_{t},x_{\ast}))+\eta_{t}\kappa \langle X_{t},\log_{x_{t}}(x_{\ast}) \rangle\cdot \sinc(\sqrt{\kappa} d(x_{t},x_{\ast}))\\ &-\dfrac{1}{2}\eta^{2}_{t}\kappa\|X_{t}\|^{2}\Delta,
\end{align*}
where $\Delta\topdoteq\left(\sin(\tilde{\sigma})\sin(\sqrt{\kappa}d(x_{t},x_{\ast}))\cos(\vartheta)+\cos(\sigma)\cos(\sqrt{\kappa}d(x_{t},x_{\ast}))\right)$, and $\sinc(x)=\sin(x)/x$. Since $|\Delta|<2$, we then have
\begin{align*}
&\eta_{t}\kappa \langle X_{t},\log_{x_{t}}(x_{\ast}) \rangle\cdot\sinc(\sqrt{\kappa} d(x_{t},x_{\ast}))\\ &\hspace{10mm} \leq \cos(\sqrt{\kappa}d(x_{t+1},x_{\ast}))-\cos(\sqrt{\kappa}d(x_{t},x_{\ast}))+\eta^{2}_{t}\kappa\|X_{t}\|^{2}.
\end{align*}
From the definition $X_{t}=e_{t}-\grad_{x_{t}}\mathfrak{L}(x_{t},\lambda_{t})$, we obtain
\begin{align}
\nonumber
-\eta_{t}\kappa \langle \grad_{x_{t}}\mathfrak{L}&(x_{t},\lambda_{t}),\log_{x_{t}}(x_{\ast}) \rangle\cdot\sinc(\sqrt{\kappa} d(x_{t},x_{\ast}))\\ \nonumber &\leq \cos(\sqrt{\kappa}d(x_{t+1},x_{\ast}))-\cos(\sqrt{\kappa}d(x_{t},x_{\ast}))+\eta^{2}_{t}\kappa\|X_{t}\|^{2}\\  \label{Eq:Right_hand_side223} &\hspace{4mm} -\eta_{t}\kappa \langle e_{t},\log_{x_{t}}(x_{\ast}) \rangle\cdot\sinc(\sqrt{\kappa} d(x_{t},x_{\ast})).
\end{align}
Now, from the geodesic convexity inequality in eq. \eqref{Eq:geodesic_convex_Ref} we have
\begin{align*}
\eta_{t}\kappa \left(\mathfrak{L}(x_{t},\lambda_{t})-\mathfrak{L}(x_{\ast},\lambda_{t}) \right) \rangle\cdot\sinc(\sqrt{\kappa} d(x_{t},x_{\ast}))\leq  \text{r.h.s. of  eq.}\ \eqref{Eq:Right_hand_side223}.
\end{align*}
By expanding the left hand side we compute
\begin{align}
\label{Eq:drop}
\eta_{t}\kappa \left(f(x_{t})-f(x_{\ast})+\langle\lambda_{t},h(x_{t})\rangle \right)\cdot\sinc(\sqrt{\kappa} d(x_{t},x_{\ast}))\leq  \text{r.h.s. of  eq.}\ \eqref{Eq:Right_hand_side223},
\end{align}
where we removed the term $-\langle \lambda_{t},h(x_{\ast}) \rangle\cdot \sinc(\sqrt{\kappa} d(x_{t},x_{\ast}))$ in the left hand side since $0\preceq \lambda_{t}$, $h(x_{\ast})\preceq 0$, and $\sinc(\sqrt{\kappa} d(x_{t},x_{\ast}))\geq 0$ due to the fact that $d(x_{t},x_{\ast})< \pi/\sqrt{\kappa}$ (cf. Thm. \ref{Thm:Meyer} in Appendix \ref{App:Proof_of_Prop_Triangle}).

We now proceed by analyzing the following two cases:
\vspace{4mm}
\begin{itemize}[leftmargin=*]
\item  \textit{Case of} $\langle \lambda_{t}, h(x_{t})\rangle \geq 0$:
 
\vspace{3mm}
We noted earlier $\sinc(\sqrt{\kappa}d(x_{t},x_{\ast}))\geq 0$. Therefore, $\langle \lambda_{t},h(x_{t})\rangle\sinc(\sqrt{\kappa} d(x_{t},x_{\ast}))\geq 0$ in this case and we can drop this term from the left hand side of \eqref{Eq:drop}. After taking the expectation and dividing both sides by $\kappa$ we obtain
\begin{align}
\nonumber
\eta_{t} \expect &\Big[\Big(f(x_{t})-f(x_{\ast})\Big)\sinc(\sqrt{\kappa}d(x_{t},x_{\ast}))\Big]\\ \label{Eq:Tight1}
&\leq \dfrac{1}{\kappa}\expect\Big[\left( \cos(\sqrt{\kappa}d(x_{t+1},x_{\ast}))-\cos(\sqrt{\kappa}d(x_{t},x_{\ast}))\right)\Big]+\eta^{2}_{t}\expect[\|X_{t}\|^{2}].
\end{align}

\item \textit{Case of}  $\langle \lambda_{t}, h(x_{t})\rangle < 0$: 
\vspace{3mm}

To analyze this case, we begin from our earlier result in eq. \eqref{Eq:sum2} for the dual variables which is also valid in the elliptic case,
\begin{align}
\nonumber
\|\lambda_{t+1}-\lambda\|^{2}\leq \|\lambda_{t}-\lambda\|^{2}+2\eta_{t}\langle \grad_{\lambda_{t}}\mathfrak{L}(x_{t},\lambda_{t}) ,\lambda_{t}-\lambda \rangle  +\eta_{t}^{2}\|\grad_{\lambda_{t}}\mathfrak{L}(x_{t},\lambda_{t})\|^{2}.
\end{align}
Therefore, after using the concavity of $\mathfrak{L}(x_{t},\cdot)$, we derive
\begin{align}
2\eta_{t}(\mathfrak{L}(x_{t},\lambda)-\mathfrak{L}(x_{t},\lambda_{t}))\leq \|\lambda_{t}-\lambda\|^{2}-\|\lambda_{t+1}-\lambda\|^{2}+\eta_{t}^{2}\|\grad_{\lambda_{t}}\mathfrak{L}(x_{t},\lambda_{t})\|^{2}.
\end{align}
Expanding the left hand side with $\lambda=0$ yields
\begin{align*}
2\eta_{t}\left(-\langle \lambda_{t},h(x_{t})\rangle+\dfrac{\alpha}{2}\|\lambda_{t}\|^{2}\right)\leq \|\lambda_{t}\|^{2}-\|\lambda_{t+1}\|^{2}+\eta_{t}^{2}\|\grad_{\lambda_{t}}\mathfrak{L}(x_{t},\lambda_{t})\|^{2}.
\end{align*}
Hence, by dropping the second term on the left hand side, we get
\begin{align}
\label{Eq:ineqality1}
-\eta_{t}\langle \lambda_{t},h(x_{t})\rangle \leq \dfrac{1}{2}\|\lambda_{t}\|^{2}-\dfrac{1}{2} \|\lambda_{t+1}\|^{2}+\dfrac{1}{2}\eta_{t}^{2}\|\grad_{\lambda_{t}}\mathfrak{L}(x_{t},\lambda_{t})\|^{2}.
\end{align}

Combining the inequalities \eqref{Eq:ineqality1} and \eqref{Eq:drop} results in
\begin{align*}
&\eta_{t}\kappa\Big(f(x_{t})-f(x_{\ast})\Big)\sinc(\sqrt{\kappa} d(x_{t},x_{\ast}))+\eta_{t}\langle\lambda_{t},h(x_{t})\rangle\cdot (\sinc(\sqrt{\kappa} d(x_{t},x_{\ast}))-1) \\
&\leq \cos(\sqrt{\kappa}d(x_{t+1},x_{\ast}))-\cos(\sqrt{\kappa}d(x_{t},x_{\ast}))+\eta^{2}_{t}\kappa\|X_{t}\|^{2}+\dfrac{1}{2}\|\lambda_{t}\|^{2}-\dfrac{1}{2} \|\lambda_{t+1}\|^{2}\\  &\hspace{4mm} -\eta_{t}\kappa \langle e_{t},\log_{x_{t}}(x_{\ast}) \rangle\cdot\sinc(\sqrt{\kappa} d(x_{t},x_{\ast}))+\dfrac{1}{2}\eta_{t}^{2}\|\grad_{\lambda_{t}}\mathfrak{L}(x_{t},\lambda_{t})\|^{2}.
\end{align*}
The second term on the left hand side is non-negative since $\langle h(x_{t}),\lambda_{t}\rangle < 0$ under the current case, and $\sinc(\sqrt{\kappa} d(x_{t},x_{\ast}))\leq 1$. Therefore, after eliminating this term and then taking the expectation, we arrive at
\begin{align}
\nonumber
\eta_{t}\kappa\expect[(f(x_{t})-f(x_{\ast}))&\sinc(\sqrt{\kappa} d(x_{t},x_{\ast}))]\\ \nonumber
&\leq \expect[\cos(\sqrt{\kappa}d(x_{t+1},x_{\ast}))-\cos(\sqrt{\kappa}d(x_{t},x_{\ast}))]+\eta^{2}_{t}\kappa\expect[\|X_{t}\|^{2}]\\  \label{Eq:strict_inequality}  &\hspace{4mm} + \dfrac{1}{2}\expect[\|\lambda_{t}\|^{2}]-\dfrac{1}{2}\expect[ \|\lambda_{t+1}\|^{2}]+\dfrac{1}{2}\eta_{t}^{2}\expect[\|\grad_{\lambda_{t}}\mathfrak{L}(x_{t},\lambda_{t})\|^{2}].
\end{align}
\end{itemize}

\begin{lemma}
\label{Eq:Proof_of_In}
Under the condition $\langle h(x_{t}),\lambda_{t} \rangle<0$, the following inequality holds
\begin{align}
\dfrac{1}{2}\|\lambda_{t}\|^{2}-\dfrac{1}{2} \|\lambda_{t+1}\|^{2}+\dfrac{1}{2}\eta_{t}^{2}\|\grad_{\lambda_{t}}\mathfrak{L}(x_{t},\lambda_{t})\|^{2}\geq 0.
\end{align}
\end{lemma}

From the preceding lemma, we see that the upper bound \eqref{Eq:strict_inequality} for the case $\langle h(x_{t}),\lambda_{t} \rangle<0$ is larger than the upper bound \eqref{Eq:Tight1} for the case $\langle h(x_{t}),\lambda_{t} \rangle\geq 0$. As a result, we write a unified upper bound in both cases
\begin{align}
\nonumber
&\eta_{t}\expect[(f(x_{t})-f(x_{\ast}))\sinc(\sqrt{\kappa} d(x_{t},x_{\ast}))]\\ \nonumber
&\leq \dfrac{1}{\kappa}\expect\left[\left( \cos(\sqrt{\kappa}d(x_{t+1},x_{\ast}))-\cos(\sqrt{\kappa}d(x_{t},x_{\ast}))\right)\right]+\dfrac{1}{2}\expect\left[\left(\|\lambda_{t}\|^{2}-\|\lambda_{t+1}\|^{2}\right)\right]\\
&\hspace{4mm} +\eta_{t}^{2}\expect[\|\grad_{x_{t}}\mathfrak{L}(x_{t},\lambda_{t})\|^{2}]+\eta_{t}^{2}\expect[\|\grad_{\lambda_{t}}\mathfrak{L}(x_{t},\lambda_{t})\|^{2}]+\eta_{t}^{2}\expect[\|e_{t}\|^{2}],
\end{align}
where here we used the fact that $\|X_{t}\|^{2}\leq 2\|e_{t}\|^{2}+2\|\grad_{x_{t}}\mathfrak{L}(x_{t},\lambda_{t})\|^{2}$. Now, we leverage the upper bounds on the gradients that we computed in Lemma \ref{Lemma:3}. After taking the sum over $t=0,1,\cdots,T-1$ we obtain
\begin{align}
\nonumber 
\sum_{t=0}^{T-1}\eta_{t}&\expect[(f(x_{t})-f(x_{\ast}))\sinc(\sqrt{\kappa} d(x_{t},x_{\ast}))]\\ \nonumber
&\stackrel{\rm{(d)}}{\leq} \dfrac{1}{\kappa}\left(\expect[\cos(\sqrt{\kappa}d(x_{T+1},x_{\ast}))]-\cos(\sqrt{\kappa}d(x_{0},x_{\ast}))\right)\\ \nonumber &\hspace{4mm} +\dfrac{1}{2}\|\lambda_{0}\|^{2}-\dfrac{1}{2}\expect[\|\lambda_{T}\|^{2}] +\sum_{t=0}^{T-1}(A+\expect[\|e_{t}\|^{2})\eta_{t}^{2}  \\ \nonumber
&\stackrel{\rm{(e)}}{\leq} \dfrac{1}{\kappa}\left(1-\cos(\sqrt{\kappa}d(x_{0},x_{\ast}))\right)+\sum_{t=0}^{T-1}(A+\expect[\|e_{t}\|^{2})\eta_{t}^{2}\\ \nonumber 
&\leq \dfrac{2}{\kappa}\sin^{2}\left(\dfrac{\sqrt{\kappa}}{2}d(x_{0},x_{\ast})\right)+\sum_{t=0}^{T-1}(A+\expect[\|e_{t}\|^{2})\eta_{t}^{2},
\end{align}
where in ${\rm(d)}$, $A$ is the constant we defined in eq. \eqref{Eq:Def_A}, and in ${\rm{(e)}}$ we used the fact that $\lambda_{0}=0$. Now, using the following inequality
\begin{align*}
\Big(\min_{t\in [T]}&\expect[(f(x_{t})-f(x_{\ast}))\sinc(\sqrt{\kappa} d(x_{t},x_{\ast}))]\Big)\sum_{t=0}^{T-1}\eta_{t}\\   &\hspace{10mm}\leq \sum_{t=0}^{T-1}\eta_{t}\expect[(f(x_{t})-f(x_{\ast}))\sinc(\sqrt{\kappa} d(x_{t},x_{\ast}))],
\end{align*}
completes the proof.

\subsection{Proof of Theorem \ref{Thm:3} (Elliptic Manifold-High Probability Bound)}
\label{subsec:Proof_of_Thm3}

Here, we use the concentration inequalities for random variables with bounded Orlicz norm to obtain a high probability bound for the convergence of the primal-dual algorithm, cf. \cite{duchi2012randomized},\cite{lan2012validation}. Throughout this subsection, recall that $\mathfrak{F}_{t}=\sigma(\xi_{0},\cdots,\xi_{t-1})$.

Following the analysis in Section \ref{subsection:proof_of_Thm2} without taking the expectation, we obtain 
\begin{align*}
\sum_{t=0}^{T-1}\eta_{t}(f(x_{t})-f(x_{\ast}))\sinc(\sqrt{\kappa} d(x_{t},x_{\ast}))&\leq \dfrac{2}{\kappa}\sin^{2}\left(\dfrac{\sqrt{\kappa}}{2}d(x_{0},x_{\ast})\right)+\sum_{t=0}^{T-1}(A+\|e_{t}\|^{2})\eta_{t}^{2}\\ \nonumber
&-\kappa\sum_{t=0}^{T-1}\eta_{t}\langle e_{t},\log_{x_{t}}(x_{\ast}) \rangle\cdot \sinc(\sqrt{\kappa}d(x_{t},x_{\ast})).
\end{align*}

Consider the term $-\kappa\sum_{t=0}^{T-1}\eta_{t}\langle e_{t},\log_{x_{t}}(x_{\ast}) \rangle\cdot \sinc(\sqrt{\kappa}d(x_{t},x_{\ast}))$. Due to Theorem \ref{Thm:Meyer} in Appendix \ref{App:Proof_of_Prop_Triangle}, we have $d(x_{t},x_{\ast})\leq \pi/\sqrt{\kappa}$ and hence $\sinc(\sqrt{\kappa}d(x_{t},x_{\ast}))\geq 0$. Therefore, \newline 
\begin{align}
\nonumber
|\langle e_{t},\log_{x_{t}}(x_{\ast}) \rangle\cdot \sinc(\sqrt{\kappa}d(x_{t},x_{\ast}))| &\leq  \|e_{t}\|\cdot \|\log_{x_{t}}(x_{\ast}) \| \cdot \sinc(\sqrt{\kappa}d(x_{t},x_{\ast})) \\ \nonumber
&= \|e_{t}\| d(x_{t},x_{\ast})\cdot \sinc(\sqrt{\kappa}d(x_{t},x_{\ast}))\\ \nonumber
&=\dfrac{1}{\sqrt{\kappa}} \|e_{t}\| \sin(\sqrt{\kappa}d(x_{t},x_{\ast}))\\ \label{Eq:Upper_Bound}
&\leq \dfrac{1}{\sqrt{\kappa}}\|e_{t}\|.
\end{align}
Due to a bounded Orlicz norm of error $\vertiii{e_{t}}_{\psi_{\nu}}=\beta<\infty$ given $\mathfrak{F}_{t}$, we have
\begin{align*}
&\expect\big[\psi_{\nu} (\sqrt{\kappa}\beta^{-1}|{\langle e_{t},\log_{x_{t}}(x_{\ast}) \rangle\cdot \sinc(\sqrt{\kappa}d(x_{t},x_{\ast}))}|)|\mathfrak{F}_{t} \big]\\
&\stackrel{\rm{(a)}}{\leq} \expect\big[\psi_{\nu}(\beta^{-1}\|e_{t}\|)|\mathfrak{F}_{t}\big]\leq 1,
\end{align*}
where $\rm{(a)}$ is due to the the upper bound \eqref{Eq:Upper_Bound} and the fact that the Young-Orlicz moduli $\psi$ is non-decreasing (cf. Definition \ref{Def:Orlicz}). Hence,
\begin{align*}
\|{\langle e_{t},\log_{x_{t}}(x_{\ast}) \rangle\cdot \sinc(\sqrt{\kappa}d(x_{t},x_{\ast}))})\|_{\psi_{\nu}}\leq \beta/\sqrt{\kappa}.
\end{align*}

We use Lemma \ref{Lemma:1} in Appendix \ref{sec:Heavy_Tail} with $Y_{t}\topdoteq \langle e_{t},\log_{x_{t}}(x_{\ast}) \rangle \sinc(\sqrt{\kappa}d(x_{t},x_{\ast}))$ and $\nu\geq 2$, where $\zeta=\beta/\sqrt{\kappa}$. Recall from eq. \eqref{Eq:Inequality_Martingale} that $\expect[Y_{t}|\mathfrak{F}_{t}]=0$. Then, for any $\delta>0$ we have
\begin{align*}
&\prob\left[\left|\sum_{t=0}^{T-1}\eta_{t}\langle e_{t},\log_{x_{t}}(x_{\ast}) \rangle\cdot \sinc(\sqrt{\kappa}d(x_{t},x_{\ast}))\right|\geq \delta \right] \\ &\hspace{30mm} \leq 2\exp\left(-c\cdot \max\left\{\dfrac{\delta^{2}\kappa}{\beta^{2}\|\eta\|_{2}^{2}}, \dfrac{\delta^{\nu}\kappa^{{\nu\over 2}}}{\beta^{\nu}\|\eta\|_{{\nu \over{\nu -1}}}^{\nu}}\right\}\right)\\
&\hspace{30mm}\leq 2\exp\left(-c\cdot \dfrac{\delta^{2}\kappa}{\beta^{2}\|\eta\|_{2}^{2}}\right),
\end{align*}
where $\eta\topdoteq (\eta_{0},\cdots,\eta_{T-1})$. Hence, with the probability of at least $1-\varrho$ we have
\begin{align}
\nonumber
-\kappa\sum_{t=0}^{T-1}\eta_{t}\langle e_{t},&\log_{x_{t}}(x_{\ast}) \rangle\cdot \sinc(\sqrt{\kappa}d(x_{t},x_{\ast}))\\&\hspace{10mm}\leq \nonumber
\kappa\bigg|\sum_{t=0}^{T-1}\eta_{t}\langle e_{t},\log_{x_{t}}(x_{\ast}) \rangle\cdot \sinc(\sqrt{\kappa}d(x_{t},x_{\ast}))\bigg|\\ \label{Eq:Ineq:11}
&\hspace{10mm} \leq {\beta \sqrt{\kappa} \|\eta\|_{2}}\left(\dfrac{1}{c}\log\left({2\over \varrho}\right) \right)^{1\over 2}.
\end{align}

Now, consider the term $\sum_{t=0}^{T-1}\eta_{t}^{2}\|e_{t}\|^{2}$. We define 
\begin{align*}
Z_{t}\topdoteq\|e_{t}\|^{2}-\expect[\|e_{t}\|^{2}|\mathfrak{F}_{t}].
\end{align*}
The following lemma is proved in Appendix \ref{Proof of Lemma:2}:
\begin{lemma}
\label{Lemma:2}
Let $\vertiii{e_{t}}_{\psi_{\nu}}=\beta$ for $0<\nu$. Then, $\|Z_{t}\|_{\psi_{{\nu\over 2 }}}\leq 3^{{2\over \nu}}2\beta^{2}$.
\end{lemma}

Using Lemma \ref{Lemma:2} along with Lemma \ref{Lemma:1} in Appendix \ref{sec:Heavy_Tail} provides the following tail bound for all $\nu\geq 2$,
\begin{align*}
\prob\Bigg[ \Bigg|\sum_{t=1}^{T-1}\eta_{t}^{2} &\|e_{t}\|^{2}-\sum_{t=1}^{T-1}\eta_{t}^{2}\expect[\|e_{t}\|^{2}|\mathfrak{F}_{t}]\Bigg|\geq \delta\Bigg] \\ &\leq 2\exp\left(-c\cdot \left(\dfrac{\delta^{2}}{3^{{4\over \nu}}4\beta^{4}\|\eta^{2}\|_{2}^{2}}\lor^{{\nu}} \dfrac{\delta^{\nu\over 2}}{3\cdot2^{\nu\over 2}\beta^{\nu}\|\eta^{2}\|_{\nu\over {\nu-2} }^{\nu\over 2}}\right)\right),
\end{align*}
where $\eta^{2}\topdoteq(\eta_{1}^{2},\eta_{2}^{2},\cdots,\eta_{T-1}^{2})$, and we define $x\lor^{\nu} y\topdoteq\min\{x,y\}$ if $2\leq \nu\leq 4$ and $x\lor^{\nu} y\topdoteq\max\{x,y\}$ if $\nu\geq 4$.
Moreover, due to the upper bound \eqref{Eq:Ineq:2} in Appendix \ref{Proof of Lemma:2} we have
\begin{align*}
\expect[\|e_{t}\|^{2}|\mathfrak{F}_{t}]\leq  \dfrac{4\beta^{2}}{\nu}\Gamma\left({2\over \nu}\right).
\end{align*}
Consequently, with the probability of at least $1-\varrho$ we derive the following inequality for $\nu\geq 2$,
\begin{align}
\nonumber
\sum_{t=1}^{T-1}\eta_{t}^{2} \|e_{t}\|^{2}
&\leq \dfrac{4\beta^{2}}{\nu}\Gamma\left({2\over \nu}\right)\sum_{t=1}^{T-1}\eta_{t}^{2}\\ \label{Eq:Ineq:22}
&\hspace{4mm}+2\beta^{2}3^{2\over \nu}\left(\|\eta^{2}\|_{2}\left(\dfrac{1}{c}\log\left({2\over \varrho}\right) \right)^{1\over 2}\land^{\nu} \|\eta^{2}\|_{\nu\over {\nu-2}}\left(\dfrac{1}{c}\log\left({2\over \varrho}\right) \right)^{2\over \nu}\right),
\end{align}
where we define $x\land^{\nu} y\topdoteq\max\{x,y\}$ if $2\leq \nu\leq 4$ and $x\land^{\nu} y\ \dot{=} \ \min\{x,y\}$ if $\nu\geq 4$.

Based on the inequalities \eqref{Eq:Ineq:11}-\eqref{Eq:Ineq:22} and the union bound, we derive that with the probability of at least $1-2\varrho$ the following inequality holds,
\begin{align}
\nonumber
-\kappa \sum_{t=0}^{T-1}\eta_{t}\langle e_{t}&,\log_{x_{t}}(x_{\ast})\rangle\cdot \sinc(\sqrt{\kappa}d(x_{t},x_{\ast}))+\sum_{t=1}^{T-1}2\eta_{t}^{2} \|e_{t}\|^{2}\leq \\ \nonumber
&+\dfrac{8\beta^{2}}{\nu}\Gamma\left({2\over \nu}\right)\sum_{t=1}^{T-1}\eta_{t}^{2}+{2\sqrt{\kappa}\beta\|\eta\|_{2}}\left(\dfrac{1}{c^{1\over 2}}\log\left({2\over \varrho}\right) \right)^{1\over 2}\\ \label{Eq:This_upper_bound}
&+4\beta^{2}3^{2\over \nu}\left(\|\eta^{2}\|_{2}\left(\dfrac{1}{c}\log\left({2\over \varrho}\right) \right)^{1\over 2}\land^{\nu} \|\eta^{2}\|_{\nu\over {\nu-2}}\left(\dfrac{1}{c}\log\left({2\over \varrho}\right) \right)^{2\over \nu}\right).
\end{align}
Using the upper bound \eqref{Eq:This_upper_bound} in conjunction with the proof of Section \ref{subsection:proof_of_Thm2} gives us the desired high probability convergence bound. Note that the convergence bounds we obtained in Theorem \ref{Thm:3} is a special cases of our general analysis here when $\nu=2$.

\subsection{Proof of Theorem \ref{Thm:5} (Asymptotically Elliptic)}
 
Here, we only present a sketch of the proof. The complete proof is similar to the proof of Theorem \ref{Thm:1} in Section \ref{subsec:Hyperbolic}. The only difference here is that a new upper bound on the Hessian of the quadratic distance function in \eqref{Eq:inequality0} is derived. The following lemma characterizes such an upper bound:
\begin{lemma}
\label{Lemma:5}
Consider the condition of Assumption \ref{Assumption:5} with the reference point $q=x_{\ast}$, \textit{i.e.},  $K(X,Y)\geq -\kappa/d^{2}(x_{t},x_{\ast})$ for all $x_{t}\in \mathcal{M}$ and all $X,Y\in T_{x_{t}}\mathcal{M}$. Then,
\begin{align*}
\nabla_{X}\left(\dfrac{d^{2}(x_{t},x_{\ast})}{2}\right)&\leq {1\over 2}(1+\sqrt{1+4\kappa^{2}}) \|X\|^{2}, \quad \text{for all}\ \ X\in T_{x_{t}}\mathcal{M}.
\end{align*} 
\end{lemma}
\begin{proof}
The proof is due to Villani \cite[pp. 276-279]{villani2008optimal}. For completeness, we present the proof here. Given the points $x_{t},x_{\ast}\in \mathcal{M}$, let $\gamma(t):[0,d(x_{t},x_{\ast})]\rightarrow \mathcal{M}$ be a minimizing geodesic that connects these two points, where $\gamma(0)=x_{t}$ and $\gamma(d(x_{t},x_{\ast}))=x_{\ast}$. Let $(e_{1}(t),e_{2}(t),\cdots,e_{n}(t))$ with $e_{1}(t)\topdoteq \dot{\gamma}(t)$ be an orthogonal basis for $T_{\gamma(t)}\mathcal{M}$. We note that $e_{1}(t)=\dot{\gamma}(t)$ is the eigenvector of the Hessian with the identity eigenvalue $\nabla_{e_{1}(t)}\left(\dfrac{d^{2}(x_{t},x_{\ast})}{2}\right)=\|e_{1}(t)\|^{2}$. Therefore, we focus on other eigenvalues in the tangent space $T_{\gamma(t)}\mathcal{M}$.

Let $\lambda_{i}(t)\topdoteq (\nabla^{2}_{e_{i}(t)}(d(x_{t},x_{\ast})))_{\gamma(t)}$ denotes the Hessian of the distance function (not quadratic). Then, the following inequality holds 
\begin{align}
\label{Eq:goal}
\dot{\lambda}_{i}(t)+\lambda_{i}^{2}(t)+K(e_{1}(t),e_{i}(t))\leq 0,
\end{align}
where $K$ is the sectional curvature. Now, we compute the Hessian of the map $x\mapsto d^{2}(x,x_{\ast})$,
\begin{align*}
\nabla^{2}\left(\dfrac{d^{2}(x,x_{\ast})}{2}\right)&=d(x,x_{\ast})\nabla^{2} d(x,x_{\ast})+\nabla d(x,x_{\ast})\otimes \nabla d(x,x_{\ast}).
\end{align*}
For $x=\gamma(t)$ we further have that
\begin{align*}
\nabla^{2}\left(\dfrac{d^{2}(\gamma(t),x_{\ast})}{2}\right)
&=d(\gamma(t),x_{\ast})\nabla^{2} d(\gamma(t),x_{\ast})+\dot{\gamma}(t)\otimes \dot{\gamma}(t)\\
&=t\nabla^{2} d(\gamma(t),x_{\ast})+e_{1}(t)\otimes e_{1}(t).
\end{align*}
Hence, for any $e_{i}(t),i=2,3,\cdots,n$ we have that
\begin{align}
\label{Eq:identity}
\delta_{i}(t)\topdoteq \nabla_{e_{i}(t)}^{2}\left(\dfrac{d^{2}(x_{t},x_{\ast})}{2}\right)&=t\nabla_{e_{i}(t)}^{2} d(\gamma(t),x_{\ast})+\langle e_{1}(t), e_{i}(t)\rangle^{2}=t\lambda_{i}(t).
\end{align}
From \eqref{Eq:goal} and the identity \eqref{Eq:identity}, we establish
\begin{align}
t\dot{\delta}_{i}(t)-\delta_{i}(t)+\delta_{i}^{2}(t)\leq -t^{2}K(e_{1}(t),e_{i}(t)), \quad i=2,3,\cdots,n.
\end{align}
Now, due to Assumption \ref{Assumption:5}, we have the inequality $K(e_{1}(t),e_{i}(t))\geq -\kappa/d(\gamma(t),x_{*})= -\kappa/t^{2}$. Therefore,
\begin{align}
t\dot{\delta}_{i}(t)\leq \delta_{i}(t)-\delta_{i}^{2}(t) +\kappa, \quad i=2,3,\cdots,n.
\end{align}
Now, if $\delta_{i}(t)\geq (1+\sqrt{1+4\kappa^{2}})/2$, the right hand side becomes negative and hence $\dot{\delta}_{i}(t)\leq 0$. We thus conclude that $\delta_{i}(t)\leq (1+\sqrt{1+4\kappa^{2}})/2$. 
\end{proof} 
From the Taylor series expansion in \eqref{Eq:Taylor_Rules} and Lemma \ref{Lemma:5}, we now obtain
\begin{align*}
\nonumber
{1\over 2} d^{2}(x_{t+1},x_{\ast})=&{1\over 2} d^{2}(x_{t},x_{\ast})+\eta_{t}\langle \log_{x_t}(x_{\ast}) ,X_{t} \rangle+{1\over 8}(1+\sqrt{1+4\kappa^{2}}) \eta_{t}^{2}\|X_{t}\|^{2}.
\end{align*}

The rest of the proof now is similar to Appendix \ref{app:Proof_of_Theorem}.

\section{Discussion and Final Remarks}
\label{Sec:discussion}

In this paper, we have developed a primal-dual method for optimizing functions over the Riemannian sub-manifolds. We analyzed the convergence of the proposed algorithm on manifolds with the positive (elliptic) and non-positive (hyperbolic) sectional curvature with a lower bound. We also analyzed the convergence rate of the algorithm on asymptotically elliptic manifolds, where the sectional curvature at each given point on the manifold is lower bounded by the inverse of the quadratic distance function. We outlined a few practical problems for which the primal-dual algorithm is applicable.

In the primal-dual algorithm we proposed, the minimizing geodesics must be computed at each algorithm iteration. This amounts to solving a system of coupled ordinary differential equations (cf. eq. \eqref{Eq:Computing_Geodesics}) which can be computationally expensive in the absence of symmetries to exploit. One remedy for this issue is to replace the exponential map in Algorithm \ref{alg:1} with a \textit{retraction} \cite{absil2009optimization}. Conceptually, a retraction can be thought of as an approximation of the exponential map that provides flexibility in terms of computation. Fortunately, every manifold that admits a Riemannian metric also admits a retraction defined by the exponential mapping. We note that the theoretical analysis we presented in this paper for the exponential map can be generalized to obtain similar results for the retractions.

In this work, we established a connection between the convergence rate of the primal-dual algorithm and the sectional curvature as a geometric invariant of smooth manifolds. As a future research, it is interesting to obtain performance guarantees in terms of the topological invariants such as the Betti or Euler numbers. In fact, the Gauss-Bonnet theorem and its higher dimensional variant, the Chern-Gauss-Bonnet theorem, establish a relation between the Riemannian curvature (a geometric invariant) and the Euler number (a topological invariant).  This relation as well as the geometric analysis of this paper can be used to obtain performance guarantees in terms of the Euler number.  

\appendix
\section{Proof of Lemma \ref{Lemma:3}}
\label{App:Proof_Lemma}
From the update rule of $\lambda_{t}$ in Algorithm \ref{alg:1} we have
\begin{align}
\nonumber
\|\lambda_{t+1}\|^{2}&= \left\|\Pi_{\real_{+}^{m}}\left(\lambda_{t}+\eta_{t}\grad_{\lambda_{t}}\mathfrak{L}(x_{t},\lambda_{t}) \right)\right\|^{2} \\ \nonumber
&\leq  \left\|\lambda_{t}+\eta_{t}\grad_{\lambda_{t}}\mathfrak{L}(x_{t},\lambda_{t})\right\|^{2}\\ \nonumber
&=\left\|\lambda_{t}+\eta_{t} h(x_{t})-\alpha\eta_{t}\lambda_{t}\right\|^{2}\\ \nonumber
&\stackrel{\rm{(a)}}{\leq} (1+\delta_{t})(1-\eta_{t}\alpha)^{2}\|\lambda_{t}\|^{2}+(1+\delta_{t}^{-1})\|\eta_{t} h(x_{t})\|^{2}\\ \label{Eq:Note01}
&\stackrel{\rm{(b)}}{\leq} (1+\delta_{t})(1-\eta_{t}\alpha)^{2}\|\lambda_{t}\|^{2}+(1+\delta_{t}^{-1})m\eta_{t}^{2}G^{2},
\end{align}
where $\rm{(a)}$ follows from Young's inequality with $\delta_{t}>0$, and the last step follows by Assumption \ref{Assumption:1}. Recall that $\eta_{t}\alpha\leq 1$ for all $t\in [T]$. Therefore, we can choose $\delta_{t}=\dfrac{\varepsilon_{t}}{(1-\eta_{t}\alpha)^{2}}-1$, where $\varepsilon_{t}\in ((1-\eta_{t}\alpha)^{2},1) $. This choice of $\delta_{t}$ results in $(1+\delta_{t})(1-\eta_{t}\alpha)^{2}=\varepsilon_{t}$ and we can proceed from eq. \eqref{Eq:Note01} as below
\begin{align}
\nonumber
\|\lambda_{t+1}\|^{2} &\leq \varepsilon_{t} \|\lambda_{t}\|^{2}+\dfrac{\varepsilon_{t}\eta_{t}^{2}}{\varepsilon_{t}-(1-\eta_{t}\alpha)^{2}}mG^{2}.
\end{align}
From this recursion with $\lambda_{0}=0$, we compute
\begin{align}
\nonumber
\|\lambda_{t+1}\|^{2}&\leq  mG^{2}\sum_{\ell=0}^{t}\dfrac{\eta_{\ell}^{2}\prod_{k=\ell}^{t}\varepsilon_{k}}{\varepsilon_{\ell}-(1-\eta_{\ell}\alpha)^{2}}  \\
 \ \nonumber
&\stackrel{\rm{(c)}}{=}\dfrac{mG^{2}}{\alpha^{2}}\sum_{\ell=0}^{t}\eta_{\ell}\alpha\prod_{k=\ell+1}^{t}(1-\eta_{k}\alpha)  \\ \label{Eq:Return_point}
&\stackrel{\rm{(d)}}{\leq} \dfrac{mG^{2}}{\alpha^{2}},
\end{align}
where $\rm{(c)}$ follows by choosing $\varepsilon_{k}=(1-\eta_{k}\alpha)$ for all $k=0,1,\cdots,t$, and $\rm{(d)}$ is due to the following inequality,
\begin{align}
\label{Eq:claim}
\sum_{\ell=0}^{t}\eta_{\ell}\alpha\prod_{k=\ell+1}^{t}(1-\eta_{k}\alpha)\leq 1,
\end{align}
which holds when $\eta_{t}\alpha\leq 1$ for all $t=0,1,\cdots,T-1$. For now, we proceed with the proof, returning to establish the inequality \eqref{Eq:claim} later.

From eq. \eqref{Eq:Return_point} and Assumption \ref{Assumption:1} with $M\topdoteq\max\{M_{f},M_{g}\}$, we now obtain
\begin{align}
\nonumber
\|\grad_{x_{t}}\mathfrak{L}(x_{t},\lambda_{t})\|&\leq \|\grad_{x_{t}}f(x_{t})\|+\|\lambda_{t}\|\cdot\|\grad_{x_{t}}h(x_{t})\|\\ \label{Eq:Lagrange_bound1}
&\leq M\left(1+\dfrac{\sqrt{m}G}{\alpha}\right).
\end{align}
Similarly,
\begin{align}
\nonumber
\|\grad_{\lambda_{t}}\mathfrak{L}(x_{t},\lambda_{t})\|&\leq \|h(x_{t})\|+\alpha\|\lambda_{t}\|\\ \nonumber
&\leq G\left(1+\sqrt{m}\right)\\
\label{Eq:Lagrange_bound2}
&\leq 2\sqrt{m}G.
\end{align}

We now prove the inequality \eqref{Eq:claim}. For ease of notation let $\theta_{t}\topdoteq \eta_{t}\alpha\in (0,1]$. First, suppose $\theta_{k}\not= 1$ for all $k=1,2,\cdots,t$. In this case, the sum of products is monotone increasing in $\theta_{0}$ and we have
\begin{align}
\nonumber
&\sum_{\ell=0}^{t}\theta_{\ell}\prod_{k=\ell+1}^{t}(1-\theta_{k})\\ \nonumber
&= \Big(\theta_{0}(1-\theta_{1})(1-\theta_{2})\cdots(1-\theta_{t})\Big)+\Big(\theta_{1}(1-\theta_{2})\cdots(1-\theta_{t})\Big)+\cdots+\theta_{t}\\
\nonumber
&\stackrel{(\rm{e})}{\leq} \Big((1-\theta_{1})(1-\theta_{2})\cdots(1-\theta_{t})\Big)+\Big(\theta_{1}(1-\theta_{2})\cdots(1-\theta_{t})\Big)+\cdots+\theta_{t}\\
\nonumber
&\stackrel{(\rm{f})}{=}\Big((1-\theta_{1}+\theta_{1})(1-\theta_{2})\cdots (1-\theta_{t})\Big)+\Big(\theta_{2}(1-\theta_{3})\cdots(1-\theta_{t})\Big)+\cdots+\theta_{t}\\
\nonumber
&\stackrel{(\rm{g})}{=}\big((1-\theta_{2}+\theta_{2})(1-\theta_{3})\cdots (1-\theta_{t})\big)+\cdots+\theta_{t}\\ \nonumber
\hspace{5mm}\vdots \\ \label{Eq:expansion_reduction}
&=(1-\theta_{t})+\theta_{t}=1,
\end{align}
where $(\rm{f})$ follows by combining the first and second parentheses in $(\rm{e})$, $(\rm{g})$ follows by combining the first and second parentheses in $(\rm{f})$, and so on.

Now, consider that $\theta_{k}=1$ for some indices $k\in \{1,2,\cdots,t\}$ and let $j$ be the largest index among those indices, \textit{i.e.}, $\theta_{j}=1$ and $\theta_{k}<1$ for all $k\in \{j+1,\cdots,t\}$. In this case, we simply have
\begin{align*}
\sum_{\ell=0}^{t}\theta_{\ell}\prod_{k=\ell+1}^{t}(1-\theta_{k})=\sum_{\ell=j}^{t}\theta_{\ell}\prod_{k=\ell+1}^{t}(1-\theta_{k})= 1,
\end{align*}
where the last equality can be proved using the same approach we used to derive \eqref{Eq:expansion_reduction}. 

\section{Proof of Proposition \ref{Prop:Triangle Inequality on Elliptic Manifolds}}
\label{App:Proof_of_Prop_Triangle}

First, we state a definition. 

\begin{definition}\textsc{(\cite[Def. 2.1]{meyer1989toponogov})}
A geodesic hinge $\gamma$,$\gamma_{0}$,$\vartheta$ in $\mathcal{M}$ consists of two non-constant geodesic
segments $\gamma$,$\gamma_{0}$ with the same initial point making the angle $\vartheta$. A minimal connection
$\gamma_{1}$ between the endpoints of $\gamma$ and $\gamma_{0}$ is called a closing edge of the hinge.
\end{definition}

In the following theorem $\mathcal{M}_{\kappa}$ denotes the model space that has a constant sectional curvature $\kappa$ everywhere on the manifold. 
\begin{theorem}\label{Thm:Toponogov}\textsc{({Toponogov's Comparison Theorem}\ \cite[Thm. 2.2]{meyer1989toponogov})} Let $\mathcal{M}$ be a complete Riemannian manifold with the sectional curvature $K\geq \kappa$.

Let $\gamma,\gamma_{0},\vartheta$ be a hinge in $\mathcal{M}$ with $\gamma_{0}$ minimal and $L[\gamma]\leq {\pi\over \sqrt{\kappa}}$ in case $\kappa>0$, and $\gamma_{1}$ a closing edge. Then, the closing edge of $\tilde{\gamma}_{1}$ of any hinge $\tilde{\gamma},\tilde{\gamma}_{0},\vartheta$ in $\mathcal{M}_{\kappa}$ with $L[\tilde{\gamma}]=L[\gamma]$, $L[\tilde{\gamma}_{0}]=L[\gamma_{0}]$ satisfies 
\begin{align}
\label{Eq:Inequality_Top}
L[\tilde{\gamma}_{1}]\geq L[\gamma_{1}].
\end{align}
\end{theorem}

Before we proceed with the proof, we state a few technical remarks about the minimality condition of the geodesics $\gamma_{0},\gamma_{1}$ in Theorem \ref{Thm:Toponogov}. When the sectional curvature lower bound is positive $\kappa>0$, the minimality condition of the geodesics $\gamma_{0},\gamma_{1}$ puts a restriction on their lengths. This is due to Myers' theorem that is stated below:
\begin{theorem}\textsc{(\cite[Thm. 19.4]{milnor2016morse})}
\label{Thm:Meyer}
Suppose the Ricci curvature satisfies
\begin{align*}
\text{Ric}(X,X)\geq (n-1)/r^{2},
\end{align*}
for every unit vector $X\in \mathfrak{X}(\mathcal{M})$, where $r$
is a positive constant. Then every geodesic $\gamma$ on $\mathcal{M}$ of
length $L[\gamma]> \pi r$ contains conjugate points; and hence is not
minimal.
\end{theorem}
Notice that the Ricci tensor is related to the sectional curvature. In particular, given a unit-length vector $X\in T_{p}\mathcal{M}$, we obtain the Ricci curvature of $X$ by extending $X=U_{n}$ to an orthonormal basis $U_{1},\cdots,U_{n}$. Then,
\begin{align*}
\text{Ric}(X,X)=\text{Tr}(\langle  R(U_{n},U_{i})U_{n},U_{i} \rangle)
&\stackrel{\rm{(a)}}{=}\sum_{i=1}^{n-1}\langle R(U_{n},U_{i})U_{n},U_{i}\rangle\\
&=\sum_{i=1}^{n-1}K(U_{n},U_{i}),
\end{align*}
where $\rm{(a)}$ is due to the fact that $\langle R(U_{n},U_{n})U_{n},U_{n}\rangle=0$. Now, the fact that $K\geq \kappa>0$ implies that $\text{Ric}\geq (n-1)\kappa$. Therefore, by Theorem \ref{Thm:Meyer}, any minimal geodesics $\gamma_{0},\gamma_{1}$ must satisfy $L[\gamma_{0}],L[\gamma_{1}]\leq {\pi\over \sqrt{\kappa}}$.

We now use Theorem \ref{Thm:Toponogov}, where here the model space $\mathcal{M}_{\kappa}$ is a sphere of radius $1/\sqrt{\kappa}$.  
We consider a geodesic hinge where $x_{t}$ is the joint of this hinge. Further, the end points of the geodesics $\gamma$ and $\gamma_{0}$ are respectively $x_{\ast}$ and $x_{t+1}$. Since our derivations for the sphere is based on the assumption that $L[\tilde{\gamma}]=d(x_{t},x_{\ast})<\dfrac{\pi}{\sqrt{\kappa}}$ and $L[\gamma]=L[\tilde{\gamma}]$, the condition $L[\gamma]\leq {\pi\over \sqrt{\kappa}}$ is satisfied. We now have the following chain of inequalities, 
\begin{align}
\cos(\sqrt{\kappa}L[\gamma_{1}])&\stackrel{\rm{(a)}}{\geq} \cos(\sqrt{\kappa} L[\tilde{\gamma}_{1}])\\
&\stackrel{\rm{(b)}}{=}\cos(\sqrt{\kappa}L[\tilde{\gamma}])\cos(\sqrt{\kappa}L[\tilde{\gamma}_{0}])+\sin(\sqrt{\kappa}L[\tilde{\gamma}]) \sin(\sqrt{\kappa}L[\tilde{\gamma}_{0}])\cos\vartheta\\ \label{Eq:Inequality_L}
&\stackrel{\rm{(c)}}{=}\cos(\sqrt{\kappa}L[\gamma])\cos(\sqrt{\kappa}L[\gamma_{0}])+\sin(\sqrt{\kappa}L[\gamma]) \sin(\sqrt{\kappa}L[\gamma_{0}])\cos\vartheta,
\end{align}
where $\rm{(a)}$ follows from inequality \eqref{Eq:Inequality_Top} and the fact that $\cos(x)$ is a monotone decreasing function on the interval $x\in [0,\pi]$, $\rm{(b)}$ follows from the law of cosine on the sphere \cite{meyer1989toponogov}, and $\rm{(c)}$ follows from the fact that $L[\tilde{\gamma}]=L[\gamma]$, $L[\tilde{\gamma}_{0}]=L[\gamma_{0}]$ due to Theorem \ref{Thm:Toponogov}. Hence, by putting $L[\gamma_{1}]=d(x_{t+1},x_{\ast})$, $L[\gamma]=d(x_{t},x_{\ast})$, and $L[\gamma_{0}]=d(x_{t+1},x_{t})$ in inequality \eqref{Eq:Inequality_Martingale}, we have that following inequality for any elliptic manifolds with $K\geq \kappa>0$,
\begin{align}
\nonumber
\cos(\sqrt{\kappa}d(x_{t+1},x_{\ast}))\geq &\cos(\sqrt{\kappa}d(x_{t},x_{\ast}))\cos(\sqrt{\kappa}d(x_{t},x_{t+1}))\\ \label{Eq:we_obtain_2}
&+\sin(\sqrt{\kappa}d(x_{t},x_{\ast})) \sin(\sqrt{\kappa}d(x_{t},x_{t+1}))\cos\vartheta.
\end{align}

Notice that here $d(x_{t},x_{t+1})$ is well-defined since we assumed that the step size satisfies $\eta_{t}\leq i(\mathcal{M})/\|X_{t}\|$.

\section{Proof Lemma \eqref{Eq:Proof_of_In}}
\label{Eq:Proof_of_Inequality_In}
Recall the update rule for $\lambda_{t+1}$ from Algorithm \ref{alg:1}. We have
\begin{align*}
\|\lambda_{t+1}\|^{2}&= \|\Pi_{\real_{+}^{m}}\left(\lambda_{t}+\eta_{t}\grad_{\lambda}\mathfrak{L}(x_{t},\lambda_{t}) \right)\|^{2}\\
&\leq \|\lambda_{t}+\eta_{t}\grad_{\lambda}\mathfrak{L}(x_{t},\lambda_{t})\|^{2}\\
&=\|\lambda_{t}\|^{2}+\eta_{t}^{2}\|\grad_{\lambda}\mathfrak{L}(x_{t},\lambda_{t})\|^{2}+2\eta_{t}\langle \grad_{\lambda}\mathfrak{L}(x_{t},\lambda_{t}),\lambda_{t} \rangle\\
&\leq \|\lambda_{t}\|^{2}+\eta^{2}_{t}\|\grad_{\lambda}\mathfrak{L}(x_{t},\lambda_{t})\|^{2},
\end{align*}
where the last inequality is due to the fact that
\begin{align*}
\langle \grad_{\lambda}\mathfrak{L}(x_{t},\lambda_{t}),\lambda_{t} \rangle= \langle h(x_{t})-\alpha\lambda_{t},\lambda_{t}\rangle=\langle h(x_{t}),\lambda_{t} \rangle-\alpha \|\lambda_{t}\|^{2}<0,
\end{align*}
since $\langle h(x_{t}),\lambda_{t} \rangle<0$.

\section{Proof of Lemma \ref{Lemma:2}}
\label{Proof of Lemma:2}
Let $p\in {\rm I\!N}$, and $\vertiii{e_{t}}_{\psi_{\nu}}=\beta$. Recall that 
$Z_{t}\topdoteq\|e_{t}\|^{2}-\expect[\|e_{t}\|^{2}|\mathfrak{F}_{t}]$, where $\mathfrak{F}_{t}$ is the $\sigma$-field 
of all the random variables $\xi_{0},\cdots,\xi_{t-1}$.

The $p$-th moment of $Z_{t}$ is upper bounded as follows 
\begin{align}
\nonumber
\expect\Big[|Z_{t}|^{p}|\mathfrak{F}_{t}\Big]&=\expect\Big[\big|\|e_{t}\|^{2}-\expect[\|e_{t}\|^{2}|\mathfrak{F}_{t}]\big|^{p}|\mathfrak{F}_{t}\Big]\\ \nonumber &\stackrel{\rm{(a)}}{\leq} 2^{p-1}\expect\left[\|e_{t}\|^{2p}|\mathfrak{F}_{t}\right]+2^{p-1}\expect\left[\big|\expect[\|e_{t}\|^{2}|\mathfrak{F}_{t}]\big|^{p}|\mathfrak{F}_{t}\right]\\ \nonumber
&\stackrel{\rm{(b)}}{\leq}  2^{p-1}\expect\left[\|e_{t}\|^{2p}|\mathfrak{F}_{t}\right]+2^{p-1}\expect\left[\expect[\|e_{t}\|^{2p}|\mathfrak{F}_{t}]\big|\mathfrak{F}_{t}\right]\\ \label{Eq:Ineq:1}
&\stackrel{\rm{(c)}}{=} 2^{p}\expect\left[\|e_{t}\|^{2p}|\mathfrak{F}_{t}\right],
\end{align}
where $\rm{(a)}$ follows from the inequality $|x-y|^{p}\leq 2^{p-1}(|x|^{p}+|y|^{p})$ for $p\in {\rm I\!N}$, $\rm{(b)}$ follows from Jensen's inequality, and $\rm{(c)}$ follows from the law of total expectation. Now, we compute an upper bound for $\expect\left[\|e_{t}\|^{2p}|\mathfrak{F}_{t}\right]$ as follows
\begin{align}
\nonumber
\expect\left[\|e_{t}\|^{2p}|\mathfrak{F}_{t}\right]&=\int_{0}^{\infty}2pu^{2p-1}\prob(\|e_{t}\|\geq u)du\\ \nonumber
&=\int_{0}^{\infty}2pu^{2p-1}\prob\left(e^{\|e_{t}\|^{\nu}/\beta^{\nu}}\geq e^{u^{\nu}/\beta^{\nu}}\right)du\\ \nonumber
&\stackrel{\rm{(d)}}{\leq} \int_{0}^{\infty}2pu^{2p-1}\expect\left[e^{\|e_{t}\|^{\nu}/\beta^{\nu}}\right]\cdot e^{-u^{\nu}/\beta^{\nu}} du\\ \nonumber
&\stackrel{\rm{(e)}}\leq 4p\int_{0}^{\infty}u^{2p-1}e^{-u^{\nu}/\beta^{\nu}} du\\ \label{Eq:Ineq:2}
&\stackrel{\rm{(f)}}= \dfrac{4p\beta^{2p}}{\nu}\Gamma\left({2p\over \nu}\right),
\end{align}
for $p=1,2,\cdots$, where $\rm{(d)}$ follows from Markov's inequality, $\rm{(e)}$ follows from the fact that $\vertiii{e}_{\psi_{\nu}}=\beta$ and hence $\expect\left[e^{\|e_{t}\|^{\nu}/\beta^{\nu}}\right]\leq 2$ by the monotone convergence theorem, $\rm{(f)}$ follows by first applying a change of variable $t=u^{\nu}/\beta^{\nu}$ and then using the definition of the Gamma function $\Gamma(x)=\int_{0}^{\infty}t^{x-1}\exp(-t)dt$.

By putting together the inequalities \eqref{Eq:Ineq:1} and \eqref{Eq:Ineq:2}, we obtain
\begin{align*}
\expect\left[|Z_{t}|^{p}\big| \mathfrak{F}_{t}\right]=\expect\left[\big|\|e_{t}\|^{2}-\expect[\|e_{t}\|^{2}|\mathfrak{F}_{t}]\big|^{p}| \mathfrak{F}_{t}\right]\leq 2^{p}\dfrac{4p\beta^{2p}}{\nu}\Gamma\left({2p\over \nu}\right), \quad p\in {{\rm I\!N}}.
\end{align*}
Given $\gamma>2\beta^{2}$, we now have
\begin{align}
\nonumber
\expect\left[\exp(|Z_{t}|^{\nu\over 2}/\gamma^{\nu\over 2})|\mathfrak{F}_{t} \right]&=\sum_{p=0}^{\infty}\dfrac{1}{p!\gamma^{(\nu p)/2}}\expect\left[|Z_{t}|^{(\nu p)/2}|\mathfrak{F}_{t}\right]\\ \nonumber &\leq 1+\sum_{p=1}^{\infty}\dfrac{1}{p!\gamma^{(\nu p)/2}}2^{{\nu p\over 2}+1}\beta^{\nu p}p\Gamma(p)\\  \nonumber
&\stackrel{\rm{(g)}}{=} -1+2\sum_{p=0}^{\infty} \left(\Big(\dfrac{\sqrt{2}\beta}{\sqrt{\gamma}}\Big)^{\nu}\right)^{p}\\ \label{Eq:Inequality-2}
&=\dfrac{2}{1-(\sqrt{2}\beta/\sqrt{\gamma})^{\nu}}-1,
\end{align}
where $\rm{(g)}$ follows from the identity $p\Gamma(p)= p!$. Therefore, from the upper bound \eqref{Eq:Inequality-2}, we conclude that it suffices that $\gamma\geq 2\beta^{2}3^{{2\over \nu}} $ to ensure 
\begin{align*}
\expect\left[\exp(|Z_{t}|^{\nu/2}/\gamma^{\nu/2}) |\mathfrak{F}_{t}\right]\leq 2.
\end{align*}

\section{Lie Groups and Haar Measure}
\label{sec:Lie Groups}
Here, we briefly review the important aspects of Lie groups that we require for our numerical studies. A more comprehensive treatise can be found in \cite[Chapter 5]{lee2009manifolds}.

\begin{definition}
A Lie group $(G,\ast)$ is a smooth Riemannian manifold and an abstract group with the group operation $\ast$ and the identity element $e$, such that $G\times G\rightarrow G:(g_{1},g_{2})\mapsto g_{1}\ast g_{2}$ and inverse $G\rightarrow G:g\mapsto g^{-1}$ are $C^{\infty}$ maps.
\end{definition}

\begin{definition}
A Lie algebra is a vector space $\mathfrak{g}$ over the field $\mathbb{F}$, equipped with with a skew-symmetric bilinear form $[\cdot,\cdot]:\mathfrak{g}\times \mathfrak{g}\rightarrow \mathfrak{g}$ that satisfies the Jacobi identity,
\begin{align*}
[x,[y,z]]+[z,[x,y]]+[y,[z,x]]=0,
\end{align*}
for all $x,y,z\in \mathfrak{g}$.
\end{definition}

The connection between the Lie algebras and Lie groups is established through left invariant vector fields. In particular, we define the left translation as $L_{g}:G\rightarrow G: L_{g}(h)=g\ast h$. Similarly, for the right translation $R_{g}(h)=h\ast g$. A vector field $X\in \mathfrak{X}(G)$ is left invariant iff $d_{h}L_{g}X(h)=X(g\ast h)$ for all $g,h\in G$, \textit{i.e.}, the following diagram commutes
\begin{align*}
\begin{CD}
TG @>dL_{g}>> TG\\
@AAXA @AAXA\\
G @>L_{g}>> G
\end{CD}
\end{align*}
Here, $dL_{g}$ is the differential of $L_{g}$.

Let $\mathfrak{g}$ denotes the set of all left invariant vector fields. Given a vector $v\in T_{e}G$, define the map $L^{v}(g)=d_{e}L_{g}v$ which provides a smooth left invariant vector field. Then $\mathfrak{g}\cong T_{e}G$ by the linear isomorphism $v\mapsto L^{v}$.
\begin{definition}
\label{Def:Lie_group}
The (group) exponential map on Lie groups is defined as 
\begin{align}
\text{Exp}: \mathfrak{g}\rightarrow G, \quad X\mapsto \Phi_{X}(1),
\end{align}
where $\Phi_{X}(t)$ is the integral curve of $\dot{g}(t)=X(g(t))$ with the boundary condition $g(0)=e$. 
\end{definition}
We remark that in general the group exponential map (Def. \ref{Def:Lie_group}) is different from the geodesic exponential map (Def. \ref{Def:exponential_map}). However, these two notions coincide, \textit{i.e.}, $\text{Exp}(X)=\Exp_{e}(X)$, under the condition that the Riemannian metric is both left and right invariants (bi-invariant) \cite{pennec2013bi},
\begin{align*}
\langle X,Y \rangle_{g}&=\langle d_{g}L_{h}X,d_{g}L_{h}Y\rangle, \\ 
\langle X,Y \rangle_{g}&=\langle d_{g}R_{h}X,d_{g}R_{h}Y\rangle,\quad \forall g,h\in G, \ \ X,Y\in T_{g}G. 
\end{align*}

We now define the notion of a \textit{Haar measure} for the compact Lie group $G$. Let $\mathcal{B}(G)$ be the Borel $\sigma$-algebra of $G$, \textit{i.e.,} the smallest $\sigma$-algebra that contains all the open sets of $G$. A regular Borel measure on $G$ is \textit{left} Haar measure (left invariant) if $\mu(gS)=\mu(S)$ for all $S\in \mathcal{B}(G)$ and $g\in G$. Similarly, $\mu$ is a \textit{right} Haar measure (right invariant) if $\mu(Sg)=\mu(S)$ for all $S\in \mathcal{B}(G)$ and $g\in G$. Now, we have the following theorem (cf. \cite[Thm. 1.2.1]{applebaum2014probability})

\begin{theorem}
Let $G$ be locally compact and Hausdorff. Left Haar measure and right Haar measure exist and are unique up to a positive multiplicative constant. Furthermore, such measure assigns finite mass to compact sets in $G$.
\end{theorem}

Among the example of compact Hausdorff Lie groups are ${\rm O}(n)$, ${\rm SO}(n)$, ${\rm U}(n)$, and ${\rm SU}(n)$.

\section{Review of Orlicz Spaces and Orlicz Norm}
\label{sec:Heavy_Tail}
In this appendix we present basic definitions and results about Orlicz spaces and the corresponding Orlicz norms. For a more comprehensive treatment on this subject see \cite{chen1996geometry}.

Let $(\mathcal{X};\Sigma;\mu)$ be a measure space and $\psi$ be a Young-Orlicz modulus (cf. Definition \ref{Def:Orlicz}). Let $L^{\psi}(\mathcal{X})$ denote the space of all real-valued measurable functions $f$ onto $\mathcal{X}$ such that
\begin{align*}
\|f\|_{\psi}\topdoteq\inf\{\beta:\expect_{\mu}\psi(|f|/\beta)\leq 1\}<\infty. 
\end{align*}
Define $\mathcal{L}^{\psi}(\mathcal{X})\topdoteq L^{\psi}(\mathcal{X})/\sim$ through the identification that $f\sim g$ if $f=g$ almost everywhere with respect to the measure $\mu$. Note that for the choice of Young-Orlicz modulus $\phi_{p}=x^{p}$, we recover the familiar $\mathcal{L}^{p}(\mathcal{X})$ space; the space of functions for which the $p$-th power of the absolute value is Lebesgue integrable, \textit{i.e.},
\begin{align*}
\|f\|_{p}=\left(\int_{\mathcal{X}}|f|^{p}\text{d}\mu\right)^{1/p}<\infty.
\end{align*}
Further, note that for any measure space $(\mathcal{X},\Sigma,\mu)$ and any Young-Orlicz modulus $\psi$, $\mathcal{L}^{\psi}(\mathcal{X})$ is a complete normed vector space (Banach space) \cite{dudley1999uniform}.

Now, we state a concentration inequality for the sum of random variables that we use in our derivations:

\begin{lemma}\textsc{(Talagrand \cite[Corollaries 2.8-2.9]{talagrand1994supremum})}
\label{Lemma:1}
Let $Y_{0},Y_{1},\cdots,Y_{T-1}$ be independent random variables centered at zero such that $\|Y_{i}\|_{\psi_{\nu}}\leq \zeta$.  There exists an absolute constant $c>0$ depending on $\nu$ only such that for every sequence $\eta\ \dot{=}\ (\eta_{0},\cdots,\eta_{T-1})$ and $\varepsilon>0$ we have\footnote{For the case of $1\leq \nu\leq 2$, the constant factor of $2$ on the right hand side does not appear in Corollary 2.9 of \cite{talagrand1994supremum}. However, we include this factor to have a symmetric bound in both cases of $1\leq \nu\leq 2$ and $\nu\geq 2$.}
\begin{itemize}
\item $1\leq \nu\leq 2$:
\begin{align*}
\prob\left[\left|\sum_{t=0}^{T-1}\eta_{t}Y_{t}\right|\geq \varepsilon \zeta \right]\leq 2\exp\left(-c\cdot \min\left(\dfrac{\varepsilon^{2}}{\|\eta\|_{2}^{2}},\dfrac{\varepsilon^{\nu}}{\|\eta\|_{\nu_{*}}^{\nu}}\right)\right),
\end{align*}
\item $\nu\geq 2$
\begin{align*}
\prob\left[\left|\sum_{t=0}^{T-1}\eta_{t}Y_{t}\right|\geq \varepsilon \zeta \right]\leq 2\exp\left(-c\cdot \max\left(\dfrac{\varepsilon^{2}}{\|\eta\|_{2}^{2}},\dfrac{\varepsilon^{\nu}}{\|\eta\|_{\nu_{*}}^{\nu}}\right)\right),
\end{align*}
\end{itemize}
where $\nu^{-1}+\nu_{*}^{-1}=1$.
\end{lemma}
\begin{table}
\centering{
\caption{List of symbols and their description}
\label{Table:1}
\begin{tabular}{SS} \toprule[0.25ex] \midrule
    {\text{Symbol}} & {\text{Description}} \\ \midrule
    $\mathcal{M}$   &\text{Riemannian Manifold}\\
    $T_{p}\mathcal{M}$  &\text{Tangent vector at $p\in \mathcal{M}$} \\
    $T^{\ast}_{p}\mathcal{M}$  &\text{Cotangent vector at $p\in \mathcal{M}$} \\
    $T\mathcal{M}$  &\text{Tangent bundle} \\
    $\mathfrak{X}\mathcal{(M)}$ &\text{Space of smooth vector fields on $\mathcal{M}$}\\    
    $G$  &\text{Lie group}  \\ 
    $\mathfrak{g}$  &\text{Lie algebra}   \\
    $\gamma$  &\text{Geodesic}  \\
    $\Exp_{p}$ &\text{Geodesic exponential map at $p\in \mathcal{M}$}\\
    $\Gamma_{ij}^{k}$  &\text{Christoffel symbols}\\
    $g_{p}\mathcal{(\cdot,\cdot)}$  &\text{Metric tensor at $p\in \mathcal{M}$}  \\
    $\nabla$  &\text{Levi-Civita connection}  \\ 
    $K\mathcal{(\cdot,\cdot)}$  &\text{Sectional curvature tensor}\\
    $\text{Ric} \mathcal{(\cdot,\cdot)}$ &\text{Ricci curvature}\\ \bottomrule[0.25ex]
\end{tabular}
}
\end{table}

\section*{Acknowledgments}
The research of MBK is supported by IBM PhD Fellowship. 

\bibliographystyle{siamplain}

\bibliography{references}

\end{document}